\documentclass[10pt]{article}
\usepackage{eurosym}
\usepackage{amsfonts}
\usepackage{amsmath,amssymb}
\usepackage[a4paper]{geometry}

\newtheorem{Lemma}{Lemma}[section]
\newtheorem{theorem}[Lemma]{Theorem}
\newtheorem{proposition}[Lemma]{Proposition}
\newtheorem{lemma}[Lemma]{Lemma}

\newtheorem{remark}[Lemma]{Remark}
\newtheorem{definition}[Lemma]{Definition}

\begin{document}
\begin{center}
{\Large \textbf{Optimal control for two-dimensional stochastic second grade
fluids}} \vspace{3mm}\\[0pt]

\vspace{1 cm} \textsc{\ Nikolai Chemetov}{\footnote{%
Departamento de Matem\'atica, Faculdade de Ci\^encias da Universidade de
Lisboa. E-mail: nvchemetov{\char'100}fc.ul.pt.},} \hspace{0mm} \textsc{\
Fernanda Cipriano}{\footnote{%
Centro de Matem\'atica e Aplica\c c\~oes (CMA) FCT/UNL and Departamento de
Matem\'atica, Faculdade de Ci\^encias e Tecnologia da Universidade Nova de
Lisboa. E-mail: cipriano{\char'100}fct.unl.pt.}}


\end{center}

\date{\today }

\begin{abstract}
This article deals with a stochastic control problem for certain fluids of
non-Newtonian type. More precisely, the state equation is given by the
two-dimensional stochastic second grade fluids perturbed by a multiplicative
white noise. The control acts through an external stochastic force and we
search for a control that minimizes a cost functional. We show that the
G\^ateaux derivative of the control to state map is a stochastic process
being the unique solution of the stochastic linearized state equation. The
well-posedness of the corresponding stochastic backward adjoint equation is
also established, allowing to derive the first order optimality condition.
\vspace{2mm}\newline

\vspace{3mm}\newline
\vspace{3mm}
\end{abstract}

\vspace{ 1 cm} 
\noindent


\section{Introduction}

\setcounter{equation}{0}\label{1}

\bigskip

This paper is devoted to the study of a control problem for stochastic
incompressible fluid of second grade filling a bounded two-dimensional
domain $\mathcal{O}$. The evolution equations are given by
\begin{equation}
\left\{
\begin{array}{l}
\tfrac{\partial }{\partial t}\left( Y-\alpha \Delta Y\right) =\nu \Delta Y-%
\mathrm{curl}\left( Y-\alpha \Delta Y\right) \times Y-\nabla \pi +U+G(t,Y)\,%
\dot{W}_{t},\vspace{2mm} \\
\mathrm{div}\,Y=0\qquad \qquad \qquad \qquad \qquad \qquad \qquad \qquad
\qquad \qquad \qquad \qquad \qquad \text{in }(0,T)\times \mathcal{O},\vspace{%
2mm}%
\end{array}%
\right.  \label{equation_etat_temps}
\end{equation}%
where $Y$ denotes the velocity of the fluid, $G(t,Y)\,\dot{W}_{t}$ is a
multiplicative white noise and the control is exercised through an adapted
stochastic process $U=U(\omega ,t)$. \ Due to the adaptiveness, such
stochastic external force $U$ also acts as a feedback control in the sense
of the noise.

To study the well-posedness of the state equations ({\ref%
{equation_etat_temps}), we should impose a suitable boundary condition on
the boundary $\Gamma $ of the domain $\mathcal{O}$. Considering the
classical homogeneous Dirichlet boundary condition, the authors in \cite%
{RS12} proved the existence and uniqueness of strong stochastic solutions.
Other physically relevant boundary condition is the so-called Navier-slip
boundary condition, which reads as
\begin{equation}
Y\cdot \mathrm{n}=0,\qquad (\mathrm{n}\cdot DY)\cdot \mathrm{\tau }=0\qquad %
\mbox{on}\ \Gamma ,  \label{NS}
\end{equation}%
where $\mathrm{n}=(n_{1},n_{2})$ and \textrm{$\tau $}$=(-n_{2},n_{1})$ are
the unit normal and tangent vectors, respectively, to the boundary $\Gamma $
and $DY=\frac{1}{2}\left[ \nabla Y+\nabla Y^{\top }\right] $ is the
rate-of-strain tensor. This boundary condition allows fluid particles to
move tangentially to the boundary of the domain and for the Navier-Stokes
fluid it has the special feature of being compatible with the limit
transition on vanishing viscosity (see \cite{CC_1_13}, \cite{CC_2_13}, \cite%
{CT15}, \cite{CMR98}, \cite{K06}). Let us refer that the state equations ({%
\ref{equation_etat_temps}) supplemented with the Navier-slip boundary
condition ({\ref{NS}) has been studied in the article \cite{CC16}, where the
authors obtained the existence and uniqueness results, as well as the
stability properties of the solutions. }}} Concerning the physical
justification and relevance of the second grade fluid equations, we refer to
the articles \cite{DF74}, \cite{DR95}, \cite{NT65} and \cite{RE55}. Related
non-viscous model can be found in \cite{HMR981} and \cite{HMR98}.

The main goal of the present article is to study the control problem for the
stochastic state equation ({\ref{equation_etat_temps}) with the Navier-slip
boundary condition ({\ref{NS}). As far as we know this is the first article
where the control problem is addressed for the stochastic second grade
fluids. Let us mention the three key steps to be overtaken in solving the
control problem: one of them relies on the study of the G\^{a}teau
derivative of the control to state map, which will be given by the unique
solution of the stochastic linearized state equation, the second one
consists in the determination and well-posedness of the stochastic backward
adjoint equations in a suitable functional space, and the third one is the
relation between the solution of the linearized equation and {the adjoint
state solution}, which allows to write the necessary first order optimality
condition in terms of the adjoint state. }}

In the framework of the Newtonian fluids, we refer \cite{CG07}, \cite{DaD00}%
, \cite{MS03}, \cite{Sr00} where the authors studied control problems for
the stochastic Navier-Stokes involving different techniques. In the present
work we follow the methods introduced in \cite{B99}, \cite{L00} and \cite%
{L02} for the two-dimensional Navier-Stokes equations. When comparing with
the Navier-Stokes equations, here the main difficulty relies on the fact
that the nonlinear term contains third order derivatives that substantially
complicates the analysis. Fortunately, we take advantage of the Navier-slip
boundary condition because it allows to estimate the vorticity of the fluid
on the boundary throughout the tangential velocity which roughly speaking
reduces the order of the derivatives, and consequently allows to obtain a
priori estimates for the Galerkin approximations of the forward stochastic
linearized and backward stochastic adjoint equations in the Sobolev space $%
H^2$.

The plan of the present paper is as follows. The formulation of the problem
is stated in Section \ref{2}. In Section \ref{3}, we present some auxiliary
lemmas concerning the normas of the functional spaces, relevant properties
and estimates of the nonlinear terms. In Section \ref{4}, we collect useful
results and estimates for the solution of the stochastic state equation.
Section \ref{5} deals with the well-posedness of the stochastic linearized
equation. In Section \ref{6}, we analyse the derivative of the control to
state map. Section \ref{7} is devoted to the formulation of the stochastic
backward adjoint equations and to the study of the existence and uniqueness
of the solutions. We also improve the integrability properties for the
solution of the stochastic state equation, and deduce the duality relation
between the linearized and the adjoint stochastic processes. Finally, in
Section \ref{81} we establish the main result of the article, by proving the
existence of a solution to the control problem and establishing a first
order optimality condition. 


\section{Functional spaces and formulation of the problem}

\setcounter{equation}{0}\label{2}
We consider a stochastic second grade fluid model in a bounded domain $%
\mathcal{O}$ of $\mathbb{R}^{2}$ with a sufficiently regular boundary $%
\Gamma $
\begin{equation}
\left\{
\begin{array}{cc}
d(\sigma (Y))=(\nu \Delta Y-\mathrm{curl}(\sigma (Y))\times Y-\nabla \pi
+U)\,dt+G(t,Y)\,dW_{t}, &  \\
\multicolumn{1}{l}{\mathrm{div}\,Y=0} & \multicolumn{1}{l}{\mbox{in}\
(0,T)\times \mathcal{O},} \\
\multicolumn{1}{l}{Y\cdot \mathrm{n}=0,\qquad (\mathrm{n}\cdot DY)\cdot
\mathrm{\tau }=0} & \multicolumn{1}{l}{\mbox{on}\ (0,T)\times \Gamma ,} \\
\multicolumn{1}{l}{Y(0)=Y_{0}} & \multicolumn{1}{l}{\mbox{in}\ \mathcal{O},}%
\end{array}%
\right.  \label{equation_etat}
\end{equation}%
where $\nu >0$ is a constant viscosity of the fluid, $\alpha >0$ is a
constant material modulus, $\ \Delta $ and $\nabla $ denote the Laplacian
and the gradient, $Y=(Y_{1},Y_{2})$ is a 2D velocity field and
\begin{equation*}
\sigma (Y)=Y-\alpha \Delta Y.
\end{equation*}%
The function $\pi $ represents the pressure, $U$ is a distributed mechanical
force and the term $G(t,$ $Y)\,dW_{t}$ corresponds to the stochastic
perturbation, where $W_{t}=(W^1_{t}, \dots,W^m_{t})$ is a standard $\mathbb{R%
}^{m}$-valued Wiener process defined on a complete probability space $%
(\Omega ,\mathcal{F},P)$ endowed with a filtration $\left\{ \mathcal{F}%
_{t}\right\} _{t\in \lbrack 0,T]}$ for the Wiener process.

Let us define some normed spaces. Let $X$ be a real Banach space with norm $%
\left\Vert \cdot \right\Vert _{X}.$ We denote by $L^{p}(0,T;X)$ the space of
$X$-valued measurable $p-$integrable functions defined on $[0,T]$ for $p\geq
1$.

For $p,r\geq 1$ let $L^{p}(\Omega ,L^{r}(0,T;X))$ be the space of stochastic
processes $Y=Y(\omega ,t)$ with values in $X$ defined on $\ \Omega \times
\lbrack 0,T],$ adapted to the filtration $\left\{ \mathcal{F}_{t}\right\}
_{t\in \lbrack 0,T]}$ ,\ and endowed with the norms
\begin{equation*}
\left\Vert Y\right\Vert _{L^{p}(\Omega ,L^{r}(0,T;X))}=\left( \mathbb{E}%
\left( \int_{0}^{T}\left\Vert Y\right\Vert _{X}^{r}\,dt\right) ^{\frac{p}{r}%
}\right) ^{\frac{1}{p}}\text{ }
\end{equation*}%
and%
\begin{equation*}
\left\Vert Y\right\Vert _{L^{p}(\Omega ,L^{\infty }(0,T;X))}=\left( \mathbb{E%
}\sup_{t\in \lbrack 0,T]}\left\Vert Y\right\Vert _{X}^{p}\ \right) ^{\frac{1%
}{p}}\quad \text{if }r=\infty ,
\end{equation*}%
where $\mathbb{E}$ is the mathematical expectation with respect to the
probability measure $P.$ As usual in the notation of processes $Y=Y(\omega
,t)$ we normally omit the dependence on $\omega \in \Omega .$

In equation (\ref{equation_etat}) the vector product $\times $ for 2D
vectors $u=(u_{1},u_{2})$ and $v=(v_{1},v_{2})$ \ is calculated as $u\times
v=(u_{1},u_{2},0)\times (v_{1},v_{2},0);$ the curl of the vector $u$ is
equal to $\mathrm{curl}\,u=\tfrac{\partial u_{2}}{\partial x_{1}}-\tfrac{%
\partial u_{1}}{\partial x_{2}}$ and the vector product of $\mathrm{curl}\,u$
with the vector $v$ is understood as
\begin{equation*}
\mathrm{curl}\,u\times v=(0,0,\mathrm{curl}\,u)\times (v_{1},v_{2},0).
\end{equation*}
Given two vectors $u,v\in \mathbb{R}^{2}$, $u\cdot
v=\sum_{i=1}^{2}u_{i}v_{i} $ corresponds to the usual scalar product in $%
\mathbb{R}^{2}$.

Let us introduce the following Hilbert spaces
\begin{align}
H& =\left\{ v\in L^{2}(\mathcal{O})\mid \mathrm{div}\,v=0\ \text{ in }%
\mathcal{O}\ \mbox{ and }\ v\cdot \mathrm{n}=0\ \mbox{ on }\Gamma \right\} ,%
\vspace{2mm}  \notag \\
V& =\left\{ v\in H^{1}(\mathcal{O})\mid \mathrm{div}\,v=0\ \mbox{ in }\
\mathcal{O}\mbox{ and }\ v\cdot \mathrm{n}=0\ \text{ on }\ \Gamma \right\} ,%
\vspace{2mm}  \notag \\
W& =\left\{ v\in V\cap H^{2}(\mathcal{O})\mid (\mathrm{n}\cdot Dv)\cdot
\mathrm{\tau }=0\ \ \mbox{on}\ \Gamma \right\} ,\vspace{2mm}  \notag \\
\widetilde{W}& =W\cap H^{3}(\mathcal{O}).  \label{w}
\end{align}%
We denote by $(\cdot ,\cdot )$ the inner product in $L^{2}(\mathcal{O})$ and
by $\Vert \cdot \Vert _{2}$ the associated norm. The norm in the space $%
H^{p}(\mathcal{O})$ is denoted by $\Vert \cdot \Vert _{H^{p}}$.

Let us consider the Helmholtz projector $\mathbb{P}:L^{2}(\mathcal{O}%
)\longrightarrow H$. It is well know that $\mathbb{P}$ is a linear bounded
operator being characterized by the equality $\mathbb{P}v=\tilde{v}$, where $%
\tilde{v}$ is defined by the Helmholtz decomposition
\begin{equation*}
v=\tilde{v}+\nabla \phi ,\qquad \tilde{v}\in H\qquad \mbox{and}\qquad \phi
\in H^{1}(\mathcal{O}).
\end{equation*}%
On the functional spaces $V,$ $W$\ and $\widetilde{W}$ defined in (\ref{w}),
we introduce the inner products
\begin{eqnarray}
\left( u,v\right) _{V} &=&\left( \sigma (u),v\right) =\left( u,v\right)
+2\alpha \left( Du,Dv\right) ,  \notag \\
\left( u,v\right) _{W} &=&\left( \mathbb{P}\sigma (u),\mathbb{P}\sigma
(v)\right) +\left( u,v\right) _{V},  \notag \\
\left( u,v\right) _{\widetilde{W}} &=&\left( \mathrm{curl}\sigma (u),\mathrm{%
curl}\sigma (v)\right) +\left( u,v\right) _{V}  \label{ss}
\end{eqnarray}%
and denote by $\Vert \cdot \Vert _{V},$ $\Vert \cdot \Vert _{W}$, $\ \Vert
\cdot \Vert _{\widetilde{W}}$ \ \ the corresponding norms. Let us notice
that the functional spaces $V$, $W$ and $\widetilde{W}$ defined in (\ref{w})
being subspaces of Sobolev spaces are naturally endowed with the Sobolev
norms, however it is more convenient to use the norms induced by the inner
products (\ref{ss}) that are related with the structure of the equations.

\bigskip

The main goal of this paper is to control the solution of the stochastic
model (\ref{equation_etat}) by the stochastic distributed mechanical force $%
U $, which belongs to a suitable admissible set $\mathcal{U}_{ad}^{b}$, that
will be defined later on.

Let us consider the cost functional
\begin{equation}
\displaystyle J(U,Y_{U})=\mathbb{E}\int_{0}^{T}L\left( t,U,Y_{U}\right) \,dt+%
\mathbb{E}\,h(Y_{U}(T)),\qquad U\in \mathcal{U}_{ad}^{b},  \label{cost}
\end{equation}%
which will be specified in Section \ref{81}. \ We aim to control the
solution $Y$ of the stochastic second grade fluid equations $(\ref%
{equation_etat}),$ \ minimizing the cost functional (\ref{cost}) for an
appropriate force $U\in \mathcal{U}_{ad}^{b}$. More precisely, our purpose
is to solve the following problem
\begin{equation*}
(\mathcal{P})\left\{
\begin{array}{l}
\underset{U}{\mbox{minimize}}\{J(U,Y):~U\in \mathcal{U}_{ad}^{b}\quad \text{%
and}\quad \\
\qquad \qquad \qquad \qquad Y\mbox{  is  the solution of  system }\ (\ref%
{equation_etat})\mbox{  for the minimizing   }U\in \mathcal{U}_{ad}^{b}\}.%
\end{array}%
\vspace{3mm}\quad \right.
\end{equation*}

\begin{remark}
\bigskip A standard example of the cost functional (\ref{cost})\ is
presented as
\begin{equation*}
\displaystyle J(U,Y_{U})=\;\mathbb{E}\int_{0}^{T}\int_{\mathcal{O}%
}\left\vert Y_{U}-Y_{d}\right\vert ^{2}\,dxdt+\lambda \,\mathbb{E}%
\int_{0}^{T}\int_{\mathcal{O}}\left\vert U\right\vert ^{2}\,dxdt+\mathbb{E}%
\,h(Y_{U}(T))
\end{equation*}%
with a desired target field $Y_{d}\in L^{2}(\Omega \times (0,T)\times
\mathcal{O})$ and some given $\lambda \geq 0.$ Hence the Lagrangian $L$ is
defined by
\begin{equation*}
L\left( t,U,Y_{U}\right) =\left\Vert Y_{U}-Y_{d}\right\Vert _{2}^{2}+\tfrac{%
\lambda }{2}\left\Vert U\right\Vert _{2}^{2}.
\end{equation*}
\end{remark}

\bigskip

\section{ Preliminaries}

\label{3} \setcounter{equation}{0}

\bigskip

In this section, we present some useful lemmas concerning the properties of
the operators that will be involved in the study of the stochastic forward
and backward differential equations that will be analyzed in the next
sections. They are straightforward adaptations of analogous results in the
literature, however, for the convenience of the reader, we present the
sketch of some proofs.

The first lemma is related with the boundary condition yielding a
characterization of the Navier slip boundary condition (\ref{NS}) for the
velocity field in terms of its vorticity (see Lemma 4.1 and Corollary 4.2 in
\cite{K06}). Denoting by $k$ the curvature of $\Gamma $ and parameterizing $%
\Gamma $ by arc length $s$, we have
\begin{equation*}
\tfrac{\partial \mathrm{n}}{\partial \mathrm{\tau }}=\tfrac{d\mathrm{n}}{ds}%
=k\mathrm{\tau }.
\end{equation*}
Let us recall the following relation between the anti-symmetric tensor $%
Av=\nabla v-(\nabla v)^\top$ and the vorticity operator:
\begin{equation*}
Av=\mathrm{curl}\,v{\small \left[
\begin{array}{lr}
0 & -1 \\
1 & 0%
\end{array}%
\right] .}
\end{equation*}

\begin{lemma}
Let $v\in H^{2}(\mathcal{O})$ be a vector field fulfilling the Navier slip
boundary condition (\ref{NS}). Then
\begin{equation}
\mathrm{curl}\,v=g(v)\quad \text{ with}\quad g(v)=2k\,v\cdot \mathrm{\tau }%
\quad \text{ on}\ \Gamma .  \label{curl_trace}
\end{equation}
\end{lemma}

\textbf{Proof.} The symmetry of $Dv$ and the anti-symmetry of $Av$ imply
that
\begin{equation*}
\left( Dv\right) \mathrm{\tau }\cdot \mathrm{n}=\left( Dv\right) \mathrm{n}%
\cdot \mathrm{\tau }\quad \text{and}\quad \left( Av\right) \mathrm{\tau }%
\cdot \mathrm{n}=-\left( Av\right) \mathrm{n}\cdot \mathrm{\tau }.
\end{equation*}%
It follows that $~(\nabla v)\mathrm{\tau }\cdot \mathrm{n}=\left( Dv\right)
\mathrm{n}\cdot \mathrm{\tau }-\tfrac{1}{2}\left( Av\right) \mathrm{n}\cdot
\mathrm{\tau },$ $~$which is equivalent to
\begin{equation}
\mathrm{curl}\,v=-2(\nabla v)\mathrm{\tau }\cdot \mathrm{n}+2(Dv)\mathrm{n}%
\cdot \mathrm{\tau }.  \label{navier}
\end{equation}%
Taking the derivative of the expression $v\cdot \mathrm{n}=0$ in the
direction of the tangent vector $\mathrm{\tau }$, we deduce
\begin{equation}
(\nabla v)\mathrm{\tau }\cdot \mathrm{n}=-k\,v\cdot \mathrm{\tau }.
\label{flux}
\end{equation}%
The conclusion is then a consequence of (\ref{navier}) and (\ref{flux}). $%
\hfill \hfill \blacksquare \hfill $

\bigskip

The next result follows from Lemma 5 in \cite{BI06}, and it allows to deduce
uniform estimates for the Galerkin approximations of the linearized and
adjoint equations in the Sobolev space $H^{2}.$

\begin{lemma}
\label{sigma_psigma} For each $v\in W$, we have
\begin{equation}
\left\Vert \sigma (v)-\mathbb{P}\sigma (v)\right\Vert _{2}\leq C\left\Vert
v\right\Vert _{H^{1}},  \label{sigma_psigma1}
\end{equation}%
\begin{equation}
\left\Vert \sigma (v)-\mathbb{P}\sigma (v)\right\Vert _{H^{1}}\leq
C\left\Vert v\right\Vert _{H^{2}}.  \label{sigma_psigma2}
\end{equation}
\end{lemma}

\bigskip

In the above lemma and throughout the article, we denote by $C$ a generic
\textit{positive} constant that can take different values, and may depend
only on the domain $\mathcal{O},$ \ the regularity of the boundary $\Gamma $
and the physical constant $\alpha .$

\bigskip

The following result is a direct consequence of the Korn inequality.

\begin{lemma}
\label{non_lin2} There exist some positive constants $C_{\ast },$ \ $K_{\ast
},$\ such that%
\begin{equation}
\left\Vert u\right\Vert _{H^{1}}^{2}\leq C_{\ast }\left\Vert u\right\Vert
_{V}^{2},\qquad \forall u\in V,  \label{korn}
\end{equation}%
\begin{equation}
\left\Vert u\right\Vert _{H^{1}}^{2}\leq K_{\ast }\left\Vert Du\right\Vert
_{2}^{2},\qquad \forall u\in V,  \label{kornn}
\end{equation}
\end{lemma}

\textbf{Proof.} \ The first inequality is a consequence of the Korn
inequality%
\begin{equation*}
\left\Vert u\right\Vert _{H^{1}}\leq C\left( \left\Vert Du\right\Vert
_{2}+\left\Vert u\right\Vert _{2}\right) ,\qquad \forall u\in H^{1}(\mathcal{%
O}),
\end{equation*}%
where $C$ is a positive constant just depending on $\mathcal{O}$. \ Taking
into account the results in \cite{DV02}, for non axisymmetric domains, we
also have (\ref{kornn}).$\hfill \hfill \blacksquare $

\bigskip

Now, we introduce the modified Stokes problem with the Navier-slip boundary
condition
\begin{equation}
\left\{
\begin{array}{ll}
h-\alpha \Delta h+\nabla p=f,\qquad \mathrm{div}\,h=0 & \quad \mbox{in}\
\mathcal{O},\vspace{2mm} \\
h\cdot \mathrm{n}=0,\qquad (\mathrm{n}\cdot Dh)\cdot \mathrm{\tau }=0 &
\quad \mbox{on}\ \Gamma \vspace{2mm}%
\end{array}%
\right.  \label{GS_NS}
\end{equation}%
that according to \cite{S73}, for any $f\in H^{m}(\mathcal{O})$, $m\in
\mathbb{N}$, has a unique solution $(h,p)$ which belongs to $H^{m+2}(%
\mathcal{O})\times H^{m+1}(\mathcal{O})$ and verifies
\begin{equation}
\left\Vert h\right\Vert _{H^{m+2}}\leq C\left\Vert f-h\right\Vert _{H^{m}}.
\label{h}
\end{equation}%
In what follows, we also denote the solution $h$ as
\begin{equation*}
h=\left( I-\alpha \Delta \right) ^{-1}f.
\end{equation*}


Next we state a lemma concerning the equivalence of the norms, see
Propositions 3 in \cite{BR03} and Lemma 2.1 in \cite{CG97} for similar
results.

\begin{lemma}
\label{non_lin} There exist some positive constants $C,$ $C_{\ast \ast },$\
such that%
\begin{equation}
\left\Vert u\right\Vert _{H^{2}}^{2}\leq C\left\Vert u\right\Vert
_{W}^{2},\qquad \forall u\in W,  \label{www}
\end{equation}%
\begin{equation}
\left\Vert u\right\Vert _{H^{3}}^{2}\leq C\left\Vert u\right\Vert _{%
\widetilde{W}}^{2},\qquad \forall u\in \widetilde{W},  \label{y2_sigma}
\end{equation}%
\begin{equation}
C_{\ast \ast }\left\Vert u\right\Vert _{\widetilde{W}}^{2}\leq \left\Vert
u\right\Vert _{\ast \ast }^{2},\qquad \forall u\in \widetilde{W},
\label{CCC}
\end{equation}%
where $~\left\Vert u\right\Vert _{\ast \ast }^{2}=2\alpha \left\Vert
Du\right\Vert _{2}^{2}\,+\left\Vert \mathrm{curl}\,\sigma \left( u\right)
\right\Vert _{2}^{2}.$
\end{lemma}


\textbf{Proof.} 1) Multiplying (\ref{GS_NS})$_{1}$ by $h,$ \ we derive
\begin{equation*}
\left\Vert h\right\Vert _{2}^{2}+2\alpha \left\Vert Dh\right\Vert
_{2}^{2}=(f,h)\leq \left\Vert f\right\Vert _{2}\left\Vert h\right\Vert _{2},
\end{equation*}%
that is
\begin{equation}
\left\Vert h\right\Vert _{2}\leq \left\Vert f\right\Vert _{2},\qquad
\left\Vert Dh\right\Vert _{2}\leq C\left\Vert f\right\Vert _{2}.  \label{upo}
\end{equation}%
Due to (\ref{h}) we have
\begin{equation*}
\left\Vert h\right\Vert _{H^{2}}\leq C\left\Vert f-h\right\Vert _{2}\leq
C(\left\Vert f\right\Vert _{2}+\left\Vert h\right\Vert _{2})\leq C\left\Vert
f\right\Vert _{2}.
\end{equation*}%
Hence taking $f=\mathbb{P}\sigma (u)$ in (\ref{GS_NS})$_{1}$ and using (\ref%
{upo}), we get
\begin{equation*}
\left\Vert h\right\Vert _{H^{2}}\leq C\left\Vert \mathbb{P}\sigma
(u)\right\Vert _{2}.
\end{equation*}%
By another hand if we use estimates (\ref{sigma_psigma1}), (\ref{korn}) and
apply result (\ref{h}) for $~f=\mathbb{P}\sigma (u)-\sigma (u),$ $~$we
obtain
\begin{equation*}
\left\Vert h-u\right\Vert _{H^{2}}\leq C\left\Vert \mathbb{P}\sigma
(u)-\sigma (u)\right\Vert _{2}\leq C\left\Vert u\right\Vert _{H^{1}}\leq
C\left\Vert u\right\Vert _{V}.
\end{equation*}%
Combining the above estimates for $h~$ and $~h-u,$ we deduce (\ref{www}).

2)\ Since $~\mathrm{curl}\,\sigma (v)\in L^{2}(\mathcal{O})~$ and $~\nabla
\cdot \left( \mathrm{curl}\,\sigma (v)\right) =0$,$~$ there exists a unique
vector-potential $\psi \in H^{1}(\mathcal{O})$ such that
\begin{equation*}
\left\{
\begin{array}{ll}
\mathrm{{curl}\,}\psi =\tilde{y},\qquad \nabla \cdot \psi =0 & \quad %
\mbox{in}\ \mathcal{O},\vspace{2mm} \\
\mathrm{{curl}\,}\psi =0 & \quad \mbox{on}\ \Gamma%
\end{array}%
\right.
\end{equation*}%
and
\begin{equation}
\left\Vert \psi \right\Vert _{H^{1}}\leq C\left\Vert \mathrm{curl}\,\sigma
(v)\right\Vert _{2}.  \label{sigma_phi}
\end{equation}%
It follows that ~~$\mathrm{curl}\left( v-\alpha \Delta v-\psi \right) =0$~
and there exists $~\pi \in L^{2}(\mathcal{O}),~$ such that \ ~~$v-\alpha
\Delta v-\psi +\nabla \pi =0.$ \ \ \ Hence $v$ is the solution of the Stokes
system
\begin{equation*}
\left\{
\begin{array}{ll}
v-\alpha \Delta v+\nabla \pi =\psi ,\qquad \mathrm{div}\,v=0, & \quad %
\mbox{in}\ \mathcal{O},\vspace{2mm} \\
v\cdot \mathrm{n}=0,\qquad \left( \mathrm{n}\cdot Dv\right) \cdot \mathrm{%
\tau }=0 & \quad \mbox{on}\ \Gamma .%
\end{array}%
\right.
\end{equation*}

As a consequence of (\ref{h}), we have
\begin{equation}
\left\Vert v\right\Vert _{H^{3}}\leq C\left\Vert \psi -v\right\Vert
_{H^{1}}\leq C\left( \left\Vert \psi \right\Vert _{H^{1}}+\left\Vert
v\right\Vert _{H^{1}}\right) .  \label{y_phi}
\end{equation}%
Using (\ref{sigma_phi}) we state%
\begin{equation*}
\left\Vert v\right\Vert _{H^{3}}\leq C\left( \left\Vert \mathrm{curl}%
\,\sigma (v)\right\Vert _{2}+\left\Vert v\right\Vert _{H^{1}}\right) ,
\end{equation*}%
that implies the claimed result (\ref{y2_sigma}).

\bigskip

3) To prove (\ref{CCC}) it is enough to show the existence of a positive
constant $C,$ satisfying the following estimate
\begin{equation}
\left\Vert u\right\Vert _{2}\leq C\left( \left\Vert \mathrm{curl}\,\sigma
(u)\right\Vert _{2}+\left\Vert Du\right\Vert _{2}\right) ,\qquad \forall
u\in \widetilde{W}.  \label{uk1}
\end{equation}

Let us assume that this estimate is not valid, then for any $n\in \mathbb{N}$
\ there exists a vector-function $v_{n}\in \widetilde{W},$ such that%
\begin{equation*}
\left\Vert v_{n}\right\Vert _{2}\geq n\left( \left\Vert \mathrm{curl}%
\,\sigma (v_{n})\right\Vert _{2}+\left\Vert Dv_{n}\right\Vert _{2}\right) .
\end{equation*}%
If we denote\ $z_{n}=\frac{v_{n}}{||v_{n}||_{2}},$ then the sequence $%
\left\{ z_{n}\right\} _{n=1}^{\infty }$ is uniformly bounded in $\widetilde{W%
}$
\begin{equation}
\left\Vert z_{n}\right\Vert _{2}=1,\qquad \left\Vert \mathrm{curl}\,\sigma
(z_{n})\right\Vert _{2}+\left\Vert Dz_{n}\right\Vert _{2}\leq \frac{1}{n}.
\label{oi}
\end{equation}%
Having $\widetilde{W}\subset H^{3}(\mathcal{O})$ and the compact embedding $%
H^{3}(\mathcal{O})\hookrightarrow L^{2}(\mathcal{O}),$ there exists a
subsequence $\left\{ z_{n}\right\} _{n=1}^{\infty }$ such that%
\begin{eqnarray*}
z_{n} &\rightharpoonup &z\qquad \text{ weakly in }\widetilde{W}, \\
z_{n} &\rightarrow &z\text{ \qquad\ strongly in }L^{2}(\mathcal{O}).
\end{eqnarray*}%
From (\ref{oi}) we have%
\begin{equation*}
\left\Vert z\right\Vert _{2}=1,\qquad \left\Vert \mathrm{curl}\,\sigma
(z)\right\Vert _{2}+\left\Vert Dz\right\Vert _{2}=0.
\end{equation*}%
Taking into account that $z\cdot \mathrm{n}=0\ $ on $\ \Gamma ,$ we easily
derive that $z=0,$ yielding a contradiction with $\left\Vert z\right\Vert
_{2}=1.$ Therefore estimate\ (\ref{CCC}) holds.$\hfill \hfill \blacksquare $

\bigskip

Now, we turn our attention to the nonlinear operators and introduce the
operator
\begin{equation*}
b(u,v,\phi )=\left( u\cdot \nabla v,\phi \right) ,
\end{equation*}%
which is well known in the context of Navier-Stokes equations.

The next two results deal with the properties of the nonlinear operators
that will appear in the linearized and backward adjoint equations. Due to
these results, we will be able to give a meaning to the stochastic
linearized and backward adjoint equations in the variational form, as well
as to deduce uniform estimates for the Galerkin approximations (we refer to
\cite{AC1}, \cite{AC2} and \cite{CC16} for more details).

\begin{lemma}
\label{non_lin_curl} For any $u,v\in \widetilde{W}$ and $\phi \in V$, we
have
\begin{equation}
\left( \mathrm{curl}\,\sigma (v)\times u,\phi \right) =b\left( \phi
,u,\sigma (v)\right) -b\left( u,\phi ,\sigma (v)\right) .  \label{bb1}
\end{equation}%
Moreover, if $u\in \widetilde{W}$ and $v\in W$ the following estimate holds
\begin{equation}
\left\vert \left( \mathrm{curl}\,\sigma (v)\times u,v\right) \right\vert
\leq C\left\Vert u\right\Vert _{H^{3}}\Vert v\Vert _{H^{1}}^{2}\leq
C_{1}\left\Vert u\right\Vert _{\widetilde{W}}\Vert v\Vert _{V}^{2}.
\label{rm2_lin}
\end{equation}
\end{lemma}

\textbf{Proof.} Taking into account that the vector fields $v$ and $\phi $
have zero divergence and are tangent to the boundary $\Gamma $, \ identity (%
\ref{bb1}) follows by standard computations, integrating by parts the right
hand side. In order to derive property (\ref{rm2_lin}) we set
\begin{align}
\left( \mathrm{curl}\,\sigma (v)\times u,v\right) & =b\left( v,u,\sigma
(v)\right) -b\left( u,v,\sigma (v)\right)  \notag \\
& =b(v,u,v)-b(u,v,v)-\alpha \left( b(v,u,\Delta v)-b(u,v,\Delta v)\right) .
\label{trilin_0}
\end{align}%
Integrating by parts and using the boundary conditions, we derive
\begin{align*}
b\left( u,v,\Delta v\right) & =\sum_{i,j=1}^{2}\int_{\mathcal{O}}u_{i}\tfrac{%
\partial v_{j}}{\partial x_{i}}\Delta v_{j}\ dx=\sum_{i,j,k=1}^{2}\int_{%
\mathcal{O}}u_{i}\tfrac{\partial v_{j}}{\partial x_{i}}\tfrac{\partial }{%
\partial x_{k}}\left( \tfrac{\partial v_{j}}{\partial x_{k}}-\tfrac{\partial
v_{k}}{\partial x_{j}}\right) \ dx \\
& =\sum_{i,j,k=1}^{2}\int_{\mathcal{O}}u_{i}\tfrac{\partial v_{j}}{\partial
x_{i}}\tfrac{\partial }{\partial x_{k}}(A_{jk}(v))\ dx \\
& =\sum_{i,j,k=1}^{2}\int_{\Gamma }u_{i}\tfrac{\partial v_{j}}{\partial x_{i}%
}A_{jk}(v)n_{k}\ dx-\sum_{i,j,k=1}^{2}\int_{\mathcal{O}}\tfrac{\partial }{%
\partial x_{k}}\left( u_{i}\tfrac{\partial v_{j}}{\partial x_{i}}\right)
A_{jk}(v)\ dx \\
& =\sum_{i,j=1}^{2}\int_{\Gamma }u_{i}\tfrac{\partial v_{j}}{\partial x_{i}}%
\,g(v)\tau _{j}\ dx-\sum_{i,j,k=1}^{2}\int_{\mathcal{O}}\tfrac{\partial u_{i}%
}{\partial x_{k}}\tfrac{\partial v_{j}}{\partial x_{i}}A_{jk}(v)\
dx-\sum_{i,j,k=1}^{2}\int_{\mathcal{O}}u_{i}\tfrac{\partial ^{2}v_{j}}{%
\partial x_{x}\partial x_{i}}A_{jk}(v)\ dx \\
& =I_{1}+I_{2}+I_{3}.
\end{align*}%
Again, integrating by parts, it follows that
\begin{align}
I_{1}& =\sum_{i,j=1}^{2}\int_{\Gamma }u_{i}\tfrac{\partial v_{j}}{\partial
x_{i}}\,g(v)\tau _{j}\ dx=\sum_{i,j=1}^{2}\int_{\mathcal{O}}u_{i}\tfrac{%
\partial v_{j}}{\partial x_{i}}\,g(v)\tau _{j}\,\mathrm{div}\,\mathrm{n}\
dx+\sum_{i,j,k=1}^{2}\int_{\mathcal{O}}\tfrac{\partial }{\partial x_{k}}%
\bigl(u_{i}\tfrac{\partial v_{j}}{\partial x_{i}}\,g(v)\tau _{j}\bigr)n_{k}\
dx  \notag \\
& =b\left( u,v,\mathrm{div}\,\mathrm{n}\,g(v)\mathrm{\tau }\right)
+\sum_{k=1}^{2}b\left( \tfrac{\partial u}{\partial x_{k}},v,\,n_{k}\,g(v)%
\mathrm{\tau }\right) -\sum_{k=1}^{2}b\left( u,\,n_{k}\,g(v)\mathrm{\tau },%
\tfrac{\partial v}{\partial x_{k}}\right)  \notag \\
& +\sum_{k=1}^{2}b\left( u,v,\,n_{k}\,\tfrac{\partial }{\partial x_{k}}\bigl(%
g(v)\mathrm{\tau }\bigr)\right) .  \label{GF_111m}
\end{align}

It is clear that $\displaystyle I_{2}=-\sum_{k=1}^{2}b\left( \tfrac{\partial
u}{\partial x_{k}},v,A_{\cdot \,k}(v)\right) \quad $ and $\quad I_{3}=0.$
Moreover, analogous computations can be performed in order to obtain
\begin{align*}
b\left( v,u,\Delta v\right) & =\sum_{i,j=1}^{2}\int_{\mathcal{O}}v_{i}\tfrac{%
\partial u_{j}}{\partial x_{i}}\Delta v_{j}\ dx=\sum_{i,j,k=1}^{2}\int_{%
\mathcal{O}}v_{i}\tfrac{\partial u_{j}}{\partial x_{i}}\tfrac{\partial }{%
\partial x_{k}}A_{jk}(v)\ dx \\
& =\sum_{i,j,k=1}^{2}\int_{\Gamma }v_{i}\tfrac{\partial u_{j}}{\partial x_{i}%
}A_{jk}(v)n_{k}\ dx-\sum_{i,j,k=1}^{2}\int_{\mathcal{O}}\tfrac{\partial }{%
\partial x_{k}}\bigl(v_{i}\tfrac{\partial u_{j}}{\partial x_{i}}\bigr)%
A_{jk}(v)\ dx \\
& =\sum_{i,j=1}^{2}\int_{\Gamma }v_{i}\tfrac{\partial u_{j}}{\partial x_{i}}%
g(v)\tau _{j}\ dx-\sum_{i,j,k=1}^{2}\int_{\mathcal{O}}\tfrac{\partial v_{i}}{%
\partial x_{k}}\tfrac{\partial u_{j}}{\partial x_{i}}A_{jk}(v)-%
\sum_{i,j,k=1}^{2}\int_{\mathcal{O}}v_{i}\tfrac{\partial ^{1}u_{j}}{\partial
x_{k}\partial x_{i}}A_{jk}(v)\ dx \\
& =\sum_{i,j=1}^{2}\int_{\Gamma }v_{i}\tfrac{\partial u_{j}}{\partial x_{i}}%
g(v)\tau _{j}\ dx-\sum_{k=1}^{2}b\left( \tfrac{\partial v}{\partial x_{k}}%
,u,A_{\cdot \,k}(v)\right) -\sum_{k=1}^{2}b\left( v,\tfrac{\partial u}{%
\partial x_{k}},A_{\cdot \,k}(v)\right) .
\end{align*}%
Therefore, from one hand, we have
\begin{align}
\left\vert b\left( u,v,\Delta v\right) \right\vert & \leq \left\vert b\left(
u,v,\mathrm{div}\,\mathrm{n}\,g(v)\mathrm{\tau }\right) \right\vert
+\sum_{k=1}^{2}\left\vert b\left( \tfrac{\partial u}{\partial x_{k}}%
,v,\,n_{k}\,g(v)\mathrm{\tau }\right) \right\vert +\sum_{k=1}^{2}\left\vert
b\left( u,\,n_{k}\,g(v)\mathrm{\tau },\tfrac{\partial v}{\partial x_{k}}%
\right) \right\vert  \notag \\
& +\sum_{k=1}^{2}\left\vert b\left( u,v,\,n_{k}\,\tfrac{\partial }{\partial
x_{k}}\bigl(g(v)\mathrm{\tau }\bigr)\right) \right\vert
+\sum_{k=1}^{2}\left\vert b\left( \tfrac{\partial u}{\partial x_{k}}%
,v,A_{\cdot \,k}(v)\right) \right\vert  \notag \\
& \leq C\Vert u\Vert _{H^{2}}\Vert v\Vert _{H^{1}}^{2}.  \label{lin_1}
\end{align}%
On other hand, taking into account the trace theorem $~\Vert v\Vert
_{L^{2}(\Gamma )}\leq C\Vert v\Vert _{H^{1}}$,$~$ we deduce
\begin{align}
|b\left( v,u,\Delta v\right) |& \leq \sum_{i,j=1}^{2}\int_{\Gamma }\bigl|%
v_{i}\tfrac{\partial u_{j}}{\partial x_{i}}g(v)\tau _{j}\bigr|\ dx+\
\sum_{k=1}^{2}\left( \bigl|b\left( \tfrac{\partial v}{\partial x_{k}}%
,u,A_{\cdot \,k}(v)\right) \bigr|+\bigl|b\left( v,\tfrac{\partial u}{%
\partial x_{k}},A_{\cdot \,k}(v)\right) \bigr|\right)  \notag \\
& \leq 2\Vert k\Vert _{C(\bar{\mathcal{O}})}\Vert \mathrm{n}\Vert _{C(\bar{%
\mathcal{O}})}^{2}\Vert v\Vert _{L^{2}(\Gamma )}^{2}\left\Vert \nabla
u\right\Vert _{C(\bar{\mathcal{O}})}+\left\Vert \nabla v\right\Vert
_{2}^{2}\left\Vert \nabla u\right\Vert _{\infty }+\left\Vert v\right\Vert
_{4}\left\Vert \nabla u\right\Vert _{1,4}\left\Vert \nabla v\right\Vert _{2}%
\vspace{2mm}  \notag \\
& \leq (2C\Vert k\Vert _{C(\bar{\mathcal{O}})}\Vert \mathrm{n}\Vert _{C(\bar{%
\mathcal{O}})}^{2}+1)\left\Vert \nabla v\right\Vert _{2}^{2}\left\Vert
\nabla u\right\Vert _{\infty }+\left\Vert v\right\Vert _{4}\left\Vert \nabla
u\right\Vert _{1,4}\left\Vert \nabla v\right\Vert _{2}\vspace{2mm}  \notag \\
& \leq C\left\Vert u\right\Vert _{H^{3}}\left\Vert v\right\Vert _{H^{1}}^{2}.
\label{GF_2mm}
\end{align}%
Then the claimed result is a consequence of (\ref{korn}), (\ref{y2_sigma})
and (\ref{trilin_0})-(\ref{GF_2mm}). $\hfill \hfill \blacksquare $

\bigskip

Also in the following sections we will need the following results.

\begin{lemma}
\label{prop_a_adj_2} For any $u,v\in \widetilde{W}$ and $\phi \in V$, we
have
\begin{align}
\left( \mathrm{curl}\left( \sigma \left( u\times v\right) \right) ,\phi
\right) & =b\left( \sigma (v),u,\phi \right) +b\left( u,\phi ,\sigma
(v)\right) -b\left( \sigma (u),v,\phi \right) +b\left( v,\sigma (u),\phi
\right) +b\left( u,v,\phi \right)  \notag \\
& -b\left( v,u,\phi \right) -2\alpha \sum_{i=1}^{2}\left( b(\tfrac{\partial v%
}{\partial x_{i}},\tfrac{\partial u}{\partial x_{i}},\phi )-b(\tfrac{%
\partial u}{\partial x_{i}},\tfrac{\partial v}{\partial x_{i}},\phi )\right)
,  \label{adj_for2} \\
\left( \mathrm{curl}\,\sigma \left( u\times v\right) ,\phi \right) &
=b\left( v,u,\sigma (\phi )\right) -b\left( u,v,\sigma (\phi )\right) .
\label{bb2}
\end{align}%
Moreover the following estimate hods
\begin{equation}
\left\vert \left( \mathrm{curl}\,\left( \sigma \left( u\times v\right)
\right) ,\phi \right) \right\vert \leq C\left\Vert u\right\Vert
_{H^{3}}\left\Vert v\right\Vert _{H^{2}}\left\Vert \phi \right\Vert _{H^{1}}.
\label{rm2_adj}
\end{equation}
\end{lemma}

\vspace{1mm} \textbf{Proof.} We verify equality (\ref{adj_for2}). Since $~%
\mathrm{curl}\left( \sigma \left( u\times v\right) \right) =\mathrm{curl}%
\left( u\times v-\alpha \Delta \left( u\times v\right) \right) ,~$ the
algebraic equality
\begin{equation*}
\Delta \left( u\times v\right) =u\times \Delta v-v\times \Delta
u-2\sum_{i=1}^{2}\tfrac{\partial v}{\partial x_{i}}\times \tfrac{\partial u}{%
\partial x_{i}}
\end{equation*}%
yields
\begin{equation*}
\mathrm{curl}\left( \sigma \left( u\times v\right) \right) =\mathrm{curl}%
\left( u\times \sigma (v)\right) -\mathrm{curl}\left( v\times \sigma
(u)\right) +\mathrm{curl}\left( v\times u\right) +2\alpha \sum_{i=1}^{2}%
\mathrm{curl}\left( \tfrac{\partial v}{\partial x_{i}}\times \tfrac{\partial
u}{\partial x_{i}}\right) .
\end{equation*}%
Then using the identity
\begin{equation}
\mathrm{curl}\left( \varphi \times \psi \right) =\psi \cdot \nabla \varphi
-\varphi \cdot \nabla \psi  \label{curl_prod}
\end{equation}%
for divergence free vector fields $\varphi ,$ $\psi $, \ we deduce
\begin{align*}
\mathrm{curl}\left( \sigma \left( u\times v\right) \right) & =\sigma
(v)\cdot \nabla u-u\cdot \nabla \left( \sigma (v)\right) -\sigma (u)\cdot
\nabla v+v\cdot \nabla \left( \sigma (u)\right) +u\cdot \nabla v-v\cdot
\nabla u \\
& -2\alpha \sum_{i=1}^{2}\left( \tfrac{\partial v}{\partial x_{i}}\cdot
\nabla \left( \tfrac{\partial u}{\partial x_{i}}\right) -\tfrac{\partial u}{%
\partial x_{i}}\cdot \nabla \left( \tfrac{\partial v}{\partial x_{i}}\right)
\right) ,
\end{align*}%
which implies (\ref{adj_for2}).

\bigskip

Now, we show (\ref{bb2}). We assume that $u,v\in W\cap H^{3}(\mathcal{O})$.
\ If $u,v$ are just in $H^{3}(\mathcal{O})$, the result follows from a
standard regularization procedure. We have
\begin{equation*}
\left( \mathrm{curl}\,\sigma \left( u\times v\right) ,\phi \right) =\left(
\mathrm{curl}\,\left( u\times v\right) ,\phi \right) -\alpha \left( \mathrm{%
curl}\,\Delta \left( u\times v\right) ,\phi \right) .
\end{equation*}%
Knowing that $\Delta =-\mathrm{curl}\;\mathrm{curl}$, integrating twicely by
parts and using relation (\ref{curl_prod}), we deduce
\begin{align}
\left( \mathrm{curl}\left( \Delta \left( u\times v\right) \right) ,\phi
\right) & =b(v,u,\Delta \phi )-b(u,v,\Delta \phi )  \notag \\
& +\int_{\Gamma }\left( \Delta \left( u\times v\right) \times \phi \right)
\cdot \mathrm{n}\,dx-\int_{\Gamma }\left( \mathrm{curl}\left( u\times
v\right) \times \mathrm{curl}\,\phi \right) \cdot \mathrm{n}\,dx.
\label{defBT}
\end{align}%
Algebraic calculations show that
\begin{align}
\Delta \left( u\times v\right) & =u\cdot \nabla \left( \mathrm{curl}%
\,v\right) -v\cdot \nabla \left( \mathrm{curl}\,u\right)
+2\sum_{k=1}^{2}\left( \nabla u_{k}\times \nabla v_{k}\right) ,  \notag \\
\mathrm{curl}\left( u\times v\right) & =v\cdot \nabla u-u\cdot \nabla v.
\label{12}
\end{align}%
The boundary condition (\ref{curl_trace}) on $\Gamma $ implies that
\begin{align*}
& \Delta \left( u\times v\right) =\left( 0,0,u\cdot \nabla g(v)-v\cdot
\nabla g(u)\right) +2\sum_{k=1}^{2}\left( \nabla u_{k}\times \nabla
v_{k}\right) \vspace{0mm} \\
& =2\left( 0,0,k\left( u\cdot \nabla v\right) \cdot \tau +\left( u\cdot
\nabla (k\tau )\right) \cdot v-k\left( v\cdot \nabla u\right) \cdot \tau
-\left( v\cdot \nabla (k\tau )\right) \cdot u\right) +2\sum_{k=1}^{2}\left(
\nabla u_{k}\times \nabla v_{k}\right) .
\end{align*}%
Since $v$ and $u$ are tangent to the boundary $\Gamma $, we obtain
\begin{align}
\left( \Delta \left( u\times v\right) \times \phi \right) \cdot \mathrm{n}&
=-2k\left( u\cdot \nabla v-v\cdot \nabla u\right) \cdot \tau \left( \phi
\cdot \tau \right) +2\sum_{k=1}^{2}\left( \left( \nabla u_{k}\times \nabla
v_{k}\right) \times \phi \right) \cdot \mathrm{n}  \notag \\
& -2\left\{ \left[ (u\cdot \tau )(v\cdot \tau )\left( \tau \cdot \nabla
g\right) \cdot \tau -(v\cdot \tau )(u\cdot \tau )\left( \tau \cdot \nabla
(k\tau )\right) \cdot \tau \right] \left( \phi \cdot \tau \right) \right\}
\notag \\
& =-2k\left( u\cdot \nabla v-v\cdot \nabla u\right) \cdot \tau \left( \phi
\cdot \tau \right) +2\sum_{k=1}^{2}\left( \left( \nabla u_{k}\times \nabla
v_{k}\right) \times \phi \right) \cdot \mathrm{n}.  \notag
\end{align}%
Using the same reasoning and (\ref{12}) we may verify that on $\Gamma $ the
following relation holds
\begin{align*}
\mathrm{curl}\left( u\times v\right) \times \mathrm{curl}\,\phi \cdot
\mathrm{n}& =\left( v\cdot \nabla u-u\cdot \nabla v\right) \times
(0,0,g(\phi ))\cdot \mathrm{n} \\
& =-2k\left( u\cdot \nabla v-v\cdot \nabla u\right) \left( \phi \cdot \tau
\right) \cdot \tau .
\end{align*}%
Therefore, the difference of the boundary terms is given by
\begin{equation*}
\int_{\Gamma }\left( \Delta \left( u\times v\right) \times \phi \right)
\cdot \mathrm{n}\,\ dx-\int_{\Gamma }\left( \mathrm{curl}\left( u\times
v\right) \times \mathrm{curl}\,\phi \right) \cdot \mathrm{n}\
dx=2\sum_{k=1}^{2}\int_{\Gamma }\left( \nabla u_{k}\times \nabla
v_{k}\right) \times \phi \cdot \mathrm{n}\,\ dx.
\end{equation*}%
On the other hand, the boundary condition $\left( n\cdot Dv\right) \cdot
\mathrm{\tau }=0$ and $\mathrm{div}\,v=0$ imply
\begin{equation*}
\frac{\partial v_{1}}{\partial x_{2}}+\frac{\partial v_{2}}{\partial x_{1}}=%
\frac{4n_{1}n_{2}}{n_{2}^{2}-n_{1}^{2}}\,\frac{\partial v_{1}}{\partial x_{1}%
}\quad \mbox{ on }\ \Gamma ,
\end{equation*}%
which gives
\begin{equation*}
\sum_{k=1}^{2}\left( \nabla u_{k}\times \nabla v_{k}\right) \times \phi
\cdot \mathrm{n}\big|_{\Gamma }=\left( \tfrac{\partial u_{1}}{\partial x_{1}}%
\left( \tfrac{\partial v_{1}}{\partial x_{2}}+\tfrac{\partial v_{2}}{%
\partial x_{1}}\right) -\tfrac{\partial v_{1}}{\partial x_{1}}\left( \tfrac{%
\partial u_{1}}{\partial x_{2}}+\tfrac{\partial u_{2}}{\partial x_{1}}%
\right) \right) \phi \cdot \tau =0.
\end{equation*}%
Therefore the difference of the boundary terms in (\ref{defBT}) vanishes. As
a consequence of it we derive
\begin{equation*}
\left( \mathrm{curl}\,\sigma \left( u\times v\right) ,\phi \right) =\left(
\mathrm{curl}\,\left( u\times v\right) ,\phi \right) -\alpha \left(
b(v,u,\Delta \phi )-b(u,v,\Delta \phi )\right) .
\end{equation*}%
Then using (\ref{12}), we deduce (\ref{bb2}).

Finally, estimate (\ref{rm2_adj}) results from (\ref{bb2}), by taking into
account the anti-symmetry of the operator $b$ with respect to the second and
third variables, and the Sobolev embedding results.$\hfill \hfill
\blacksquare $

\bigskip

\section{Well-posedness of the state equation}

\label{4}

\setcounter{equation}{0} The stochastic differential equation (\ref%
{equation_etat}) has been studied in \cite{CC16}. In this section, we recall
the notion of the solution, as well as the properties of the solution that
will be relevant in the study of the control problem.

Assume that the stochastic noise is represented by%
\begin{equation*}
G(t,y)\,dW_{t}=\sum_{k=1}^{m}g^{k}(t,y)dW^{k}_t,
\end{equation*}%
where%
\begin{eqnarray*}
G(t,y) &=&(g^{1}(\cdot ,y),...,g^{m}(\cdot ,y)):[0,T]\times V\rightarrow
V^{m}=\overbrace{V\times ...\times V}^{m-times}, \\
g^{k}(\cdot ,y) &\in &L^{2}(0,T;V)\qquad \text{for any \ }y\in V.
\end{eqnarray*}%
Let us define the norm $\ ~\left\Vert G(t,y)\right\Vert
_{V}=\sum_{i=1}^{m}\left\Vert g^{i}(t,y)\right\Vert _{V}.$ \ We also assume
that $G(t,y)$ is G\^{a}teaux differentiable in the variable $y$%
\begin{equation*}
\lim_{s\rightarrow 0}\frac{G(t,y+sv)-G(t,y)}{s}=\nabla _{y}G(t,y)v\quad
\text{for each } t\in \lbrack 0,T],
\end{equation*}%
and, moreover, for each $t\in \lbrack 0,T]$\ the function$\ \ \nabla
_{y}G(t,y)$\ is continuous and bounded in the second variable $y\in V,$
namely
\begin{equation*}
||\nabla _{y}G(t,x)-\nabla _{y}G(t,y)||_{V}\rightarrow 0\qquad \text{when\ \
}||x-y||_{V}\rightarrow 0
\end{equation*}%
and%
\begin{equation}
||\nabla _{y}G(t,y)v||_{V}\leq C||v||_{V},\qquad v\in V,  \label{cG}
\end{equation}%
for some positive constant $C.$

Due to Propositions A.2, A.3 of \cite{AH} we have that $G:[0,T]\times
V\rightarrow V^{m}$ \ is Lipschitz and Fr\'{e}chet differentiable in the
second variable $y\in V,$ that is there exists a positive constant $K$ such
that
\begin{align}
\left\Vert G(t,y)-G(t,z)\right\Vert _{V}& \leq K\left\Vert y-z\right\Vert
_{V},  \notag \\
G(t,x+y)-G(t,x)-\nabla _{y}G(t,x)y& =O\left( t,\left\Vert y\right\Vert
_{V}\right) ,\qquad \forall x,y\in V,\;~t\in \lbrack 0,T],  \label{LG}
\end{align}%
having the property%
\begin{equation*}
\frac{O\left( t,s\right) }{s}\rightarrow 0\qquad \text{in }V,\quad \ \text{as%
}\quad s\rightarrow 0.
\end{equation*}

\bigskip

\begin{definition}
Let $~U\in L^{2}(\Omega ,L^{2}(0,T;H^{1}(\mathcal{O})))$ \ and $\ Y_{0}\in
L^{2}(\Omega ,\widetilde{W}).$ A stochastic process $~Y\in L^{2}(\Omega
,L^{\infty }(0,T;\widetilde{W}))~$ is a strong solution of (\ref%
{equation_etat}), if \ for a.e.-$P$ and \ a.e. $t\in (0,T),$ the following
equation holds
\begin{align}
\left( \sigma (Y(t)),\phi \right) =& \int_{0}^{t}\left[ -2\nu \left(
DY(s),D\phi \right) -\left( \mathrm{curl}\,\sigma (Y(s))\times Y(s),\phi
\right) \right] \,ds  \notag \\
& +\left( \sigma (Y(0)),\phi \right) +\int_{0}^{t}\left( U(s),\phi \right)
\,ds+\int_{0}^{t}\left( G(s,Y(s)),\phi \right) \,dW_{s}
\label{var_form_state}
\end{align}%
for all $\phi \in V$, where the nonlinear term is understood in the
following sense
\begin{equation*}
\begin{array}{ll}
\left( \mathrm{curl}\,\sigma (Y(t))\times Y(t),\phi \right) & =b\left( \phi
,Y(t),Y(t)-\alpha \Delta Y(t)\right) -b\left( Y(t),\phi ,Y(t)-\alpha \Delta
Y(t)\right) \vspace{2mm} \\
& =b\left( Y(t),Y(t),\phi \right) -\alpha \left( b\left( \phi ,Y(t),\Delta
Y(t)\right) -b\left( Y(t),\phi ,\Delta Y(t)\right) \right) .%
\end{array}%
\end{equation*}%
\vspace{2mm}\newline
\end{definition}

\bigskip

The following theorem has been proved in the article \cite{CC16}.

\begin{theorem}
\label{the_1} Assume that
\begin{equation}
U\in L^{p}(\Omega ,L^{p}(0,T;H^{1}(\mathcal{O})))\quad \text{and}\quad
Y_{0}\in L^{p}(\Omega ,V)\cap L^{2}(\Omega ,\widetilde{W})\quad \text{for
some }4\leq p<\infty .  \label{au}
\end{equation}%
Then there exists a unique strong solution $Y$ to system (\ref{equation_etat}%
), which belongs to
\begin{equation*}
L^{2}(\Omega ,L^{\infty }(0,T;\widetilde{W}))\cap L^{p}(\Omega ,L^{\infty
}(0,T;V)).
\end{equation*}%
Moreover, the following estimates hold
\begin{equation*}
\mathbb{E}\sup_{s\in \lbrack 0,t]}\left\Vert Y(s)\right\Vert _{V}^{2}+8\nu
\mathbb{E}\int_{0}^{t}\left\Vert DY(s)\right\Vert _{2}^{2}\,\,ds\leq C\left(
\mathbb{E}\left\Vert Y_{0}\right\Vert _{V}^{2}+\mathbb{E}\Vert U\Vert
_{L^{2}(0,t;L^{2})}^{2}+1\right) ,
\end{equation*}%
\begin{equation*}
\mathbb{E}\sup_{s\in \lbrack 0,t]}\left\Vert \mathrm{curl}\,\sigma \left(
Y(s)\right) \right\Vert _{2}^{2}\leq C\left( \mathbb{E}\left\Vert \mathrm{%
curl}\,\sigma (Y_{0})\right\Vert _{2}^{2}+\mathbb{E}\left\Vert U\right\Vert
_{L^{2}(0,t;H^{1})}^{2}+1\right) ,
\end{equation*}%
\begin{equation}
\mathbb{E}\sup_{s\in \lbrack 0,t]}\left\Vert Y(s)\right\Vert _{V}^{p}\leq C%
\mathbb{E}\left\Vert Y_{0}\right\Vert _{V}^{p}+C\,(1+\mathbb{E}%
\int_{0}^{t}\left\Vert U(s)\right\Vert _{2}^{p}\,ds).  \label{lp1}
\end{equation}
\end{theorem}

\bigskip

In Sections \ref{5}-\ref{81} we consider that the data $U,$ $Y_{0}$ satisfy
assumptions (\ref{au}) and denote the unique solution of \ (\ref%
{var_form_state}) by
\begin{equation*}
Y\in L^{p}(\Omega ,L^{\infty }(0,T;V))\cap L^{2}(\Omega ,L^{\infty }(0,T;%
\widetilde{W}))
\end{equation*}%
which satisfies estimates (\ref{lp1}).

We also recall the stability result of solutions for the stochastic
differential system (\ref{equation_etat}) obtained in \cite{CC16}.

\begin{proposition}
\label{Lips} Assume that
\begin{equation*}
U_{1},U_{2}\in L^{p}(\Omega ,L^{p}(0,T;H^{1}(\mathcal{O})))\quad \text{for }%
4\leq p<\infty .
\end{equation*}%
Let
\begin{equation*}
Y_{1},Y_{2}\in L^{2}(\Omega ,L^{\infty }(0,T;\widetilde{W}))
\end{equation*}%
be the corresponding solutions of (\ref{equation_etat}) with the same
initial condition%
\begin{equation*}
Y_{0}\in L^{p}(\Omega ,V)\cap L^{2}(\Omega ,\widetilde{W}).
\end{equation*}

Then there exist positive constants $C$ and $C_{0},$ such that the following
estimate holds for a.e. $t\in (0,T)$
\begin{equation*}
\mathbb{E}\sup_{s\in \lbrack 0,t]}\xi _{0}(s)\left\Vert
Y_{1}(s)-Y_{2}(s)\right\Vert _{W}^{2}\leq C\,\mathbb{E}\int_{0}^{t}\xi
_{0}(s)\left\Vert U_{1}(s)-U_{2}(s)\right\Vert _{2}^{2}\,ds,
\end{equation*}%
with the function $\xi _{0}$ defined as
\begin{equation*}
\xi _{0}(t)=e^{-C_{0}\int_{0}^{t}\left( \left\Vert Y_{1}\right\Vert
_{H^{3}}+\left\Vert Y_{2}\right\Vert _{H^{3}}+1\right) ds}.
\end{equation*}
\end{proposition}

\section{Stochastic linearized state equation}

\bigskip \setcounter{equation}{0}\label{5}

\bigskip

In the study of the control problem the well-posedness of the stochastic
linearized state equation is a crucial step, since such a solution will
correspond to the G\^{a}teaux derivative of the control to state map.
However the mathematical analysis of this equation is not an easy issue.

\bigskip Let us take
\begin{equation}
\Psi \in L^{2}(\Omega \times (0,T)\times \mathcal{O})  \label{5.0}
\end{equation}%
and consider the following linear system%
\begin{equation}
\left\{
\begin{array}{ll}
d\sigma (Z)=(\nu \Delta Z-\mathrm{curl}\,\sigma (Z)\times Y-\mathrm{curl}%
\,\sigma (Y)\times Z-\nabla \pi +\Psi )\,dt\vspace{2mm} &  \\
\qquad \qquad \qquad \qquad \qquad \qquad \qquad +\nabla _{y}G(t,Y)Z\,dW_{t},
& \, \\
\mathrm{div}\,Z=0 & \quad \mbox{in}\ (0,T)\times \mathcal{O},\vspace{2mm} \\
Z\cdot \mathrm{n}=0,\qquad (\mathrm{n}\cdot DZ)\cdot \mathrm{\tau }=0 &
\quad \mbox{on}\ (0,T)\times \Gamma ,\vspace{2mm} \\
Z(0)=0 & \quad \mbox{in}\ \mathcal{O}.%
\end{array}%
\right.  \label{5.3}
\end{equation}

The difficulties are related with the regularity of the coefficients in
equation (\ref{5.3})$_{1}$. \ Firstly we should give a correct meaning to
the solution. Next, taking an appropriate basis on the functional space $%
H^{2}$ we follow the Galerkin method to construct the approximate solution.
Deriving uniform estimates for the Galerkin approximations we will be able
to pass to the limit and obtain the solvability for the stochastic
linearized state system \ (\ref{5.3}). \vspace{2mm}\newline

\begin{definition}
A stochastic process $Z\in L^{2}(0,T;W)$ is a strong solution of (\ref{5.3})
if\ \ for $P-$a.e. $\omega \in \Omega ,$ \ a.e. $t\in (0,T)$ the following
equality holds
\begin{align}
\left( \sigma \left( Z(t)\right) ,\phi \right) & =\int_{0}^{t}\{-2\nu \left(
DZ(s),D\phi \right)  \notag \\
& -b\left( \phi ,Y(s),\sigma (Z(s))\right) +b\left( Y(s),\phi ,\sigma
(Z(s))\right) \vspace{2mm}  \notag \\
& -b\left( \phi ,Z(s),\sigma (Y(s))\right) +b\left( Z(s),\phi ,\sigma
(Y(s))\right) +\left( \Psi (s),\phi \right) \,\}\ ds  \notag \\
& +\int_{0}^{t}(\nabla _{y}G(s,Y(s))Z(s),\phi )\,dW_{s}\qquad
\mbox{  for
any }\ \phi \in V.  \label{ds}
\end{align}%
\
\end{definition}

\bigskip

Let us recall that the proof of Theorem \ref{the_1} in \cite{CC16} has been
given by the Galerkin approximation method taking a basis $\{e_{i}\}\subset
\widetilde{W}$ of eigenfunctions to the injection operator $I:\widetilde{W}%
\hookrightarrow V$. Here, to construct the solution for the stochastic
differential equation (\ref{5.3}) we follow analogous strategy. We define an
appropriate basis which is different of the previous one. \

Since the injection operator $I:W\hookrightarrow V$ is a compact operator,
there exists a basis $\{h_{i}\}\subset W$ of eigenfunctions
\begin{equation}
\left( v,h_{i}\right) _{W}=\mu _{i}\left( v,h_{i}\right) _{V},\qquad \forall
v\in W,\quad i\in \mathbb{N},  \label{mu}
\end{equation}%
being simultaneously an orthonormal basis for $V$. \ The corresponding
sequence $\{\mu _{i}\}$ of eigenvalues fulfills the properties: $\mu _{i}>0$%
, $\forall i\in \mathbb{N}$, and $\mu _{i}\rightarrow \infty $ as $%
i\rightarrow \infty .$ \ We may consider $\{h_{i}\}\subset H^{4},$ taking
into account that the ellipticity of equation (\ref{mu}) increases the
regularity of their solutions (see \cite{BR03}). \vspace{2mm}

Let us consider the basis $\ \{h_{i}\}$ and introduce the Galerkin
approximations of equation (\ref{5.3}). Let $W_{n}=\mathrm{span}%
\,\{h_{1},\ldots ,h_{n}\}$ and define
\begin{equation*}
Z_{n}(t)=\sum_{i=1}^{n}\zeta _{i}(t)h_{i}\quad \text{for each }t\in \lbrack
0,T],
\end{equation*}%
as the solution of the stochastic differential equation%
\begin{equation*}
\begin{array}{l}
d\sigma (Z_{n})=(\nu \Delta Z_{n}-\mathrm{curl}\,\sigma (Z_{n})\times Y-%
\mathrm{curl}\,\sigma (Y)\times Z_{n}-\nabla \pi _{n}+\Psi )\,dt\vspace{2mm}
\\
\qquad \qquad \qquad \qquad \qquad \qquad \qquad +\nabla
_{y}G(t,Y)Z_{n}\,dW_{t}.%
\end{array}%
\end{equation*}%
This means that
\begin{align}
d\left( \sigma \left( Z_{n}\right) ,\phi \right) & =\left( \left( \nu \Delta
Z_{n}-\mathrm{curl}\,\sigma (Z_{n})\times Y-\mathrm{curl}\,\sigma (Y)\times
Z_{n}+\Psi \right) ,\phi \right) \,dt\vspace{2mm}  \notag \\
& +\left( \nabla _{y}G(t,Y)Z_{n},\phi \right) \,dW_{t}\ \qquad
\mbox{  for
each }\ \phi \in W_{n}.  \label{equation_etat_temps_n}
\end{align}%
Equation (\ref{equation_etat_temps_n}) defines a system of stochastic linear
ordinary differential equations, which has a unique solution $Z_{n}$ as an
adapted process in the space $C([0,{T}];W_{n}).$

\bigskip

Let us show the following result.

\begin{proposition}
\label{ex_uniq_lin} Assume that $\Psi $ satisfies assumption (\ref{5.0}).
Then equation (\ref{equation_etat_temps_n}) admits a unique solution $Z_{n}.$
\ Moreover, there exist positive constants $C_{1}$, $C_{2}$ and $C,$ \ which
are independent on the index $n,$ such that the following estimates hold for
each $t\in \lbrack 0,T]$
\begin{eqnarray}
\mathbb{E}\sup_{s\in \lbrack 0,t]}\xi _{1}(s)\left\Vert Z_{n}(s)\right\Vert
_{V}^{2}+\mathbb{E}\int_{0}^{t}\xi _{1}(s)\left\Vert DZ_{n}\right\Vert
^{2}\,\ ds &\leq &C\mathbb{E}\int_{0}^{t}\xi _{1}(s)\left\Vert \Psi
\right\Vert _{2}^{2}\,ds,  \notag \\
\mathbb{E}\sup_{s\in \lbrack 0,t]}\xi _{2}(s)\left\Vert Z_{n}(s)\right\Vert
_{W}^{2} &\leq &C\mathbb{E}\int_{0}^{t}\xi _{2}(s)\left\Vert \Psi
\right\Vert _{2}^{2}\,ds,  \label{in}
\end{eqnarray}%
where the functions $\xi _{1},$ $\xi _{2}$\ are defined as
\begin{equation*}
\xi _{1}(t)=e^{-2C_{1}\int_{0}^{t}\left\Vert Y\right\Vert _{\widetilde{W}%
}ds},\qquad \xi _{2}(t)=e^{-4C_{2}\int_{0}^{t}\left\Vert Y\right\Vert _{%
\widetilde{W}}ds}.
\end{equation*}
\end{proposition}

\textbf{Proof.} \textit{1st step: \ Estimate in the space }$V$\textit{\ for }%
$Z_{n}.$

In order to simplify the notation, let us introduce the function
\begin{equation*}
f(Z_{n})=\nu \Delta Z_{n}-\mathrm{curl}\,\sigma (Z_{n})\times Y-\mathrm{curl}%
\,\sigma (Y)\times Z_{n}+\Psi ,
\end{equation*}%
which belongs to $H^{1}.$ Setting $\phi =h_{i}$ in equation (\ref%
{equation_etat_temps_n}), we obtain
\begin{equation}
d\left( Z_{n},h_{i}\right) _{V}=\left( f(Z_{n}),h_{i}\right) \,dt\vspace{2mm}%
+\left( \nabla _{y}G(t,Y)Z_{n},h_{i}\right) \,dW_{t}.  \label{z_n}
\end{equation}%
By It\^{o}%
\'{}%
s formula we have%
\begin{eqnarray*}
d\left( Z_{n},h_{i}\right) _{V}^{2}=2\left( Z_{n},h_{i}\right) _{V}\left(
f(Z_{n}),h_{i}\right) \,dt &+&\left( \nabla _{y}G(t,Y)Z_{n},h_{i}\right)
^{2}\,dt \\
&+&2\left( Z_{n},h_{i}\right) _{V}\left( \nabla _{y}G(t,Y)Z_{n},h_{i}\right)
\,dW_{t}.
\end{eqnarray*}%
Summing these equalities on $i=1,...,n$ and using the property
\begin{equation*}
2\left( f(Z_{n}),Z_{n}\right) =-4\nu \left\Vert DZ_{n}\right\Vert
_{2}^{2}-2\left( \mathrm{curl}\,\sigma (Z_{n})\times Y-\Psi ,Z_{n}\right) ,
\end{equation*}%
we derive
\begin{align}
d\left\Vert Z_{n}\right\Vert _{V}^{2}& =\left( -4\nu \left\Vert
DZ_{n}\right\Vert _{2}^{2}-2(\mathrm{curl}\,\sigma (Z_{n})\times Y-\Psi
,Z_{n}\,)\right) \,dt\vspace{2mm}  \notag \\
& +\sum_{i=1}^{n}\left( \nabla _{y}G(t,Y)Z_{n},h_{i}\right) ^{2}\,dt+2\left(
\nabla _{y}G(t,Y)Z_{n},Z_{n}\right) \,dW_{t}.  \label{fa}
\end{align}

By \ (\ref{rm2_lin}) we have
\begin{equation}
|(\mathrm{curl}\,\sigma (Z_{n})\times Y,Z_{n}\,)|\leq C_{1}\left\Vert
Y\right\Vert _{\widetilde{W}}\left\Vert Z_{n}\right\Vert _{V}^{2}\vspace{2mm}%
,  \label{u}
\end{equation}%
therefore multiplying (\ref{fa}) by
\begin{equation*}
\xi _{1}(t)=e^{-2C_{1}\int_{0}^{t}\left\Vert Y\right\Vert _{\widetilde{W}}ds}
\end{equation*}%
and using It\^{o}%
\'{}%
s formula, we obtain%
\begin{align*}
d\left( \xi _{1}(t)\left\Vert Z_{n}\right\Vert _{V}^{2}\right) &
=-2C_{1}~\xi _{1}(t)\left\Vert Y\right\Vert _{\widetilde{W}}\left\Vert
Z_{n}\right\Vert _{V}^{2}\,dt \\
& -\xi _{1}(t)\left( 4\nu \left\Vert DZ_{n}\right\Vert _{2}^{2}+2\left(
\mathrm{curl}\,\sigma (Z_{n})\times Y-\Psi ,Z_{n}\right) \right) \,dt \\
& +~\xi _{1}(t)\sum_{i=1}^{n}\left( \nabla _{y}G(t,Y)Z_{n},h_{i}\right)
^{2}\,dt \\
& +2~\xi _{1}(t)\left( \nabla _{y}G(t,Y)Z_{n},Z_{n}\right) \,dW_{t}.
\end{align*}%
Estimate\ (\ref{u}) and the integration over the time variable imply
\begin{align}
~\xi _{1}(t)\left\Vert Z_{n}(t)\right\Vert _{V}^{2}& +4\nu \int_{0}^{t}~\xi
_{1}\left\Vert DZ_{n}\right\Vert _{2}^{2}\ ds\leq 2\int_{0}^{t}~\xi
_{1}\left\vert \left( \Psi ,Z_{n}\right) \right\vert \,ds  \notag \\
& +\int_{0}^{t}~\xi _{1}\sum_{i=1}^{n}\left( \nabla
_{y}G(s,Y)Z_{n},h_{i}\right) ^{2}\,ds  \notag \\
& +2\left\vert \int_{0}^{t}~\xi _{1}\left( \nabla
_{y}G(s,Y)Z_{n},Z_{n}\right) \,dW_{s}\right\vert .  \label{y22}
\end{align}

Let
\begin{equation*}
\tau _{N}(\omega )=\inf \{t\in \lbrack 0,T]:\Vert Z_{n}(t)\Vert _{W}\geq N\}
\end{equation*}%
be the stopping time for fixed $N\in \mathbb{N}.$ The Burkholder-Davis-Gundy
inequality and property (\ref{cG}) give%
\begin{align*}
\mathbb{E}\sup_{s\in \lbrack 0,\tau _{N}\wedge t]}& \left\vert
\int_{0}^{s}\xi _{1}(r)\left( \nabla _{y}G(r,Y)Z_{n},Z_{n}\right)
\,dW_{r}\right\vert \\
& \leq \mathbb{E}\left( \int_{0}^{\tau _{N}\wedge t}\xi _{1}^{2}(s)(\nabla
_{y}G(s,Y)Z_{n},Z_{n})^{2}\,ds\right) ^{\frac{1}{2}} \\
& \leq C\mathbb{E}\sup_{s\in \lbrack 0,\tau _{N}\wedge t]}\sqrt{\xi _{1}(s)}%
\left\Vert Z_{n}\right\Vert _{2}\left( \int_{0}^{\tau _{N}\wedge t}\xi
_{1}(s)\left\Vert Z_{n}\right\Vert _{V}^{2}\,ds\right) ^{\frac{1}{2}} \\
& \leq \varepsilon \,\mathbb{E}\sup_{s\in \lbrack 0,\tau _{N}\wedge t]}\xi
_{1}(s)\Vert Z_{n}\Vert _{V}^{2}+C_{\varepsilon }\mathbb{E}\int_{0}^{\tau
_{N}\wedge t}~\xi _{1}(s)\left\Vert Z_{n}\right\Vert _{V}^{2}\,ds.
\end{align*}%
Substituting this inequality with $\varepsilon =\frac{1}{2}$ in (\ref{y22}),
taking the supremum on $s\in \lbrack 0,\tau _{N}\wedge t]$ and the
expectation, we deduce
\begin{eqnarray*}
\frac{1}{2}\mathbb{E}\sup_{s\in \lbrack 0,\tau _{N}\wedge t]}\xi
_{1}(s)\Vert Z_{n}(s)\Vert _{V}^{2} &+&4\nu \mathbb{E}\int_{0}^{\tau
_{N}\wedge t}\xi _{1}\left\Vert DZ_{n}\right\Vert _{2}^{2}\ ds \\
&\leq &C\mathbb{E}\int_{0}^{\tau _{N}\wedge t}\xi _{1}\left\Vert
Z_{n}\right\Vert _{V}^{2}\,\ ds+C\mathbb{E}\int_{0}^{t}\xi _{1}\left\Vert
\Psi \right\Vert _{2}^{2}\,ds.
\end{eqnarray*}

Hence the function
\begin{equation*}
f(t)=\mathbb{E}\sup_{s\in \lbrack 0,\tau _{N}\wedge t]}\xi _{1}(s)\Vert
Z_{n}(s)\Vert _{V}^{2}
\end{equation*}%
fulfills Gronwall%
\'{}%
s type inequality%
\begin{equation*}
\frac{1}{2}f(t)\leq C\int_{0}^{t}f(s)ds+C\mathbb{E}\int_{0}^{t}\xi
_{1}(s)\left\Vert \Psi \right\Vert _{2}^{2}\,ds.
\end{equation*}%
The application of Gronwall%
\'{}%
s inequality gives
\begin{eqnarray}
\mathbb{E}\sup_{s\in \lbrack 0,\tau _{N}\wedge t]}\xi _{1}(s)\left\Vert
Z_{n}(s)\right\Vert _{V}^{2} &+&8\nu \mathbb{E}\int_{0}^{\tau _{N}\wedge
t}\xi _{1}\left\Vert DZ_{n}\right\Vert _{2}^{2}\ ds  \notag \\
&\leq &C\mathbb{E}\int_{0}^{t}\xi _{1}\left\Vert \Psi \right\Vert
_{2}^{2}\,ds.  \label{pas1}
\end{eqnarray}

\bigskip \textit{2nd step: Estimate in the space }$W$\textit{\ for }$Z_{n}.$

Let $\ ~\tilde{f}_{n}$ \ and $\ ~\widetilde{G}_{n}$ \ be the solutions of
problem (\ref{GS_NS}) for $~f=f(Z_{n})~$ and $\ ~f=\nabla _{y}G(t,Y)Z_{n},~$
\ respectively. The following relations hold
\begin{equation}
(\tilde{f}_{n},h_{i})_{V}=(f(Z_{n}),h_{i}),\qquad (\widetilde{G}%
_{n},h_{i})_{V}=(\nabla _{y}G(t,Y)Z_{n},h_{i})\qquad \text{ for each }i.\
\label{sss}
\end{equation}

Using relations (\ref{mu}) and multiplying (\ref{z_n}) by $\mu _{i}$, we
deduce%
\begin{equation*}
d\left( Z_{n},h_{i}\right) _{W}=(\tilde{f}_{n},h_{i})_{W}\,dt+(\widetilde{G}%
_{n},h_{i})_{W}\,dW_{t}.
\end{equation*}%
Hence the It\^{o} formula gives
\begin{align*}
d\left( Z_{n},h_{i}\right) _{W}^{2}& =2\left( Z_{n},h_{i}\right) _{W}(\tilde{%
f}_{n},h_{i})_{W}\,dt+(\widetilde{G}_{n},h_{i})_{W}^{2}\,dt \\
& +2\left( Z_{n},h_{i}\right) _{W}(\widetilde{G}_{n},h_{i})_{W}\,dW_{t}.
\end{align*}%
Multiplying these equalities by $\frac{1}{\mu _{i}}$ and summing over $%
i=1,\dots ,n$, we obtain
\begin{equation*}
d\left\Vert Z_{n}\right\Vert _{W}^{2}=2(\tilde{f}_{n},Z_{n})_{W}\,dt+||%
\widetilde{G}_{n}||_{W}^{2}\,dt+2(\widetilde{G}_{n},Z_{n})_{W}\,dW_{t}.
\end{equation*}%
Then, by the relation between the norms of the spaces $W$ and $V$, induced
by the inner products (\ref{ss}), we derive
\begin{align*}
d\left( \left\Vert Z_{n}\right\Vert _{W}^{2}\right) =2((\mathbb{P}\sigma (%
\tilde{f}_{n}),\mathbb{P}\sigma (Z_{n}))& +(\tilde{f}_{n},Z_{n})_{V})\,dt+||%
\widetilde{G}_{n}||_{W}^{2}\,dt \\
& +2\left( (\mathbb{P}\sigma (\widetilde{G}_{n}),\mathbb{P}\sigma (Z_{n}))+(%
\widetilde{G}_{n},Z_{n})_{V}\right) \,dW_{t}.
\end{align*}%
Using (\ref{sss}), this equality can be written as
\begin{align}
d\left( \left\Vert Z_{n}\right\Vert _{W}^{2}\right) & =2((f(Z_{n}),\mathbb{P}%
\sigma (Z_{n}))+(f(Z_{n}),Z_{n}))\,dt+||\widetilde{G}_{n}||_{W}^{2}\,dt
\notag \\
& +2\left[ (\nabla _{y}G(t,Y)Z_{n},\mathbb{P}\sigma (Z_{n}))+(\nabla
_{y}G(t,Y)Z_{n},Z_{n})\right] \,dW_{t}.  \label{s4}
\end{align}

\bigskip

First let us estimate the most difficult term $~\left( \mathrm{curl}\,\sigma
(Z_{n})\times Y,\mathbb{P}\sigma (Z_{n})\right) ,~$ being a part of the
expression $~(f(Z_{n}),\mathbb{P}\sigma (Z_{n})).$ If we denote $~R=\sigma
(Z_{n})~$ and $~Q=\mathbb{P}\sigma (Z_{n})-\sigma (Z_{n}),~$ applying (\ref%
{bb1}) we get%
\begin{align*}
\left\vert \left( \mathrm{curl}\,\sigma (Z_{n})\times Y,\mathbb{P}\sigma
(Z_{n})\right) \right\vert & \leq \left\vert \left( \mathrm{curl}\,R\times
Y,Q\right) \right\vert +\left\vert \left( \mathrm{curl}\,R\times Y,R\right)
\right\vert \\
& \leq \left\vert b\left( Q,Y,R\right) -b\left( Y,Q,R\right) \right\vert
+\left\vert b\left( R,Y,R\right) -b\left( Y,R,R\right) \right\vert \\
& \leq \left\vert b\left( Q,Y,R\right) |+|b\left( Y,Q,R\right) \right\vert
+\left\vert b\left( R,Y,R\right) \right\vert \\
& \leq C\left\Vert Y\right\Vert _{1,\infty }\left( \left\Vert Q\right\Vert
_{H^{1}}\left\Vert R\right\Vert _{2}+\left\Vert R\right\Vert _{2}^{2}\right)
.
\end{align*}%
Hence
\begin{equation*}
\left\vert \left( \mathrm{curl}\,\sigma (Z_{n})\times Y,\mathbb{P}\sigma
(Z_{n})\right) \right\vert \leq C\left\Vert Y\right\Vert _{1,\infty
}\left\Vert Z_{n}\right\Vert _{W}^{2}\leq \widetilde{C}_{2}\left\Vert
Y\right\Vert _{\widetilde{W}}\left\Vert Z_{n}\right\Vert _{W}^{2},
\end{equation*}%
by using estimates (\ref{sigma_psigma1})-(\ref{sigma_psigma2}) for$\
v=Z_{n}. $ Obviously we have
\begin{equation*}
\left\vert \left( \mathrm{curl}\,\sigma (Y)\times Z_{n},\mathbb{P}\sigma
(Z_{n})\right) \right\vert \leq C\left\Vert Y\right\Vert _{H^{3}}\left\Vert
Z_{n}\right\Vert _{\infty }\left\Vert \mathbb{P}\sigma (Z_{n})\right\Vert
_{W}^{2}\leq \overline{C}_{2}\left\Vert Y\right\Vert _{\widetilde{W}%
}\left\Vert Z_{n}\right\Vert _{W}^{2}.
\end{equation*}

Recalling that $\ f(Z_{n})=\nu \Delta Z_{n}-\mathrm{curl}\,\sigma
(Z_{n})\times Y-\mathrm{curl}\,\sigma (Y)\times Z_{n}+\Psi ~$ and (\ref{u}),
we deduce%
\begin{eqnarray}
\left\vert \left( f(Z_{n}),\mathbb{P}\sigma (Z_{n})\right) \right\vert &\leq
&C\left\Vert Z_{n}\right\Vert _{W}^{2}+C_{2}\left\Vert Y\right\Vert _{%
\widetilde{W}}\left\Vert Z_{n}\right\Vert _{W}^{2}+\left\vert \left( \Psi ,%
\mathbb{P}\sigma (Z_{n})\right) \right\vert ,  \notag \\
|\left( f(Z_{n}),Z_{n}\right) | &\leq &-2\nu \left\Vert DZ_{n}\right\Vert
_{2}^{2}+C_{1}\left\Vert Y\right\Vert _{\widetilde{W}}\left\Vert
Z_{n}\right\Vert _{V}^{2}+|\left( \Psi ,Z_{n}\right) |  \notag \\
&\leq &C_{2}\left\Vert Y\right\Vert _{\widetilde{W}}\left\Vert
Z_{n}\right\Vert _{W}^{2}+|\left( \Psi ,Z_{n}\right) |  \label{s444}
\end{eqnarray}%
for
\begin{equation*}
C_{2}=\max (\widetilde{C}_{2},\overline{C}_{2},C_{1}).
\end{equation*}

As in previous considerations, if we take the function $~~\xi
_{2}(t)=e^{-4C_{2}\int_{0}^{t}\left\Vert Y\right\Vert _{\widetilde{W}}ds}~$
and use Ito%
\'{}%
s formula, then (\ref{s4}) and (\ref{s444}) imply%
\begin{align*}
d\left( \xi _{2}(t)\left\Vert Z_{n}\right\Vert _{W}^{2}\right) & \leq C\xi
_{2}(t)\left( \left\Vert Z_{n}\right\Vert _{W}^{2}+\left\vert \left( \Psi ,%
\mathbb{P}\sigma (Z_{n})\right) \right\vert +|\left( \Psi ,Z_{n}\right)
|\right) \,dt \\
& +\xi _{2}(t)||\widetilde{G}_{n}||_{W}^{2}\,dt \\
& +2\xi _{2}(t)\left( \nabla _{y}G(t,Y)Z_{n},\mathbb{P}\sigma
(Z_{n})+Z_{n})\right) \,dW_{t}.
\end{align*}%
Therefore the integration of this inequality over the time variable gives
\begin{eqnarray}
~\xi _{2}(t)\left\Vert Z_{n}(t)\right\Vert _{W}^{2} &\leq &C\int_{0}^{t}~\xi
_{2}\left( \left\Vert Z_{n}\right\Vert _{W}^{2}+\left\Vert \Psi \right\Vert
_{2}^{2}\right) \,ds+\int_{0}^{t}~\xi _{2}||\widetilde{G}_{n}||_{W}^{2}\,ds
\notag \\
&&+2\left\vert \int_{0}^{t}~\xi _{2}\left( \nabla _{y}G(s,Y)Z_{n},\mathbb{P}%
\sigma (Z_{n})+Z_{n}\right) \,\,dW_{s}\right\vert .  \label{s3}
\end{eqnarray}

Since $\widetilde{G}_{n}$\ \ is the solution of the Stokes type equation (%
\ref{GS_NS}), then using (\ref{cG}) we have%
\begin{equation}
||\widetilde{G}_{n}||_{W}\ \leq C||Z_{n}||_{V}.  \label{s2}
\end{equation}%
Moreover the Burkholder-Davis-Gundy inequality gives
\begin{align}
\mathbb{E}\sup_{s\in \lbrack 0,\tau _{N}\wedge t]}& \left\vert
\int_{0}^{s}\xi _{2}(r)(\nabla _{y}G(r,Y)Z_{n},\mathbb{P}\sigma
(Z_{n})+Z_{n})\,\,dW_{r}\right\vert  \notag \\
& \leq \mathbb{E}\left( \int_{0}^{\tau _{N}\wedge t}\xi _{2}^{2}(s)(\nabla
_{y}G(s,Y)Z_{n},\mathbb{P}\sigma (Z_{n})+Z_{n})^{2}\,ds\right) ^{\frac{1}{2}}
\notag \\
& \leq C\mathbb{E}\sup_{s\in \lbrack 0,\tau _{N}\wedge t]}\sqrt{\xi _{2}(s)}%
\left\Vert Z_{n}\right\Vert _{W}\left( \int_{0}^{\tau _{N}\wedge t}\xi
_{2}(s)\left\Vert Z_{n}\right\Vert _{2}^{2}\,ds\right) ^{\frac{1}{2}}  \notag
\\
& \leq \varepsilon \,\mathbb{E}\sup_{s\in \lbrack 0,\tau _{N}\wedge t]}\xi
_{2}(s)\Vert Z_{n}\Vert _{W}^{2}+C_{\varepsilon }\mathbb{E}\int_{0}^{\tau
_{N}\wedge t}~\xi _{2}(s)\left\Vert Z_{n}\right\Vert _{W}^{2}\,ds.
\label{s1}
\end{align}%
Substituting inequalities (\ref{s2}) and (\ref{s1}) with chosen $\varepsilon
=\frac{1}{2}$ \ in (\ref{s3}) and taking the supremum on $s\in \lbrack
0,\tau _{N}\wedge t]$, we derive
\begin{equation*}
\frac{1}{2}\mathbb{E}\sup_{s\in \lbrack 0,\tau _{N}\wedge t]}\xi
_{2}(s)\Vert Z_{n}(s)\Vert _{W}^{2}\leq C\mathbb{E}\int_{0}^{\tau _{N}\wedge
t}\xi _{2}(s)\left\Vert Z_{n}\right\Vert _{W}^{2}\,ds+C\mathbb{E}%
\int_{0}^{t}\xi _{2}(s)\left\Vert \Psi \right\Vert _{2}^{2}\,ds.
\end{equation*}%
The application of Gronwall%
\'{}%
s inequality gives
\begin{equation}
\mathbb{E}\sup_{s\in \lbrack 0,\tau _{N}\wedge t]}\xi _{2}(s)\left\Vert
Z_{n}(s)\right\Vert _{W}^{2}\leq C\mathbb{E}\int_{0}^{t}\xi _{2}\left\Vert
\Psi \right\Vert _{2}^{2}ds.  \label{pas2}
\end{equation}%
Hence passing to the limit $N\rightarrow \infty $ in (\ref{pas1}) and (\ref%
{pas2}), we derive estimates (\ref{in})$.\hfill \hfill \blacksquare $

\bigskip

As a direct consequence of Proposition \ref{ex_uniq_lin} we derive the
following existence result for the stochastic linearized state equation (\ref%
{5.3}).

\begin{theorem}
\label{the_2} Assume that $\Psi $ satisfies assumption (\ref{5.0}).\ Then
there exists a unique solution $Z$\ to system (\ref{5.3}), such that%
\begin{equation*}
Z\in L^{\infty }(0,T;W)\quad \text{for \ \ }\ P\text{-a.e. \ }\ \omega \in
\Omega ,
\end{equation*}%
satisfying the following a priori estimates%
\begin{eqnarray*}
\mathbb{E}\sup_{s\in \lbrack 0,t]}\xi _{1}(s)\left\Vert Z(s)\right\Vert
_{V}^{2}+\mathbb{E}\int_{0}^{t}\xi _{1}\left\Vert Z\right\Vert _{V}^{2}\,\
ds &\leq &C\mathbb{E}\int_{0}^{t}\xi _{1}\left\Vert \Psi \right\Vert
_{2}^{2}\,ds\quad \text{for a.e.\ \ }\ t\in (0,T), \\
\mathbb{E}\sup_{s\in \lbrack 0,T]}\xi _{2}(s)\left\Vert Z(s)\right\Vert
_{W}^{2} &\leq &C\mathbb{E}\int_{0}^{T}\xi _{2}\left\Vert \Psi \right\Vert
_{2}^{2}\,ds.
\end{eqnarray*}
\end{theorem}

\textbf{Proof.} \vspace{2mm} Since
\begin{equation*}
\sup_{t\in \lbrack 0,T]}\left\Vert Y(t)\right\Vert _{\widetilde{W}}^{2}\leq
C(\omega )\quad \text{for all }\omega \in \Omega \backslash A,\quad \text{%
where }P(A)=0
\end{equation*}%
by (\ref{ss}) and (\ref{lp1}), therefore there exists a positive constant $%
K(\omega ),$ which depends only on $\omega \in \Omega \backslash A$ and
satisfies%
\begin{equation}
0<K(\omega )\leq \xi _{i}(t)\leq 1\quad \text{for all }\omega \in \Omega
\backslash A,\quad t\in \lbrack 0,T],\quad i=1,2.  \label{q2}
\end{equation}

By a priori estimates of Lemma \ref{ex_uniq_lin} we have%
\begin{equation*}
\mathbb{E}\sup_{t\in \lbrack 0,T]}\xi _{1}(t)\left\Vert Z_{n}(t)\right\Vert
_{V}^{2}\leq C,\qquad \mathbb{E}\sup_{t\in \lbrack 0,T]}\xi
_{2}(t)\left\Vert Z_{n}(t)\right\Vert _{W}^{2}\leq C.
\end{equation*}%
Hence using (\ref{q2}) there exists a suitable subsequence $Z_{n}$, indexed
by the same index $n$, such that%
\begin{eqnarray}
\sqrt{\xi _{1}(t)}Z_{n} &\rightharpoonup &\sqrt{\xi _{1}(t)}Z\qquad
\mbox{ *-weakly in
}\ L^{2}(\Omega ,L^{\infty }(0,T;V)),  \notag \\
S_{n}=\sqrt{\xi _{2}(t)}Z_{n} &\rightharpoonup &S=\sqrt{\xi _{2}(t)}Z\qquad
\mbox{ *-weakly in
}\ L^{2}(\Omega ,L^{\infty }(0,T;W)).  \label{lim1}
\end{eqnarray}

Since $Z_{n}$ solves (\ref{equation_etat_temps_n}), then we have%
\begin{align*}
d\left( \sigma \left( Z_{n}\right) ,\phi \right) & =\left( \left( \nu \Delta
Z_{n}-\mathrm{curl}\,\sigma (Z_{n})\times Y-\mathrm{curl}\,\sigma (Y)\times
Z_{n}+\Psi \right) ,\phi \right) \,dt\vspace{2mm} \\
& +\left( \nabla _{y}G(t,Y)Z_{n},\phi \right) \,dW_{t}\text{\qquad for each}%
\quad \text{ }\phi \in W_{n}.
\end{align*}%
\smallskip\ \ Using Ito%
\'{}%
s formula for $~\xi (t)=\sqrt{\xi _{2}(t)}=e^{-2C_{2}\int_{0}^{t}\left\Vert
Y\right\Vert _{\widetilde{W}}ds},~$ we get
\begin{align*}
d\left( \sigma \left( S_{n}\right) ,\phi \right) & =[\left( \left\{ \nu
\Delta S_{n}-\mathrm{curl}\,\sigma (S_{n})\times Y-\mathrm{curl}\,\sigma
(Y)\times S_{n}+\xi \Psi \right\} ,\phi \right) \\
& -2C_{2}\left\Vert Y\right\Vert _{\widetilde{W}}\ \ \left( \sigma \left(
S_{n}\right) ,\phi \right) ]\vspace{2mm}\,dt\vspace{2mm}+\left( \nabla
_{y}G(t,Y)S_{n},\phi \right) \,dW_{t},
\end{align*}%
that is
\begin{align*}
\left( \sigma \left( S_{n}(t)\right) ,\phi \right) &
=\int_{0}^{t}[-2C_{2}\left\Vert Y\right\Vert _{\widetilde{W}}\ \left( \sigma
\left( S_{n}\right) ,\phi \right) \ -2\nu \left( DS_{n},D\phi \right) \\
& -b\left( \phi ,Y,\sigma (S_{n})\right) +b\left( Y,\phi ,\sigma
(S_{n})\right) \vspace{2mm} \\
& -b\left( \phi ,S_{n},\sigma (Y)\right) +b\left( S_{n},\phi ,\sigma
(Y)\right) +\left( \xi \Psi ,\phi \right) ]\ ds \\
& +\int_{0}^{t}(\nabla _{y}G(s,Y)S_{n},\phi )\,dW_{s}.
\end{align*}%
Let us consider an arbitrary $\varphi =\varphi (\omega )\in L^{2}(\Omega ).$
Multiplying this equality by $\varphi $ and taking the expectation, we
derive
\begin{align*}
\mathbb{E}\varphi (\omega )\left( \sigma \left( S_{n}(t)\right) ,\phi
\right) & =\mathbb{E}\varphi (\omega )\left\{ \int_{0}^{t}[-2C_{2}\left\Vert
Y\right\Vert _{\widetilde{W}}\ \left( \sigma \left( S_{n}\right) ,\phi
\right) \ -2\nu \left( DS_{n},D\phi \right) \right. \\
& -b\left( \phi ,Y,\sigma (S_{n})\right) +b\left( Y,\phi ,\sigma
(S_{n})\right) \vspace{2mm} \\
& -b\left( \phi ,S_{n},\sigma (Y)\right) +b\left( S_{n},\phi ,\sigma
(Y)\right) +\left( \xi \Psi ,\phi \right) \,]\ ds \\
& +\left. \int_{0}^{t}(\nabla _{y}G(s,Y)S_{n}(s),\phi )\,dW_{s}\right\}
\qquad \text{for \ }\forall t\in \lbrack 0,T].
\end{align*}%
Using that the right side of the last equation is continuous in the time
variable $t\in \lbrack 0,T]$ and applying (\ref{lim1}) and Proposition A.3,
p.93 of \cite{B99}, we pass to the limit $n\rightarrow \infty $ in this
equality and deduce%
\begin{align*}
\mathbb{E}\varphi (\omega )\left( \sigma \left( S(t)\right) ,\phi \right) & =%
\mathbb{E}\varphi (\omega )\left\{ \int_{0}^{t}[-2C_{2}\left\Vert
Y\right\Vert _{\widetilde{W}}\ \left( \sigma \left( S\right) ,\phi \right) \
-2\nu \left( DS,D\phi \right) \right. \\
& -b\left( \phi ,Y,\sigma (S)\right) +b\left( Y,\phi ,\sigma (S)\right)
\vspace{2mm} \\
& -b\left( \phi ,S,\sigma (Y)\right) +b\left( S,\phi ,\sigma (Y)\right)
+\left( \xi \Psi ,\phi \right) \,]\ ds \\
& +\left. \int_{0}^{t}(\nabla _{y}G(s,Y)S,\phi )\,dW_{s}\right\} .
\end{align*}

Since $\varphi \in L^{2}(\Omega )$ is arbitrary and $~S(t)=\xi (t)Z(t)$,$~$
then we have\ the validity of the equality
\begin{align*}
\xi (t)\left( \sigma \left( Z(t)\right) ,\phi \right) & =\left\{
\int_{0}^{t}\xi \lbrack -2C_{2}\left\Vert Y\right\Vert _{\widetilde{W}}\
\left( \sigma \left( S\right) ,\phi \right) ds-2\nu \left( DZ,D\phi \right)
\right. \\
& -b\left( \phi ,Y,\sigma (Z)\right) +b\left( Y,\phi ,\sigma (Z)\right)
\vspace{2mm} \\
& -b\left( \phi ,Z,\sigma (Y)\right) +b\left( Z,\phi ,\sigma (Y)\right)
+\left( \Psi ,\phi \right) \,]\ ds \\
& +\left. \int_{0}^{t}\xi (s)(\nabla _{y}G(s,Y)Z,\phi )\,dW_{s}\right\}
\qquad \text{for \ }\forall t\in \lbrack 0,T]\text{\ \ \ and\ }\ P\text{%
-a.e. \ }\ \omega \in \Omega .
\end{align*}%
Hence%
\begin{align*}
d\left[ \xi \left( \sigma \left( Z\right) ,\phi \right) \right] & =\xi
\lbrack -2C_{2}\left\Vert Y\right\Vert _{\widetilde{W}}\ \left( \sigma
\left( Z\right) ,\phi \right) -2\nu \left( DZ,D\phi \right) \\
& -b\left( \phi ,Y,\sigma (Z)\right) +b\left( Y,\phi ,\sigma (Z)\right)
\vspace{2mm} \\
& -b\left( \phi ,Z,\sigma (Y)\right) +b\left( Z,\phi ,\sigma (Y)\right)
+\left( \Psi ,\phi \right) \,]\ dt \\
& +(\nabla _{y}G(s,Y)Z,\phi )\,dW_{s}\}\qquad \text{for a.e. }t\in \lbrack
0,T]\text{\ \ \ and\ }\ P\text{-a.e. \ }\ \omega \in \Omega .
\end{align*}%
Moreover if we use Ito's formula $\ $%
\begin{equation*}
\ d\left( \sigma \left( Z\right) ,\phi \right) =d\left[ \widehat{\xi }\xi
\left( \sigma \left( Z\right) ,\phi \right) \right] =\xi \left( \sigma
\left( Z\right) ,\phi \right) d\widehat{\xi }+\widehat{\xi }d\left[ \xi
\left( \sigma \left( Z\right) ,\phi \right) \right]
\end{equation*}%
for $\widehat{\xi }(t)=e^{2C_{2}\int_{0}^{t}\left\Vert Y\right\Vert _{%
\widetilde{W}}ds},$\ we derive
\begin{align*}
d\left( \sigma \left( Z\right) ,\phi \right) & =\bigl\{\lbrack -2\nu \left(
DZ,D\phi \right) \\
& -b\left( \phi ,Y,\sigma (Z)\right) +b\left( Y,\phi ,\sigma (Z)\right)
\vspace{2mm} \\
& -b\left( \phi ,Z,\sigma (Y)\right) +b\left( Z,\phi ,\sigma (Y)\right)
+\left( \Psi ,\phi \right) \,]\ dt \\
& +(\nabla _{y}G(t,Y)Z,\phi )\,dW_{t}\bigl\}\qquad \text{for a.e. }t\in
\lbrack 0,T]\text{,\ }\ P\text{-a.e. \ }\ \omega \in \Omega
\end{align*}%
and for each $\ \ \phi \in W.\quad $ \ Therefore the stochastic process $%
Z\in L^{2}(0,T;W)$ is a solution of (\ref{5.3}) in the sense of equality \ (%
\ref{ds}).$\hfill \hfill \blacksquare $

\bigskip

\section{G\^ateaux differentiability of the control-to-state mapping}

\bigskip \label{6}

It is well known that the G\^{a}teaux derivative of the control-to-state is
fundamental to deduce the necessary optimality conditions.

The goal of this section is to prove that the G\^{a}teaux derivative of $Y$
at the point $U$, in the direction $\Psi $, is defined as the solution of
the linear system (\ref{5.3}). We know that due to Theorem \ref{the_2} this
system has a unique solution. As a consequence of this theorem and
Proposition \ref{Lips} we have the following result.

\begin{proposition}
\label{Gat} Given $U,$ $Y_{0},$ satisfying (\ref{au}) and%
\begin{equation*}
\Psi \in L^{p}(\Omega ,L^{p}(0,T;H^{1}(\mathcal{O})))\quad \text{for \ }%
4\leq p<\infty ,
\end{equation*}%
let us consider%
\begin{equation*}
U_{\rho }=U+\rho \Psi ,\qquad \forall \rho \in (0,1).
\end{equation*}%
If $Y$ and $Y_{\rho }$ are the solutions of (\ref{equation_etat})
corresponding to $(U,Y_{0})$ and $(U_{\rho },Y_{0}),$ then the following
representation holds
\begin{equation}
Y_{\rho }=Y+\rho \,Z+\rho \,\delta _{\rho }\quad \mbox{   with  }\quad
\lim_{\rho \rightarrow 0}\sup_{t\in \lbrack 0,T]}\left\Vert \delta _{\rho
}\right\Vert _{V}^{2}=0\quad P\text{-a.e. \ }\ \omega \in \Omega ,
\label{gateau_1}
\end{equation}%
where
\begin{equation*}
Z\in L^{\infty }(0,T;W)\quad \text{for \ \ }\ P\text{-a.e. \ }\ \omega \in
\Omega ,
\end{equation*}%
is the solution of (\ref{5.3}), satisfying the estimates of Theorem \ref%
{the_2}.
\end{proposition}


\textbf{Proof.} It is straightforward to verify that $~Z_{\rho }=\frac{%
Y_{\rho }-Y}{\rho }~$ satisfies the equation
\begin{align*}
d\sigma \left( Z_{\rho }\right) & =\left( \nu \Delta Z_{\rho }-\mathrm{curl}%
\,\sigma \left( Z_{\rho }\right) \times Y-\mathrm{curl}\,\sigma \left(
Y_{\rho }\right) \times Z_{\rho }-\nabla \pi _{\rho }+\Psi \right) \,dt%
\vspace{2mm} \\
& +\frac{1}{\rho }\left( G(t,Y_{\rho })-G(t,Y)\right) \,dW_{t}.
\end{align*}

If we consider the solution $Z$ of (\ref{5.3}), then \ $\delta _{\rho
}=Z_{\rho }-Z$ is the solution of the equation
\begin{align}
d\sigma \left( \delta _{\rho }\right) & =\left( \nu \Delta \delta _{\rho }-%
\mathrm{curl}\,\sigma \left( \delta _{\rho }\right) \times Y-\mathrm{curl}%
\,\sigma \left( Y_{\rho }\right) \times \delta _{\rho }-\mathrm{curl}%
\,\sigma \left( Y_{\rho }-Y\right) \times Z-\nabla \left( \pi _{\rho }-\pi
\right) \right) \,dt  \notag \\
& +R\,dW_{t},  \label{G_2}
\end{align}%
where
\begin{equation}
R=\left[ \frac{1}{\rho }\left( G(t,Y_{\rho })-G(t,Y)\right) -\nabla
_{y}G(t,Y)Z\right] .  \label{RR}
\end{equation}

We apply the operator $\left( I-\alpha \mathbb{P}\Delta \right) ^{-1}$ to
equation (\ref{G_2}) and deduce a stochastic differential equation for \ $%
\delta _{\rho }$, then the It\^{o} formula gives
\begin{align}
d(\left\Vert \delta _{\rho }\right\Vert _{V})+2\nu \left\Vert D\delta _{\rho
}\right\Vert _{2}^{2}& =-2\left( \left( \mathrm{curl}\,\sigma \left( \delta
_{\rho }\right) \times Y,\delta _{\rho }\right) +\left( \mathrm{curl}%
\,\sigma \left( Y_{\rho }\right) \times \delta _{\rho },\delta _{\rho
}\right) \right) \,dt  \notag \\
& -2\left( \mathrm{curl}\,\sigma \left( Y_{\rho }-Y\right) \times Z,\delta
_{\rho }\right) \,dt+2\left( R,\delta _{\rho }\right) \,dW_{t}  \notag \\
& +\Vert \widetilde{G}_{\rho }-\widetilde{G}\Vert _{V}\,dt,  \label{energy_r}
\end{align}%
where $\widetilde{G}_{\rho }$ and $\widetilde{G}$ are the solutions of the
modified Stokes problem (\ref{GS_NS}) with $f$ replaced by $\frac{1}{\rho }%
\left( G(t,Y)-G(t,Y_{\rho })\right) $ and $\nabla _{y}G(t,Y)Z$, respectively.

Applying (\ref{LG}) we have
\begin{equation*}
R=\nabla _{y}G(t,Y)\delta _{\rho }+o\left( t,\left\Vert Z+\delta _{\rho
}\right\Vert _{V}\right) \text{\qquad with }o\left( t,\left\Vert Z+\delta
_{\rho }\right\Vert _{V}\right) \rightarrow 0\qquad \text{in }V,\quad \
\text{as}\quad \rho \rightarrow 0
\end{equation*}%
for $~t\in \lbrack 0,T]~$ and obtain
\begin{equation*}
\Vert \widetilde{G}_{\rho }-\widetilde{G}\Vert _{V}^{2}\leq \Vert R\Vert
_{2}^{2}\leq C\left\Vert \delta _{\rho }\right\Vert _{V}^{2}+o(t,\left\Vert
Z+\delta _{\rho }\right\Vert _{V}^{2}).
\end{equation*}%
Using property (\ref{bb1}), we have%
\begin{equation*}
\left( \mathrm{curl}\,\sigma \left( Y_{\rho }\right) \times \delta _{\rho
},\delta _{\rho }\right) =b\left( \delta _{\rho },\delta _{\rho },\sigma
\left( Y_{\rho }\right) \right) -b\left( \delta _{\rho },\delta _{\rho
},\sigma \left( Y_{\rho }\right) \right) =0
\end{equation*}%
and applying estimate (\ref{rm2_lin}), we get%
\begin{equation*}
\left\vert \left( \mathrm{curl}\,\sigma \left( \delta _{\rho }\right) \times
Y,\delta _{\rho }\right) \right\vert \leq C\left\Vert Y\right\Vert _{%
\widetilde{W}}\left\Vert \delta _{\rho }\right\Vert _{V}^{2}.
\end{equation*}%
In addition, using the Sobolev embedding results, \ we estimate the term
\begin{align*}
|\left( \mathrm{curl}\,\sigma \left( Y_{\rho }-Y\right) \times Z,\delta
_{\rho }\right) |& =\left\vert b\left( \delta _{\rho },Z,\sigma \left(
Y_{\rho }-Y\right) \right) -b\left( Z,\delta _{\rho },\sigma \left( Y_{\rho
}-Y\right) \right) \right\vert \\
& \leq \left( \left\Vert \delta _{\rho }\right\Vert _{4}\left\Vert \nabla
Z\right\Vert _{4}+\left\Vert Z\right\Vert _{\infty }\left\Vert \nabla \delta
_{\rho }\right\Vert _{2}\right) \left\Vert \sigma \left( Y_{\rho }-Y\right)
\right\Vert _{2} \\
& \leq C\left\Vert Z\right\Vert _{H^{2}}\left\Vert \delta _{\rho
}\right\Vert _{V}\left\Vert \sigma \left( Y_{\rho }-Y\right) \right\Vert _{2}
\\
& \leq C\left\Vert Z\right\Vert _{W}^{2}\left\Vert \delta _{\rho
}\right\Vert _{V}^{2}+C\left\Vert \sigma \left( Y_{\rho }-Y\right)
\right\Vert _{2}^{2}.
\end{align*}%
Introducing the above deduced relations in (\ref{energy_r}), we obtain
\begin{align*}
d(\left\Vert \delta _{\rho }\right\Vert _{V}^{2})& \leq C_{4}\left(
\left\Vert Y\right\Vert _{\widetilde{W}}^{2}+\left\Vert Z\right\Vert
_{W}^{2}+1\right) \left\Vert \delta _{\rho }\right\Vert _{V}^{2}+C\left\Vert
\sigma \left( Y_{\rho }-Y\right) \right\Vert _{2}^{2}\,dt\vspace{2mm} \\
& +\left( R,\delta _{\rho }\right) \,dW_{t}.
\end{align*}%
Considering the function
\begin{equation*}
\beta (t)=e^{-\widetilde{C}_{4}\int_{0}^{t}\left( \left\Vert Y\right\Vert _{%
\widetilde{W}}^{2}+\left\Vert Y_{\rho }\right\Vert _{\widetilde{W}%
}^{2}+\left\Vert Z\right\Vert _{W}^{2}+1\right) ds}\qquad \text{for }%
\widetilde{C}_{4}=\max \left\{ C_{0},C_{4}\right\} ,
\end{equation*}%
the It\^{o} formula yields
\begin{equation}
\beta (t)\left\Vert \delta _{\rho }(t)\right\Vert _{V}^{2}\leq
C\int_{0}^{t}\beta (s)\left\Vert \sigma \left( Y_{\rho }-Y\right)
\right\Vert _{2}^{2}\,ds+\int_{0}^{t}\beta (s)\left( R,\delta _{\rho
}\right) \,dW_{s}.  \label{G_est}
\end{equation}%
Moreover the Burkholder-Davis-Gundy inequality, the Young inequality \ and (%
\ref{RR}) give%
\begin{align*}
\mathbb{E}\sup_{s\in \lbrack 0,t]}& \left\vert \int_{0}^{s}\beta
(r)(R,\delta _{\rho })\,dW_{r}\right\vert \vspace{2mm}\leq \mathbb{E}\left(
\int_{0}^{t}\beta (r)^{2}\left( R,\delta _{\rho }\right) ^{2}\,ds\right)
^{1/2}\vspace{2mm} \\
& \leq \varepsilon \mathbb{E}\sup_{s\in \lbrack 0,t]}\beta (s)\left\Vert
\delta _{\rho }\right\Vert _{V}^{2}+C_{\varepsilon }\mathbb{E}%
\int_{0}^{t}\beta (s)\left[ \left\Vert \delta _{\rho }\right\Vert
_{V}^{2}+o(s,\left\Vert Z+\delta _{\rho }\right\Vert _{V}^{2})\,\right] ds.
\end{align*}%
Substituting this inequality with $\varepsilon =\frac{1}{2}$ in (\ref{G_est}%
), taking the supremum on the time interval $[0,t]$ and the expectation, we
\ obtain
\begin{eqnarray*}
\frac{1}{2}\mathbb{E}\sup_{s\in \lbrack 0,t]}\beta (s)\left\Vert \delta
_{\rho }(s)\right\Vert _{V}^{2} &\leq &C\mathbb{E}\int_{0}^{t}\beta
(s)\left\Vert \delta _{\rho }\right\Vert _{V}^{2}\,ds+C\int_{0}^{t}\beta
(s)\left\Vert \sigma \left( Y_{\rho }-Y\right) \right\Vert _{2}^{2}\,ds \\
&&+C\mathbb{E}\int_{0}^{t}\beta (s)o(s,\rho \left\Vert Z+\delta _{\rho
}\right\Vert _{V}^{2})\,ds.
\end{eqnarray*}%
Applying Gronwall's inequality, we deduce%
\begin{eqnarray*}
\mathbb{E}\sup_{s\in \lbrack 0,t]}\beta (s)\left\Vert \delta _{\rho
}\right\Vert _{V}^{2} &\leq &C\,\mathbb{E}\int_{0}^{t}\beta (s)\left\Vert
\sigma \left( Y_{\rho }-Y\right) \right\Vert _{2}^{2}\,ds+C\mathbb{E}%
\int_{0}^{t}\beta (s)o(s,\rho \left\Vert Z+\delta _{\rho }\right\Vert
_{V}^{2})\,ds \\
&\leq &CT\,\mathbb{E}\sup_{s\in \lbrack 0,t]}\xi _{0}(s)\left\Vert \sigma
\left( Y_{\rho }-Y\right) \right\Vert _{2}^{2}\,+C\mathbb{E}%
\int_{0}^{t}\beta (s)o(s,\rho \left\Vert Z+\delta _{\rho }\right\Vert
_{V}^{2})\,ds,
\end{eqnarray*}%
where%
\begin{equation*}
\xi _{0}(t)=e^{-C_{0}\int_{0}^{t}\left( \left\Vert Y_{\rho }\right\Vert _{%
\widetilde{W}}^{2}+\left\Vert Y\right\Vert _{\widetilde{W}}^{2}\right) ds}
\end{equation*}%
and $C_{0}$ is defined in Proposition \ref{Lips}. The result of Proposition %
\ref{Lips} is valid for $Y_{1}=Y_{\rho },$ $Y_{2}=Y.$ \ Therefore we obtain%
\begin{eqnarray*}
\mathbb{E}\sup_{s\in \lbrack 0,t]}\beta (s)\left\Vert \delta _{\rho
}\right\Vert _{V}^{2} &\leq &CT\,\mathbb{E}\sup_{s\in \lbrack 0,t]}\xi
_{0}(s)\left\Vert \sigma \left( Y_{\rho }-Y\right) \right\Vert _{2}^{2}\,ds
\\
&\leq &\rho C\,\mathbb{E}\int_{0}^{t}\xi _{0}\left\Vert \Psi \right\Vert
_{2}^{2}\,ds+C\mathbb{E}\int_{0}^{t}o(s,\rho \left\Vert Z+\delta _{\rho
}\right\Vert _{V}^{2})\,ds.
\end{eqnarray*}%
Since the right hand side in this inequality converges to zero as $\rho
\rightarrow 0$, then applying the Lebesgue theorem, we deduce (\ref{gateau_1}%
).$\hfill \hfill \blacksquare $

\bigskip

\bigskip

\bigskip Let us assume that the Lagrangian $~L:[0,T]\times V\times
\widetilde{W}\rightarrow \mathbb{R}^{+}~$ and the function$\ ~h:\widetilde{W}%
\rightarrow \mathbb{R}^{+},~$ involved in the cost functional (\ref{cost}),
satisfy the following hypotheses:

H1) $L(t,u,y)$ \ and $h(y)$ \ are G\^{a}teaux differentiable on $y$ for any
fixed $u\in V$, $\ \;t\in \lbrack 0,T];$

H2) there exist positive constants $C,$ such that
\begin{align}
||\nabla _{u}L(t,u,y)||_{2}+||\nabla _{y}L\left( t,u,y\right) ||_{2}& \leq
C(1+\left\Vert u\right\Vert _{V}+\left\Vert y\right\Vert _{\widetilde{W}}),
\notag \\
||\nabla _{y}h\left( t,y\right) ||_{2}& \leq C(1+\left\Vert y\right\Vert _{%
\widetilde{W}}),\qquad \forall u\in V,\;y\in \widetilde{W},\;t\in \lbrack
0,T];  \label{lip1}
\end{align}

H3) the $L^{1}-$integrability on the time for $L$ and $h$%
\begin{equation}
L(\cdot ,u,y),~h(y)\in L^{1}(0,T),\qquad \forall u\in V,\;y\in \widetilde{W}.
\label{lip2}
\end{equation}

\bigskip

As a direct consequence of Proposition \ref{Gat} and hypotheses H1)-H3) we
easily derive the following result on the variation for the cost functional (%
\ref{cost}).

\begin{proposition}
\label{Gat1} Let $U,$ $Y_{0},$ $\Psi $\ and $U_{\rho }=U+\rho \Psi $ \
satisfy the conditions of Proposition \ref{Gat}. Let the cost functional (%
\ref{cost}) fulfill hypotheses H1)-H3), then
\begin{equation*}
\begin{array}{ll}
J\left( U_{\rho },Y_{\rho }\right) = & J\left( U,Y\right) +\rho \,\mathbb{E}%
\int_{0}^{T}\left\{ \left( \nabla _{u}L(t,U,Y),\Psi \right) +\left( \nabla
_{y}L(t,U,Y),Z\right) \right\} \,dt\vspace{2mm} \\
& +\mathbb{E}\left( \nabla h(Y_{T}),Z\right) _{V}+o(\rho ),%
\end{array}%
\end{equation*}%
where $Y$, $Y_{\rho }$ are the solutions of (\ref{equation_etat}),
corresponding to $(U,Y_{0})$, $(U_{\rho },Y_{0})$ and $Z$ is the solution of
(\ref{5.3}).
\end{proposition}

\bigskip

\bigskip

\bigskip

\begin{remark}
As a consequence of (\ref{lip1}) we have that $L\left( t,\cdot ,\cdot
\right) $ and $h$ are Lipschitz continuous
\begin{align*}
|L(t,u_{1},y_{1})-L(t,u_{2},y_{2})|& \leq C\left( 1+\max_{i=1,2}\left\Vert
u_{i}\right\Vert _{V}+\max_{i=1,2}\left\Vert y_{i}\right\Vert _{\widetilde{W}%
}\right) \left( \left\Vert u_{1}-u_{2}\right\Vert _{V}+\left\Vert
y_{1}-y_{2}\right\Vert _{\widetilde{W}}\right) , \\
|h(y_{1})-h(y_{2})|& \leq C\left( 1+\max_{i=1,2}\left\Vert y_{i}\right\Vert
_{\widetilde{W}}\right) \left( \left\Vert y_{1}-y_{2}\right\Vert _{%
\widetilde{W}}\right) ,\qquad \forall u_{i}\in V,\;y_{i}\in \widetilde{W}%
,\;t\in \lbrack 0,T].
\end{align*}
\end{remark}

\section{Stochastic backward adjoint equation}

\setcounter{equation}{0}

\label{7}

\bigskip

The aim of this section is to prove the existence of the adjoint stochastic
process $(p,q)$, which is related to the G\^{a}teaux derivative of the
control to state mapping through a duality condition. Let us consider the
following backward stochastic system
\begin{equation}
\left\{
\begin{array}{ll}
d\left( \sigma \left( p\right) ,\phi \right) =\left( -\nu \Delta p-\mathrm{%
curl}\,\sigma (Y)\times p+\mathrm{curl}\left( \sigma \left( Y\times p\right)
\right) ,\phi \right) \,dt\vspace{2mm} &  \\
\quad \qquad \qquad -\displaystyle\left( \nabla _{y}L\left( t,U,Y\right)
+\nabla _{y}G(t,Y)^{T}q,\phi \right) \,dt+\left( \sigma \left( q\right)
,\phi \right) \,dW_{t} & \quad \text{for each \ }\phi \in W, \\
\mathrm{div}\,p=0 & \quad \mbox{in}\ (0,T)\times \mathcal{O},\vspace{2mm} \\
p\cdot \mathrm{n}=0,\qquad (\mathrm{n}\cdot Dp)\cdot \mathrm{\tau }=0 &
\quad \mbox{on}\ (0,T)\times \Gamma ,\vspace{2mm} \\
p(T)=\nabla h\left( Y(T)\right) & \quad \mbox{in}\ \mathcal{O},%
\end{array}%
\right.  \label{p_n}
\end{equation}%
where $Y$ is the unique solution of (\ref{var_form_state}).

\bigskip In all this section we assume that the cost functional (\ref{cost})
fulfill hypotheses H1)-H3).

\subsection{Solvability of the adjoint equation}


\label{101}

\bigskip

The solution of (\ref{p_n}) is understood in the following sense.

\begin{definition}
A stochastic process
\begin{equation*}
(p,q)\in L^{\infty }(0,T;W)\times L^{2}(0,T;W)\qquad P-\text{a.e. in }\Omega
\end{equation*}%
is a solution of (\ref{p_n}) if\ \ for $P-$a.e. $\omega \in \Omega $ \ and\
a.e. $t\in (0,T)$ the following equality holds
\begin{align}
\left( \sigma \left( \nabla h\left( Y(T)\right) \right) ,\phi \right)
-\left( \sigma \left( p(t)\right) ,\phi \right) & =\int_{t}^{T}\bigl\{2\nu
\left( Dp(s),D\phi \right)  \notag \\
& -b\left( \phi ,p(s),\sigma (Y(s))\right) +b\left( p(s),\phi ,\sigma
(Y(s))\right) \vspace{2mm}  \notag \\
& +b\left( p(s),Y(s),\sigma (\phi (s))\right) -b\left( Y(s),p(s),\sigma
(\phi (s))\right) \vspace{2mm}  \notag \\
& +\left( \nabla _{y}L\left( t,U,Y\right) +\nabla _{y}G(t,Y)^{T}q,\phi
\right) \bigr\}\,ds\vspace{2mm}  \notag \\
& +\int_{t}^{T}\left( \sigma \left( q\right) ,\phi \right) \,dW_{s}\text{%
\qquad for each\quad }\phi \in W.  \label{equation_etat1}
\end{align}
\end{definition}

\bigskip

We construct the solution of system (\ref{p_n}), using Galerkin
approximations. \ Let us consider the basis $\ \{h_{i}\},$ defined in (\ref%
{mu}),\ and consider the space $W_{n}=\mathrm{span}\,\{h_{1},\ldots
,h_{n}\}. $ Let
\begin{equation}
p_{n}=\sum_{i=1}^{n}\mathfrak{p}_{i}(t)h_{i}\quad \text{and}\quad
q_{n}=\sum_{i=1}^{n}\mathfrak{q}_{i}(t)h_{i}  \label{pnqn}
\end{equation}%
be the solution of the Galerkin approximations of system (\ref{p_n}), being
the following backward stochastic system%
\begin{equation}
\left\{
\begin{array}{ll}
d\left( \sigma \left( p_{n}\right) ,\phi \right) =\left( -\nu \Delta p_{n}-%
\mathrm{curl}\,\sigma (Y)\times p_{n}+\mathrm{curl}\left( \sigma \left(
Y\times p_{n}\right) \right) ,\phi \right) \,dt\vspace{2mm} &  \\
\quad \qquad \qquad -\displaystyle\left( \nabla _{y}L\left( t,U,Y\right)
+\nabla _{y}G(t,Y)^{T}q_{n},\phi \right) \,dt+\left( \sigma \left(
q_{n}\right) ,\phi \right) \,dW_{t} & \text{\quad for each\quad }\phi \in
W_{n}. \\
\mathrm{div}\,p_{n}=0 & \quad \mbox{in}\ (0,T)\times \mathcal{O},\vspace{2mm}
\\
p_{n}\cdot \mathrm{n}=0,\qquad (\mathrm{n}\cdot Dp_{n})\cdot \mathrm{\tau }=0
& \quad \mbox{on}\ (0,T)\times \Gamma ,\vspace{2mm} \\
p_{n}(T)=\nabla h\left( Y(T)\right) & \quad \mbox{in}\ \mathcal{O}.%
\end{array}%
\right.  \label{pn}
\end{equation}

In the next proposition we show the existence of the pair $(p_{n},q_{n}).$

\begin{proposition}
\label{est_pq} There exists a unique solution
\begin{equation*}
(p_{n},q_{n})\in L^{\infty }(0,T;W)\times L^{2}(0,T;W)\qquad P-\text{a.e. in
}\Omega
\end{equation*}%
of the backward stochastic system (\ref{pn}). Moreover there exists a
positive constant $C_{3}$, such that the pair $(p_{n},q_{n})$ satisfies the
following estimate for a.e. $t\in (0,T)$
\begin{eqnarray}
\frac{1}{2}\mathbb{E}\,\sup_{s\in \lbrack t,T]}\xi _{3}(s)\left\Vert
p_{n}(s)\right\Vert _{W}^{2} &+&\mathbb{E}\int_{t}^{T}\xi _{3}(s)(\frac{\nu
}{\alpha }\left\Vert \mathbb{P}\sigma (p_{n})\right\Vert _{2}^{2}+4\nu
\left\Vert Dp_{n}\right\Vert _{2}^{2}+\frac{1}{2}\left\Vert q_{n}\right\Vert
_{W}^{2})\,ds  \notag \\
&\leq &C(\mathbb{E}\,\left\Vert \nabla h\left( Y(T)\right) \right\Vert
_{W}^{2}+\mathbb{E}\,\int_{t}^{T}\xi _{3}(s)\left\Vert \nabla _{y}L\left(
s,U,Y\right) \right\Vert _{2}^{2}\,ds)  \label{pnn}
\end{eqnarray}%
with the function%
\begin{equation*}
\xi _{3}(t)=e^{-C_{3}\int_{t}^{T}\left\Vert Y\right\Vert _{\widetilde{W}}ds}.
\end{equation*}
\end{proposition}

\textbf{Proof.} \textit{1st step: Existence of approximate solutions $\left(
p_{n},q_{n}\right) .$}\ \ Equation (\ref{pn}) defines a system of stochastic
backward linear ordinary differential equations, which has a unique solution
$(p_{n},q_{n})$ \ as an adapted process in the space $C([{0,T}];W_{n}).$

\textit{2nd step: Estimate in the space }$\mathit{W}$\textit{\ for the
approximate solutions $\left( p_{n},q_{n}\right) .$} \ Setting $\phi =h_{i}$
in equation (\ref{pn})$_{1}$, we obtain
\begin{equation}
d\left( p_{n},h_{i}\right) _{V}=\left( f(p_{n}),h_{i}\right) \,dt\vspace{2mm}%
-\left( \nabla _{y}G(t,Y)^{T}q_{n},h_{i}\right) \,dt+\left(
q_{n},h_{i}\right) _{V}\,dW_{t}  \label{p_1n}
\end{equation}%
with
\begin{equation}
f(p_{n})=-\nu \Delta p_{n}-\mathrm{curl}\,\sigma (Y)\times p_{n}+\mathrm{curl%
}\left( \sigma \left( Y\times p_{n}\right) \right) -\nabla _{y}L\left(
t,U,Y\right) .  \label{D_f}
\end{equation}

Let $\tilde{f}_{n},$ $\widetilde{G}_{n}$ \ be the solutions of the modified
Stokes problem (\ref{GS_NS}) for $f=f(p_{n}),$ $f=\nabla
_{y}G(t,Y)^{T}q_{n}, $ \ respectively. The following relations hold
\begin{equation}
(\tilde{f}_{n},h_{i})_{V}=(f(p_{n}),h_{i}),\qquad (\widetilde{G}%
_{n},h_{i})_{V}=(\nabla _{y}G(t,Y)^{T}q_{n},h_{i})\qquad \text{ for each}\ i.
\label{sssp}
\end{equation}%
Using relations (\ref{mu}) and multiplying (\ref{p_1n}) by $\mu _{i}$, we
deduce%
\begin{equation*}
d\left( p_{n},h_{i}\right) _{W}=\left( \tilde{f}_{n},h_{i}\right)
_{W}\,dt-\left( \widetilde{G}_{n},h_{i}\right) _{W}\,dt+\left(
q_{n},h_{i}\right) _{W}\,dW_{t}.
\end{equation*}%
On the other hand, the It\^{o} formula gives
\begin{align*}
d\left( p_{n},h_{i}\right) _{W}^{2}& =2\left( p_{n},h_{i}\right) _{W}(\tilde{%
f}_{n},h_{i})_{W}\,dt-2\left( p_{n},h_{i}\right) _{W}(\widetilde{G}%
_{n},h_{i})_{W}\,dt \\
& +2\left( p_{n},h_{i}\right) _{W}\left( q_{n},h_{i}\right)
_{W}\,dW_{t}+\left( q_{n},h_{i}\right) _{W}^{2}\,dt.
\end{align*}%
Multiplying these equalities by $\frac{1}{\mu _{i}}$ and summing over $%
i=1,\dots ,n$, we obtain
\begin{equation*}
d\left\Vert p_{n}\right\Vert _{W}^{2}=2(\tilde{f}_{n},p_{n})_{W}\,dt-2(%
\widetilde{G}_{n},p_{n})_{W}\,dt+2\left( q_{n},p_{n}\right)
_{W}\,dW_{t}+||q_{n}||_{W}^{2}\,dt.
\end{equation*}%
Then, by the relation between the inner products for the spaces $W$ and $V,$
we derive%
\begin{align*}
d\left\Vert p_{n}\right\Vert _{W}^{2}=2((\mathbb{P}\sigma (\tilde{f}_{n}),%
\mathbb{P}\sigma (p_{n}))& +(\tilde{f}_{n},p_{n})_{V})\,dt-2\left[ (\mathbb{P%
}\sigma (\widetilde{G}_{n}),\mathbb{P}\sigma (p_{n}))+(\widetilde{G}%
_{n},p_{n})_{V}\right] \,dt \\
& +2\left( (\mathbb{P}\sigma (q_{n}),\mathbb{P}\sigma
(p_{n}))+(q_{n},p_{n})_{V}\right) \,dW_{t}+||q_{n}||_{W}^{2}\,dt.
\end{align*}%
Using (\ref{sssp}), this equality can be written as
\begin{align*}
d\left\Vert p_{n}\right\Vert _{W}^{2}& =2\left[ (f(p_{n}),\mathbb{P}\sigma
(p_{n}))+(f(p_{n}),p_{n})\,\right] dt-2[(\nabla _{y}G(t,Y)^{T}q_{n},\mathbb{P%
}\sigma (p_{n})) \\
& +(\nabla _{y}G(t,Y)^{T}q_{n},p_{n})]\,dt+2\left( \sigma \left(
q_{n}\right) ,\mathbb{P}\sigma (p_{n})+p_{n})\right)
\,dW_{t}+||q_{n}||_{W}^{2}\,dt.
\end{align*}

\bigskip In what follows we estimate the terms of the right hand side in
this equality. From (\ref{D_f}) we have that
\begin{equation*}
\left( f(p_{n}),p_{n}\right) =2\nu \left\Vert Dp_{n}\right\Vert
_{2}^{2}+\left( \mathrm{curl}\left( \sigma \left( Y\times p_{n}\right)
\right) ,p_{n}\right) -(\nabla _{y}L\left( t,U,Y\right) ,p_{n})
\end{equation*}%
and%
\begin{align*}
\left( f(p_{n}),\mathbb{P}\sigma \left( p_{n}\right) \right) & =\ \frac{\nu
}{\alpha }\left\Vert \mathbb{P}\sigma \left( p_{n}\right) \right\Vert
_{2}^{2}-\frac{\nu }{\alpha }\left( p_{n},\mathbb{P}\sigma \left(
p_{n}\right) \right) -\left( \mathrm{curl}\,\sigma (Y)\times p_{n},\mathbb{P}%
\sigma \left( p_{n}\right) \right) \\
& +\left( \mathrm{curl}\left( \sigma \left( Y\times p_{n}\right) \right) ,%
\mathbb{P}\sigma \left( p_{n}\right) \right) -\left( \nabla _{y}L\left(
t,U,Y\right) ,\mathbb{P}\sigma \left( p_{n}\right) \right) .
\end{align*}%
Since
\begin{align}
\left\vert \left( \mathrm{curl}\,\sigma (Y)\times p_{n},\mathbb{P}\sigma
(p_{n})\right) \right\vert & \leq \left\Vert \mathrm{curl}\,\sigma
(Y)\right\Vert _{2}\left\Vert p_{n}\right\Vert _{\infty }\left\Vert \sigma
(p_{n})\right\Vert _{2}\leq C\left\Vert \mathrm{curl}\,\sigma (Y)\right\Vert
_{2}\left\Vert p_{n}\right\Vert _{H^{2}}^{2}  \notag \\
& \leq C\left\Vert Y\right\Vert _{\widetilde{W}}\left\Vert p_{n}\right\Vert
_{W}^{2}  \label{+}
\end{align}%
and taking into account Lemma \ref{prop_a_adj_2}, we obtain
\begin{align}
\left\vert \left( \mathrm{curl}\left( \sigma \left( Y\times p_{n}\right)
\right) ,\mathbb{P}\sigma \left( p_{n}\right) \right) \right\vert & \leq
\left\vert b\left( \sigma (p_{n}),Y,\mathbb{P}\sigma \left( p_{n}\right)
\right) \right\vert +\left\vert b\left( Y,\mathbb{P}\sigma (p_{n})-\sigma
\left( p_{n}\right) ,\sigma (p_{n}\right) \right\vert \vspace{2mm}  \notag \\
& +\left\vert b\left( \sigma \left( Y\right) ,p_{n},\mathbb{P}\sigma \left(
p_{n}\right) \right) +b\left( p_{n},\sigma (Y),\mathbb{P}\sigma \left(
p_{n}\right) \right\vert \right) \vspace{2mm}  \notag \\
& +\left\vert b\left( Y,p_{n},\mathbb{P}\sigma \left( p_{n}\right) \right)
-b\left( p_{n},Y,\mathbb{P}\sigma \left( p_{n}\right) \right) \right\vert
\vspace{1mm}  \notag \\
& +\left\vert 2\alpha \displaystyle\sum_{i=1}^{2}\left( b(\tfrac{\partial
p_{n}}{\partial x_{i}},\tfrac{\partial Y}{\partial x_{i}},\mathbb{P}\sigma
\left( p_{n}\right) )-b(\tfrac{\partial Y}{\partial x_{i}},\tfrac{\partial
p_{n}}{\partial x_{i}},\mathbb{P}\sigma \left( p_{n}\right) )\right)
\right\vert \vspace{1mm}  \notag \\
& \leq \left\Vert \sigma (p_{n})\right\Vert _{2}\left\Vert \nabla
Y\right\Vert _{\infty }\left\Vert \mathbb{P}\sigma \left( p_{n}\right)
)\right\Vert _{2}  \notag \\
& +\left\Vert Y\right\Vert _{\infty }\left\Vert \mathbb{P}\sigma \left(
p_{n}\right) -\sigma (p_{n})\right\Vert _{H^{1}}\left\Vert \sigma
(p_{n})\right\Vert _{2}\vspace{1mm}  \notag \\
& +\left\Vert \sigma (Y)\right\Vert _{4}\left\Vert \nabla p_{n}\right\Vert
_{4}\left\Vert \mathbb{P}\sigma \left( p_{n}\right) \right\Vert
_{2}+\left\Vert p_{n}\right\Vert _{\infty }\left\Vert \nabla \sigma
(Y)\right\Vert _{2}\left\Vert \mathbb{P}\sigma \left( p_{n}\right)
\right\Vert _{2}\vspace{1mm}  \notag \\
& +\left\Vert Y\right\Vert _{\infty }\left\Vert \nabla p_{n}\right\Vert
_{2}\left\Vert \mathbb{P}\sigma \left( p_{n}\right) \right\Vert
_{2}+\left\Vert p_{n}\right\Vert _{2}\left\Vert \nabla Y\right\Vert _{\infty
}\left\Vert \mathbb{P}\sigma \left( p_{n}\right) \right\Vert _{2}\vspace{1mm}
\notag \\
& +2\alpha \displaystyle\sum_{i=1}^{2}\left[ \left\Vert \tfrac{\partial p_{n}%
}{\partial x_{i}}\right\Vert _{4}\left\Vert \nabla (\tfrac{\partial Y}{%
\partial x_{i}})\right\Vert _{4}\left\Vert \mathbb{P}\sigma \left(
p_{n}\right) \right\Vert _{2}\right.  \notag \\
& \left. +\left\Vert \tfrac{\partial Y}{\partial x_{i}}\right\Vert _{\infty
}\left\Vert \nabla \left( \tfrac{\partial p_{n}}{\partial x_{i}}\right)
\right\Vert _{2}\left\Vert \mathbb{P}\sigma \left( p_{n}\right) \right\Vert
_{2}\right]  \notag \\
& \leq C\left\Vert Y\right\Vert _{\widetilde{W}}\left\Vert p_{n}\right\Vert
_{W}^{2}.  \label{8}
\end{align}%
By \ (\ref{rm2_adj}) we also have
\begin{equation}
|\left( \mathrm{curl}\,\left( \sigma \left( Y\times p_{n}\right) \right)
,p_{n}\right) |\leq C\left\Vert Y\right\Vert _{\widetilde{W}}\left\Vert
p_{n}\right\Vert _{V}^{2}\vspace{2mm}.  \label{up}
\end{equation}%
As a consequence of estimates (\ref{+})-(\ref{up}), there exists a fixed
positive constant $C_{3},$ satisfying
\begin{align*}
2\left( f(p_{n}),\mathbb{P}\sigma (p_{n})\right) +2\left(
f(p_{n}),p_{n}\right) & \geq \frac{2\nu }{\alpha }\left\Vert \mathbb{P}%
\sigma (p_{n})\right\Vert _{W}^{2}+4\nu \left\Vert Dp_{n}\right\Vert
_{2}^{2}-\frac{2\nu }{\alpha }\left( p_{n},\mathbb{P}\sigma \left(
p_{n}\right) \right) \\
& -C_{3}\left\Vert Y\right\Vert _{\widetilde{W}}\left\Vert p_{n}\right\Vert
_{W}^{2}-2\left( \nabla _{y}L\left( t,U,Y\right) ,\mathbb{P}\sigma
(p_{n})+p_{n}\right) .
\end{align*}%
Reasoning as above in Lemma \ref{ex_uniq_lin}, if we take the function $~\xi
_{3}(t)=e^{-C_{3}\int_{t}^{T}\left\Vert Y\right\Vert _{\widetilde{W}}ds}$,
the Ito formula yields
\begin{align*}
d\left( \xi _{3}(t)\left\Vert p_{n}\right\Vert _{W}^{2}\right) & \geq \xi
_{3}(t)\left( \frac{2\nu }{\alpha }\left\Vert \mathbb{P}\sigma
(p_{n})\right\Vert _{2}^{2}+4\nu \left\Vert Dp_{n}\right\Vert _{2}^{2}-\frac{%
2\nu }{\alpha }\left( p_{n},\mathbb{P}\sigma \left( p_{n}\right) \right)
\right) \,dt \\
& -2\xi _{3}(t)\left( \nabla _{y}L\left( t,U,Y\right) ,\mathbb{P}\sigma
(p_{n})+p_{n}\right) \,dt \\
& -2\xi _{3}(t)\left( \nabla _{y}G(t,Y)^{T}q_{n},\mathbb{P}\sigma
(p_{n})+p_{n}\right) \,dt \\
& +\xi _{3}(t)\left( 2\left( \sigma (q_{n}),\mathbb{P}\sigma
(p_{n})+p_{n})\right) \,dW_{t}+||q_{n}||_{W}^{2}\,dt\right) .
\end{align*}%
Integrating this inequality over the time interval $(t,T),$ we obtain
\begin{align*}
\xi _{3}(t)\left\Vert p_{n}(t)\right\Vert _{W}^{2}& +\int_{t}^{T}\xi
_{3}\left( \frac{2\nu }{\alpha }\left\Vert \mathbb{P}\sigma
(p_{n})\right\Vert _{2}^{2}+4\nu \left\Vert Dp_{n}\right\Vert
_{2}^{2}+\left\Vert q_{n}\right\Vert _{W}^{2}\right) \,ds \\
& \leq \left\Vert p_{n}(T)\right\Vert _{W}^{2}-2\int_{t}^{T}\xi _{3}\left(
\sigma \left( q_{n}\right) ,\mathbb{P}\sigma (p_{n})+p_{n})\right) \,dW_{s}
\\
& +\frac{2\nu }{\alpha }\int_{t}^{T}\xi _{3}\left( p_{n},\mathbb{P}\sigma
\left( p_{n}\right) \right) \,ds \\
& +2\int_{t}^{T}\xi _{3}\left( \nabla _{y}L\left( s,U,Y\right) ,\mathbb{P}%
\sigma (p_{n})+p_{n}\right) \,ds \\
& +2\int_{t}^{T}\xi _{3}\left( \nabla _{y}G(s,Y)^{T}q_{n},\mathbb{P}\sigma
(p_{n})+p_{n}\right) \,ds.
\end{align*}%
We can apply Young's inequality $2ab\leq \varepsilon a^{2}+b^{2}/\varepsilon
$ to the last three terms. \ Therefore for an appropriate chosen of $%
\varepsilon $ in each of these three terms, accounting (\ref{cG}), we easily
derive
\begin{align}
\xi _{3}(t)\left\Vert p_{n}(t)\right\Vert _{W}^{2}& +\int_{t}^{T}\xi
_{3}\left( \frac{\nu }{\alpha }\left\Vert \mathbb{P}\sigma
(p_{n})\right\Vert _{2}^{2}+4\nu \left\Vert Dp_{n}\right\Vert _{2}^{2}+\frac{%
1}{2}\left\Vert q_{n}\right\Vert _{W}^{2}\right) \,ds  \notag \\
& \leq \,\left\Vert p_{n}(T)\right\Vert _{W}^{2}+C\,\int_{t}^{T}\xi
_{3}\left\Vert p_{n}\right\Vert _{W}^{2}\,ds+C\,\int_{t}^{T}\xi
_{3}\left\Vert \nabla _{y}L\left( s,U,Y\right) \right\Vert _{2}^{2}ds  \notag
\\
& -2\int_{t}^{T}\xi _{3}\left( \sigma \left( q_{n}\right) ,\mathbb{P}\sigma
(p_{n})+p_{n})\right) \,\,dW_{s}.  \label{qe}
\end{align}

Let us define the stopping time%
\begin{equation*}
{\tau _{N}(\omega )=\sup \{t\in \lbrack 0,T]:\min \{\Vert p_{n}(t)\Vert
_{W},~\left\Vert q_{n}\right\Vert _{W}\}\geq N\}}\text{\quad\ for {fixed}}{\
\ N\in \mathbb{N}}.
\end{equation*}%
Let $~\mathbf{1}_{[\tau _{N}\vee t,T]}~$ be the characteristic function of
the time interval $~[\tau _{N}\vee t,T]~$ for each $\omega \in \Omega $. \ \
\ The application of the expectation in (\ref{qe}) implies that the function
$f~(t)=\mathbb{E}\mathbf{1}_{[\tau _{N}\vee t,T]}\xi _{3}(t)\left\Vert
p_{n}(t)\right\Vert _{W}^{2}~$ \ \ fulfills the Gronwall type inequality%
\begin{eqnarray*}
f(t) &+&\int_{t}^{T}\mathbb{E}\mathbf{1}_{[\tau _{N}\vee t,T]}\xi
_{3}(s)\left( \frac{\nu }{\alpha }\left\Vert \mathbb{P}\sigma
(p_{n})\right\Vert _{2}^{2}+4\nu \left\Vert Dp_{n}\right\Vert _{2}^{2}+\frac{%
1}{2}\left\Vert q_{n}\right\Vert _{W}^{2}\right) \,ds \\
&\leq &\int_{t}^{T}f(s)\ ds+\mathbb{E}\left\Vert p_{n}(T)\right\Vert
_{W}^{2}+C\mathbb{E}\int_{t}^{T}~\xi _{3}(s)\left\Vert \nabla _{y}L\left(
s,U,Y\right) \right\Vert _{2}^{2}\,ds,
\end{eqnarray*}%
hence%
\begin{eqnarray}
f(t) &+&\int_{t}^{T}\mathbb{E}\mathbf{1}_{[\tau _{N}\vee t,T]}\xi
_{3}(s)\left( \frac{\nu }{\alpha }\left\Vert \mathbb{P}\sigma
(p_{n})\right\Vert _{2}^{2}+4\nu \left\Vert Dp_{n}\right\Vert _{2}^{2}+\frac{%
1}{2}\left\Vert q_{n}\right\Vert _{W}^{2}\right) \,ds  \notag \\
&\leq &C(\mathbb{E}\left\Vert \nabla h(T)\right\Vert _{W}^{2}+\mathbb{E}%
\int_{t}^{T}~\xi _{3}(s)\left\Vert \nabla _{y}L\left( s,U,Y\right)
\right\Vert _{2}^{2}\,ds)\qquad \text{for \ }t\in \lbrack 0,T].  \label{qe1}
\end{eqnarray}

The Burkholder-Davis-Gundy inequality gives
\begin{align}
2\mathbb{E}\sup_{s\in \lbrack \tau _{N}\vee t,T]}& \left\vert
\int_{s}^{T}\xi _{3}(r)(\sigma \left( q_{n}\right) ,\mathbb{P}\sigma
(p_{n})+p_{n})\,\,dW_{r}\right\vert  \notag \\
& \leq 2\mathbb{E}\left( \int_{\tau _{N}\vee t}^{T}\xi _{3}^{2}(s)|(q_{n},%
\mathbb{P}\sigma (p_{n})+p_{n})|^{2}\,ds\right) ^{\frac{1}{2}}  \notag \\
& \leq 2\mathbb{E}\sup_{s\in \lbrack \tau _{N}\vee t,T]}\sqrt{\xi _{3}(s)}%
\left\Vert \mathbb{P}\sigma (p_{n})+p_{n}\right\Vert _{2}\left( \int_{\tau
_{N}\vee t}^{T}\xi _{3}(s)\Vert \sigma \left( q_{n}\right) \Vert
_{2}^{2}\,ds\right) ^{\frac{1}{2}}  \notag \\
& \leq \frac{1}{2}\,\mathbb{E}\sup_{s\in \lbrack \tau _{N}\vee t,T]}\xi
_{3}(s)\Vert p_{n}\Vert _{W}^{2}+2\mathbb{E}\int_{t}^{T}~\xi _{3}(s)\Vert
q_{n}\Vert _{W}^{2}\,ds.  \label{qe2}
\end{align}%
Therefore taking the supremum on the time interval $~s\in \lbrack \tau
_{N}\vee t,T]$~ in (\ref{qe}), the expectation and applying inequalities (%
\ref{qe1})-(\ref{qe2}), we see that the function~ $g(t)=\sup_{s\in \lbrack
\tau _{N}\vee t,T]}\xi _{3}(s)\left\Vert p_{n}\right\Vert _{W}^{2}$ \
satisfies the Gronwall type inequality
\begin{equation*}
g(t)\leq C\int_{t}^{T}g(s)ds+C\left( \mathbb{E}\left\Vert \nabla
h(T)\right\Vert _{W}^{2}+\mathbb{E}\int_{t}^{T}~\xi _{3}(s)\left\Vert \nabla
_{y}L\left( s,U,Y\right) \right\Vert _{2}^{2}\,ds\right) ,
\end{equation*}%
that implies
\begin{equation*}
\frac{1}{2}\mathbb{E}\,\sup_{s\in \lbrack \tau _{N}\vee t,T]}\xi
_{3}(s)\left\Vert p_{n}(s)\right\Vert _{W}^{2}\leq C(\mathbb{E}\,\left\Vert
p_{n}(T)\right\Vert _{W}^{2}+\mathbb{E}\,\int_{t}^{T}~\xi _{3}(s)\left\Vert
\nabla _{y}L\left( s,U,Y\right) \right\Vert _{2}^{2}\,ds).
\end{equation*}%
Hence, combining this result with (\ref{qe1}), we derive
\begin{eqnarray*}
\frac{1}{2}\mathbb{E}\,\sup_{s\in \lbrack \tau _{N}\vee t,T]}\xi
_{3}(s)\left\Vert p_{n}(s)\right\Vert _{W}^{2} &+&\int_{\tau _{N}\vee
t}^{T}\xi _{3}(s)\left( \frac{\nu }{\alpha }\left\Vert \mathbb{P}\sigma
(p_{n})\right\Vert _{2}^{2}+4\nu \left\Vert Dp_{n}\right\Vert _{2}^{2}+\frac{%
1}{2}\left\Vert q_{n}\right\Vert _{W}^{2}\right) \,ds \\
&\leq &C(\mathbb{E}\,\left\Vert p_{n}(T)\right\Vert _{W}^{2}+\mathbb{E}%
\,\int_{t}^{T}~\xi _{3}(s)\left\Vert \nabla _{y}L\left( s,U,Y\right)
\right\Vert _{2}^{2}\,ds).
\end{eqnarray*}%
Passing to the limit $N\rightarrow \infty $ in this last inequality, we
deduce the claimed result.

\hfill $\hfill \hfill \blacksquare $

\bigskip

\bigskip

As a direct consequence of Proposition \ref{est_pq} we obtain the following
existence result.

\begin{theorem}
\label{the_23} There exists a unique solution
\begin{equation*}
(p,q)\in L^{\infty }(0,T;W)\times L^{2}(0,T;W)\qquad P-\text{a.e. in }\Omega
\end{equation*}%
of the backward stochastic system (\ref{p_n}). Moreover there exists a
positive constant $C_{3}$, such that the pair $(p,q)$ satisfies the
following estimate for a.e. $t\in (0,T)$
\begin{eqnarray*}
\frac{1}{2}\mathbb{E}\,\sup_{s\in \lbrack t,T]}\xi _{3}(s)\left\Vert
p(s)\right\Vert _{W}^{2} &+\mathbb{E}&\int_{t}^{T}\xi _{3}(\frac{\nu }{%
\alpha }\left\Vert \mathbb{P}\sigma (p)\right\Vert _{2}^{2}+4\nu \left\Vert
Dp\right\Vert _{2}^{2}+\frac{1}{2}\left\Vert q\right\Vert _{W}^{2})\,ds \\
&\leq &C(\mathbb{E}\,\left\Vert \nabla h\left( Y(T)\right) \right\Vert
_{W}^{2}+\mathbb{E}\,\int_{t}^{T}~\xi _{3}\left\Vert \nabla _{y}L\left(
s,U,Y\right) \right\Vert _{2}^{2}\,ds)
\end{eqnarray*}%
for the function%
\begin{equation*}
\xi _{3}(t)=e^{-C_{3}\int_{t}^{T}\left\Vert Y\right\Vert _{\widetilde{W}}ds}.
\end{equation*}
\end{theorem}

\textbf{Proof.} \vspace{2mm} By (\ref{lp1}) we have%
\begin{equation*}
\sup_{t\in \lbrack 0,T]}\left\Vert Y(t)\right\Vert _{\widetilde{W}}^{2}\leq
C(\omega )\quad \text{for all }\omega \in \Omega \backslash A,\quad \text{%
where }P(A)=0,\text{ }
\end{equation*}%
that is, there exists a positive constant $K(\omega ),$ dependent only on $%
\omega \in \Omega \backslash A,$ satisfying%
\begin{equation}
0<K(\omega )\leq \xi _{3}(t)\leq 1\quad \text{for all }\omega \in \Omega
\backslash A,\quad t\in \lbrack 0,T].  \label{kk}
\end{equation}

By a priori estimates of Proposition \ref{est_pq} we have%
\begin{equation*}
\mathbb{E}\sup_{t\in \lbrack 0,T]}\xi _{3}(t)\left\Vert p_{n}(t)\right\Vert
_{W}^{2}\leq C,\qquad \mathbb{E}\int_{0}^{T}\xi _{3}(t)\left\Vert
q_{n}(t)\right\Vert _{W}^{2}\ dt\leq C,
\end{equation*}%
hence using (\ref{kk}) there exists a subsequence of the pairs $\left(
p_{n},q_{n}\right) $, such that%
\begin{eqnarray}
S_{n} &=&\sqrt{\xi _{3}}p_{n}\rightharpoonup S=\sqrt{\xi _{3}}p\qquad
\mbox{ *-weakly in
}\ L^{2}(\Omega ,L^{\infty }(0,T;W)),  \notag \\
Q_{n} &=&\sqrt{\xi _{3}}q_{n}\rightharpoonup Q=\sqrt{\xi _{3}}q\qquad
\mbox{  weakly in
}\ L^{2}(\Omega ,L^{2}(0,T;W)).  \label{lim}
\end{eqnarray}

Moreover the pair $\left( p_{n},q_{n}\right) $ solves%
\begin{align}
\left( \sigma \left( \nabla h\left( Y(T)\right) \right) ,\phi \right)
-\left( \sigma \left( p_{n}(t)\right) ,\phi \right) & =\int_{t}^{T}\bigl\{%
2\nu \left( Dp_{n}(s),D\phi \right)  \notag \\
& -b\left( \phi ,p_{n}(s),\sigma (Y(s))\right) +b\left( p_{n}(s),\phi
,\sigma (Y(s))\right) \vspace{2mm}  \notag \\
& +b\left( p_{n}(s),Y(s),\sigma (\phi (s))\right) -b\left(
Y(s),p_{n}(s),\sigma (\phi (s))\right) \vspace{2mm}  \notag \\
& +\int_{t}^{T}\left( \nabla _{y}L\left( s,U,Y\right) +\nabla
_{y}G(s,Y)^{T}q_{n},\phi \right) \bigr\}\,ds\vspace{2mm}  \notag \\
& +\int_{t}^{T}\left( \sigma \left( q_{n}(s)\right) ,\phi \right)
\,dW_{s}\qquad \text{for each \ }\phi \in V,  \label{ap_n}
\end{align}%
therefore we can apply similar arguments as in the proof of Theorem \ref%
{the_2} and demonstrate that the pair $(p,q)$ satisfies equation (\ref%
{equation_etat1}).$\hfill \blacksquare $\newline

\bigskip

\subsection{Exponential integrability for the solution of equation (\protect
\ref{equation_etat})}


\label{102}

\bigskip

As in the previous sections we consider that the data $U,$ $Y_{0}$ satisfy
assumptions (\ref{au}) and denote the unique solution of the stochastic
differential equation\ (\ref{var_form_state}) by
\begin{equation*}
Y\in L^{p}(\Omega ,L^{\infty }(0,T;V))\cap L^{2}(\Omega ,L^{\infty }(0,T;%
\widetilde{W})),
\end{equation*}%
satisfying estimates (\ref{lp1}). In what follows, we assume an additional
assumptions on the data in order to improve the integrability of the
stochastic process $Y$. We will show that under these additional
assumptions, the stochastic process $Y$ \ is exponential integrable.

Let $U$ be a distributed mechanical force belonging to the admissible set $%
\mathcal{U}_{ad}^{b}$, which is defined as the set of all adapted stochastic
processes
\begin{equation}
U=U(\omega ,t)\in L^{\infty }(\Omega ,L^{2}(0,T;V))\cap L^{p}(\Omega
,L^{p}(0,T;H^{1}(\mathcal{O})))\quad \text{for some }4\leq p<\infty ,
\label{bound0}
\end{equation}%
uniformly bounded in $L^{2}(0,T;V)$, this means that there exists a positive
constant $M$, independent of $\omega$, such that
\begin{equation}
\int_{0}^{T}\left\Vert U(\omega ,t)\right\Vert _{V}^{2}\,dt\leq M\qquad
\text{for }P-\text{a.e. }\omega \in \Omega .  \label{bound}
\end{equation}

We also introduce additional hypothesis on the diffusion operator $G$,
namely $G$ is bounded by a positive constant $L$ in the space $V$\
\begin{equation}
||G\left( t,y\right) ||_{V}^{2}\leq L\quad \quad \text{for all \ }t\in
\lbrack 0,T],\quad y\in V.  \label{g}
\end{equation}%
Let us set $C_{\max }=\max \{C_{1},C_{2},C_{3}\}$, where $C_{i}$, $i=1,2,3$,
are the constants defined in Propositions \ref{ex_uniq_lin} and \ref{est_pq}%
, and introduce the following two conditions:
\vspace{2mm}\\
\textit{the first condition} is given by \textbf{\ }
\begin{equation}
A=\frac{1}{2\theta _{1}}\geq C_{\max }  \label{g1}
\end{equation}%
and \textit{the second one} reads as
\begin{equation}
\text{the domain\textit{\ }}\mathcal{O}\text{\textit{\ is non axisymmetric }%
and}\quad B=\frac{\gamma _{2}^{2}}{2\theta _{2}}\geq C_{\max }.  \label{g2}
\end{equation}%
Here the constants \ $\theta _{1},$ $\theta _{2},$ $\gamma _{1},$ $\bar{%
\gamma}_{1}$, $\gamma _{2}$ and $\ \tau $ are defined by%
\begin{eqnarray}
\theta _{1} &=&4L\tau ^{2}(1+\gamma _{1}^{2}),\qquad \theta
_{2}=2Le^{\varepsilon T}(1+2\bar{\gamma}_{1}^{2}),  \notag \\
\gamma _{1} &=&\frac{C_{\ast }\nu }{\alpha }\tau ,\qquad \bar{\gamma}_{1}=%
\frac{K_{\ast }}{4\alpha }\tau ,\qquad \gamma _{2}=\frac{C_{\ast \ast }\nu }{%
\alpha }\quad \text{and}\quad \tau =Te^{\varepsilon T},  \label{con}
\end{eqnarray}%
with the constants $C_{\ast },$ $K_{\ast }$\ and $C_{\ast \ast }$ \
introduced in (\ref{korn}), (\ref{kornn}) and (\ref{CCC}).

Now, we consider the strong solution $Y$ of the state system (\ref%
{equation_etat}) which exists by Theorem \ref{the_1}, and our purpose is to
show the exponential integrability of the state process $Y$. Let us mention
that the two main arguments to show this result rely on the structure of
equation $(\ref{equation_etat})_{1}$ for $Y$, and on the martingale property
of the exponential process that appears in the right hand side of
inequalities \eqref{caw3}-\eqref{caw4}.

\begin{proposition}
\label{pop} Assume that $U\in \mathcal{U}_{ad}^{b},$ satisfies (\ref{bound0}%
)-(\ref{bound}) and the initial condition
\begin{equation}
Y_{0}\in L^{p}(\Omega ,V)\cap L^{\infty }(\Omega ,\widetilde{W})\bigskip .
\label{inbound}
\end{equation}%
Also we admit that condition (\ref{g}) and one of two conditions (\ref{g1})
or (\ref{g2}) hold, then there exists a positive constant $C$, such that the
following estimate is valid
\begin{equation}
\mathbb{E}\,\exp \biggl\{C_{\max }\int_{0}^{t}\left\Vert Y(s)\right\Vert _{%
\widetilde{W}}^{2}\,ds\biggr\}\leq C\qquad \text{for a.e. }t\in (0,T).
\label{yy}
\end{equation}
\end{proposition}

\textbf{Proof.} \ \ \textit{1st step: \ Estimates \ for }$Y$\textit{\ in the
space }$V.$

Let $\ \widetilde{G}$ be the solution of (\ref{GS_NS}) for $f=G(t,Y)$. \
Using the fact that $\widetilde{G}$ solves the elliptic type problem (\ref%
{GS_NS}) for $f=G(t,Y),$ (\ref{upo}) and assumption (\ref{g}), we have
\begin{equation*}
||\widetilde{G}||_{V}^{2}\leq ||G\left( t,Y\right) ||_{2}^{2}\leq L.
\end{equation*}

Applying the operator $\left( I-\alpha \mathbb{P}\Delta \right) ^{-1}$ to
equation (\ref{equation_etat})$_{1},$ we deduce a stochastic differential
equation for $Y$, then with the help of \ It\^{o}'s formula, as it was done
in the article \cite{CC16}, we obtain
\begin{equation*}
d\left\Vert Y\right\Vert _{V}^{2}=2\left( -2\nu \left\Vert DY\right\Vert
_{2}^{2}+\left( U,Y\right) \right) \,dt+2\left( G(t,Y),Y\right) \,dW_{t}+||%
\widetilde{G}||_{V}^{2}\,dt.
\end{equation*}%
Using Young%
\'{}%
s inequality $2ab\leq \frac{a^{2}}{\varepsilon }+\varepsilon b^{2}$ and
integrating over the time interval $(0,t)$, we obtain
\begin{equation}
\left\Vert Y(t)\right\Vert _{V}^{2}+4\nu \int_{0}^{t}\left\Vert
DY\right\Vert _{2}^{2}ds\leq \varepsilon \int_{0}^{t}\left\Vert Y\right\Vert
_{2}^{2}\,ds+(C_{\varepsilon }+g_{1}(t)),  \label{qq}
\end{equation}%
where
\begin{eqnarray*}
C_{\varepsilon } &=&\left\Vert Y(0)\right\Vert _{V}^{2}+\frac{1}{\varepsilon
}\int_{0}^{T}\left\Vert U\right\Vert _{2}^{2}\,dt+LT, \\
g_{1}(t) &=&\int_{0}^{t}f_{1}(s)\,dW_{s},\qquad f_{1}(s)=2\left(
G(s,Y),Y\right) .
\end{eqnarray*}%
Taking $z(t)=\int_{0}^{t}\left\Vert Y\right\Vert _{V}^{2}\,ds,$ expression (%
\ref{qq}) gives the differential inequality $z%
{\acute{}}%
\leq \varepsilon z+(C_{\varepsilon }+g_{1}(t)),$ which can be integrated by
Gronwall's lemma. Hence
\begin{eqnarray}
\int_{0}^{t}\left\Vert Y\right\Vert _{V}^{2}\,ds &=&z(t)\leq e^{\varepsilon
t}\int_{0}^{t}e^{-\varepsilon s}(C_{\varepsilon }+g_{1}(s))\ ds  \label{w1}
\\
&=&C_{\varepsilon }\left( \frac{e^{\varepsilon t}-1}{\varepsilon }\right)
+e^{\varepsilon t}\int_{0}^{t}e^{-\varepsilon s}g_{1}(s)\ ds.  \notag
\end{eqnarray}%
Since the inequality $|\frac{e^{x}-1}{x}|\leq e^{|x|},$ we have $|\frac{%
e^{\varepsilon t}-1}{\varepsilon }|\leq te^{\varepsilon t}\leq
Te^{\varepsilon T}$. \ Applying Fubini's theorem, we derive
\begin{equation}
e^{\varepsilon t}\int_{0}^{t}e^{-\varepsilon s}g_{1}(s)\
ds=\int_{0}^{t}\left( \frac{e^{\varepsilon (t-s)}-1}{\varepsilon }\right) \
f_{1}(s)\,dW_{s}.  \label{fub}
\end{equation}%
Hence
\begin{equation}
\int_{0}^{t}\left\Vert Y\right\Vert _{V}^{2}\,ds\leq C_{\varepsilon
}Te^{\varepsilon T}+\int_{0}^{t}\left( \frac{e^{\varepsilon (t-s)}-1}{%
\varepsilon }\right) f_{1}(s)\ dW_{s}.  \label{211}
\end{equation}

\bigskip Now, using (\ref{w1}) to estimate the right hand side of (\ref{qq}%
), we deduce
\begin{align*}
& \left\Vert Y(t)\right\Vert _{V}^{2}+4\nu \int_{0}^{t}\left\Vert
DY\right\Vert _{2}^{2}\,ds\leq \varepsilon \int_{0}^{t}\left\Vert
Y\right\Vert _{2}^{2}\,ds+(C_{\varepsilon }+g_{1}(t)) \\
& \leq \varepsilon \left[ e^{\varepsilon t}\int_{0}^{t}e^{-\varepsilon
s}(C_{\varepsilon }+g_{1}(s))\ ds\right] +(C_{\varepsilon }+g_{1}(t)) \\
& =C_{\varepsilon }e^{\varepsilon t}+\varepsilon e^{\varepsilon
t}\int_{0}^{t}e^{-\varepsilon s}g_{1}(s)\ ds+g_{1}(t) \\
& =C_{\varepsilon }e^{\varepsilon t}+\int_{0}^{t}e^{\varepsilon (t-s)}\
f_{1}(s)\,dW_{s}
\end{align*}%
by \eqref{fub}. Therefore%
\begin{equation}
\left\Vert Y(t)\right\Vert _{V}^{2}+4\nu \int_{0}^{t}\left\Vert
DY\right\Vert _{2}^{2}\,ds\leq C_{\varepsilon }e^{\varepsilon
t}+\int_{0}^{t}e^{\varepsilon (t-s)}\ f_{1}(s)\,dW_{s}.  \label{21}
\end{equation}%
\bigskip

\bigskip

\textit{2nd step: \ Estimates \ for }$\mathrm{curl}\,Y$\textit{\ in the
space }$L^{2}.$ Using (\ref{korn}) and (\ref{211}),\ we obtain
\begin{equation}
\frac{\nu }{\alpha }\int_{0}^{t}\left\Vert \mathrm{curl}\,Y\right\Vert
_{2}^{2}\ ds\leq \frac{C_{\ast }\nu }{\alpha }\left\{ C_{\varepsilon
}Te^{\varepsilon T}+\int_{0}^{t}\left( \frac{e^{\varepsilon (t-s)}-1}{%
\varepsilon }\right) f_{1}(s)\,dW_{s}\right\} .  \label{curlest}
\end{equation}%
Under the assumption that the domain $\mathcal{O}$ is not axisymmetric, we
may apply inequality (\ref{kornn}), in order to deduce as alternative
estimate for $\mathrm{curl}\,Y$, namely, we have
\begin{equation*}
\frac{\nu }{\alpha }\int_{0}^{t}\left\Vert \mathrm{curl}\,Y\right\Vert
_{2}^{2}\ ds\leq \frac{\nu }{\alpha }\int_{0}^{t}\left\Vert Y\right\Vert
_{H^{1}}^{2}\ ds\leq \frac{K_{\ast }\nu }{\alpha }\int_{0}^{t}\left\Vert
DY\right\Vert _{2}^{2}\ ds
\end{equation*}%
and considering (\ref{21}), we obtain
\begin{equation}
\frac{\nu }{\alpha }\int_{0}^{t}\left\Vert \mathrm{curl}\,Y\right\Vert
_{2}^{2}\ ds\leq \frac{K_{\ast }}{4\alpha }\left\{ C_{\varepsilon
}e^{\varepsilon T}+\int_{0}^{t}e^{\varepsilon (t-s)}f_{1}(s)\,dW_{s}\right\}
.  \label{21a}
\end{equation}

\bigskip

\textit{3d step: \ Two estimates\ for }$\mathrm{curl}\,\sigma \left(
Y\right) $\textit{\ in the space }$L^{2}.$ \ \ Applying the operator $%
\mathrm{curl}$ to equation (\ref{equation_etat})$_{1}$, we deduce
\begin{align*}
d\left( \mathrm{curl}\,\sigma \left( Y\right) \right) & =\frac{\nu }{\alpha }%
\left( -\,\mathrm{curl}\,\sigma (Y)+\mathrm{curl}\,Y\right) \,dt+\mathrm{curl%
}\,U\,dt \\
& -(Y\cdot \nabla )\mathrm{curl}\,\sigma \left( Y\right) \,dt+\mathrm{curl}%
\,G(t,Y)\,dW_{t}.
\end{align*}%
The \ It\^{o}'s formula yields
\begin{align*}
& d\Vert \mathrm{curl}\,\sigma \left( Y\right) \Vert _{2}^{2}=\frac{2\nu }{%
\alpha }\left( -\,\mathrm{curl}\,\sigma (Y)+\mathrm{curl}\,Y,\mathrm{curl}%
\,\sigma (Y)\right) \,dt+2\left( \mathrm{curl}\,U,\mathrm{curl}\,\sigma
(Y)\right) \,dt \\
& -2\left( (Y\cdot \nabla )\mathrm{curl}\,\sigma \left( Y\right) ,\mathrm{%
curl}\,\sigma (Y)\right) dt+2(\mathrm{curl}\,G(t,Y),\mathrm{curl}\,\sigma
(Y))\,dW_{t}+||\mathrm{curl}\,G(t,Y)||_{2}^{2}\,dt.
\end{align*}%
Since
\begin{equation*}
\left( (Y\cdot \nabla )\mathrm{curl}\,\sigma \left( Y\right) ,\mathrm{curl}%
\,\sigma (Y)\right) =0,
\end{equation*}%
we derive
\begin{align}
d\Vert \mathrm{curl}\,\sigma \left( Y\right) \Vert _{2}^{2}& =\frac{2\nu }{%
\alpha }\left( -\,\mathrm{curl}\,\sigma (Y)+\mathrm{curl}\,Y,\mathrm{curl}%
\,\sigma (Y)\right) \,dt+2\left( \mathrm{curl}\,U,\mathrm{curl}\,\sigma
(Y)\right) \,dt  \notag \\
& +2(\mathrm{curl}\,G(t,Y),\mathrm{curl}\,\sigma (Y))\,dW_{t}+||\mathrm{curl}%
\,G(t,Y)||_{2}^{2}\,dt.  \label{21:29}
\end{align}%
From assumption (\ref{g}) we have
\begin{equation}
||\mathrm{curl}\,G(t,Y)||_{2}^{2}\leq C||G\left( t,Y\right) ||_{V}^{2}\leq
CL,  \notag
\end{equation}%
in addition, the relation $2(-y+x)y\leq- y^{2}+x^{2}$, $\forall x,y\in
\mathbb{R}$, gives
\begin{equation*}
2\left( -\,\mathrm{curl}\,\sigma (Y)+\mathrm{curl}\,Y,\mathrm{curl}\,\sigma
(Y)\right) \leq -\left\Vert \mathrm{curl}\,\sigma \left( Y\right)
\right\Vert _{2}^{2}+\left\Vert \mathrm{curl}\,Y\right\Vert _{2}^{2}.
\end{equation*}%
Hence, introducing these estimates in \eqref{21:29} and integrating over $%
[0,t]$, we obtain
\begin{align*}
\left\Vert \mathrm{curl}\,\sigma \left( Y\right) \right\Vert _{2}^{2}& +%
\frac{\nu }{\alpha }\int_{0}^{t}\left\Vert \mathrm{curl}\,\sigma \left(
Y\right) \right\Vert _{2}^{2}\,ds\leq \left\Vert \mathrm{curl}\,\sigma
\left( Y(0)\right) \right\Vert _{2}^{2}+\frac{\nu }{\alpha }%
\int_{0}^{t}\left\Vert \mathrm{curl}\,Y\right\Vert _{2}^{2}\ ds+CLT \\
& +2\int_{0}^{t}\left( \mathrm{curl}\,U,\mathrm{curl}\,\sigma (Y)\right)
\,ds+2\int_{0}^{t}(\mathrm{curl}\,G(s,Y),\mathrm{curl}\,\sigma (Y))\,dW_{s}.
\end{align*}%
Therefore, if we use Young%
\'{}%
s inequality $2ab\leq \frac{a^{2}}{\varepsilon }+\varepsilon b^{2}$ and (\ref%
{curlest}), we deduce
\begin{equation}
\left\Vert \mathrm{curl}\,\sigma \left( Y\right) \right\Vert _{2}^{2}+\frac{%
\nu }{\alpha }\int_{0}^{t}\left\Vert \mathrm{curl}\,\sigma \left( Y\right)
\right\Vert _{2}^{2}\,dt\leq \varepsilon \int_{0}^{t}\left\Vert \mathrm{curl}%
\,\sigma (Y)\right\Vert ^{2}\ ds+(\widetilde{C}_{\varepsilon }+g_{2}(t)),
\notag
\end{equation}%
where%
\begin{eqnarray*}
\widetilde{C}_{\varepsilon } &=&\frac{1}{\varepsilon }\int_{0}^{T}\left\Vert
\mathrm{curl}\,U\right\Vert _{2}^{2}\,ds+\left\Vert \mathrm{curl}\,\sigma
\left( Y(0)\right) \right\Vert _{2}^{2}+CLT+\frac{C_{\ast }\nu }{\alpha }%
Te^{\varepsilon T}, \\
g_{2}(t) &=&\int_{0}^{t}f_{2}(s)\,dW_{s}, \\
f_{2}(s) &=&2(\mathrm{curl}\,G(s,Y),\mathrm{curl}\,\sigma (Y))\,+\frac{%
C_{\ast }\nu }{\alpha }\left( \frac{e^{\varepsilon (t-s)}-1}{\varepsilon }%
\right) 2\left( G(s,Y),Y\right) .
\end{eqnarray*}%
On the other hand, due to (\ref{21a}) we also have
\begin{equation}
\left\Vert \mathrm{curl}\,\sigma \left( Y\right) \right\Vert _{2}^{2}+\frac{%
\nu }{\alpha }\int_{0}^{t}\left\Vert \mathrm{curl}\,\sigma \left( Y\right)
\right\Vert _{2}^{2}\,dt\leq \varepsilon \int_{0}^{t}\left\Vert \mathrm{curl}%
\,\sigma (Y)\right\Vert ^{2}\ ds+(\bar{C}_{\varepsilon }+\bar{g}_{2}(t)),
\notag
\end{equation}%
where%
\begin{eqnarray*}
\bar{C}_{\varepsilon } &=&\frac{1}{\varepsilon }\int_{0}^{T}\left\Vert
\mathrm{curl}\,U\right\Vert _{2}^{2}\,ds+\left\Vert \mathrm{curl}\,\sigma
\left( Y(0)\right) \right\Vert _{2}^{2}+CLT+\frac{K_{\ast }}{4\alpha }%
e^{\varepsilon T}, \\
\bar{g}_{2}(t) &=&\int_{0}^{t}\bar{f}_{2}(s)\,dW_{s}, \\
\bar{f}_{2}(s) &=&2(\mathrm{curl}\,G(s,Y),\mathrm{curl}\,\sigma (Y))\,+\frac{%
K_{\ast }}{2\alpha }e^{\varepsilon (t-s)}\left( G(s,Y),Y\right) .
\end{eqnarray*}

Now, we proceed as in the\textit{\ 1st step,} \ defining $\
~z(t)=\int_{0}^{t}\left\Vert \mathrm{curl}\,\sigma \left( Y\right)
\right\Vert _{2}^{2}\,ds$. In the former case, we get the inequality $~z%
{\acute{}}%
\leq \varepsilon z+(\widetilde{C}_{\varepsilon }+g_{2}(t)),~$ which may be
integrated to obtain
\begin{equation}
\int_{0}^{t}\left\Vert \mathrm{curl}\,\sigma \left( Y\right) \right\Vert
_{2}^{2}\,ds=z(t)\leq \widetilde{C}_{\varepsilon }Te^{\varepsilon
T}+\int_{0}^{t}\ \left( \frac{e^{\varepsilon (t-s)}-1}{\varepsilon }\right)
f_{2}(s)\,dW_{s}.  \label{f}
\end{equation}%
Applying the same reasoning, the second situation gives
\begin{equation}
\left\Vert \mathrm{curl}\,\sigma \left( Y\right) \right\Vert _{2}^{2}+\frac{%
\nu }{\alpha }\int_{0}^{t}\left\Vert \mathrm{curl}\,\sigma \left( Y\right)
\right\Vert _{2}^{2}\,ds\leq \bar{C}_{\varepsilon }e^{\varepsilon
t}+\int_{0}^{t}e^{\varepsilon (t-s)}\bar{f}_{2}(s)\,dW_{s}.  \label{f2}
\end{equation}%
For details we refer to the deduction of inequalities (\ref{211}) and (\ref%
{21}).

\bigskip

\bigskip \textit{4th step: \ Two exponential estimates\ for }$Y.$ \ Taking
the sum of inequalities (\ref{211}) and (\ref{f}) we get
\begin{equation}
\int_{0}^{t}\left\Vert Y\right\Vert _{\widetilde{W}}^{2}\,ds\leq
(C_{\varepsilon }+\widetilde{C}_{\varepsilon })Te^{\varepsilon
T}+\int_{0}^{t}\left( \frac{e^{\varepsilon (t-s)}-1}{\varepsilon }\right) \ %
\left[ f_{1}(s)+f_{2}(s)\right] \,dW_{s}  \label{caw1}
\end{equation}%
and the sum of inequalities (\ref{21}) and (\ref{f2}) yields
\begin{equation*}
\frac{\nu }{\alpha }\int_{0}^{t}\left( 2\alpha \left\Vert DY\right\Vert
_{2}^{2}\,+\left\Vert \mathrm{curl}\,\sigma \left( Y\right) \right\Vert
_{2}^{2}\right) \ ds\leq (\frac{1}{2}C_{\varepsilon }+\bar{C}_{\varepsilon
})e^{\varepsilon t}+\int_{0}^{t}e^{\varepsilon (t-s)}\ \left[ \frac{1}{2}%
f_{1}(s)+\bar{f}_{2}(s)\right] \,dW_{s},
\end{equation*}%
then applying (\ref{CCC}), we derive
\begin{equation}
\frac{C_{\ast \ast }\nu }{\alpha }\int_{0}^{t}\left\Vert \,Y\right\Vert _{%
\widetilde{W}}^{2}\ ds\leq \frac{1}{2}(C_{\varepsilon }+\bar{C}_{\varepsilon
})e^{\varepsilon t}+\int_{0}^{t}e^{\varepsilon (t-s)}\ \left[ \frac{1}{2}%
f_{1}(s)+\bar{f}_{2}(s)\right] \,dW_{s}.  \label{caw2}
\end{equation}%
\bigskip

Let us denote by
\begin{equation*}
Z_{1}(s)=\left( \frac{e^{\varepsilon (t-s)}-1}{\varepsilon }\right) \left[
f_{1}(s)+f_{2}(s)\right] ,\qquad {Z}_{2}(s)=e^{\varepsilon (t-s)}\left[
\frac{1}{2}f_{1}(s)+\bar{f}_{2}(s)\right]
\end{equation*}%
and the constants
\begin{equation*}
H_{1}=\exp \left( (C_{\varepsilon }+\widetilde{C}_{\varepsilon
})Te^{\varepsilon T}\right) ,\qquad H_{2}=\exp \left( \frac{1}{2}%
C_{\varepsilon }+\bar{C}_{\varepsilon })e^{\varepsilon T}\right) \qquad
\text{for any fixed \ }\varepsilon >0.
\end{equation*}

Therefore if we multiply (\ref{caw1}) and (\ref{caw2}) by $\lambda ,$
applying the exponential function and the expectation, we deduce
\begin{align}
& \mathbb{E}\,\exp \biggl\{\left( \lambda -\frac{1}{2}\lambda ^{2}\theta
_{1}\right) \int_{0}^{t}\left\Vert Y\right\Vert _{\widetilde{W}}^{2}\,\ ds%
\biggr\}  \notag \\
& \leq \mathbb{E}\,\exp \biggl\{\lambda \int_{0}^{t}\left\Vert Y\right\Vert
_{\widetilde{W}}^{2}\,\ ds-\frac{\lambda ^{2}}{2}\int_{0}^{t}Z_{1}(s)^{2}\,%
\,ds\biggr\}  \notag \\
& \leq H_{1}\mathbb{E}\,\exp \biggl\{\int_{0}^{t}\lambda Z_{1}(s)\,dW_{s}-%
\frac{1}{2}\int_{0}^{t}\left( \lambda Z_{1}(s)\right) ^{2}\,\,ds\biggr\}
\label{caw3}
\end{align}%
and
\begin{align}
& \mathbb{E}\,\exp \biggl\{\left( \gamma _{2}\lambda -\frac{1}{2}\lambda
^{2}\theta _{2}\right) \int_{0}^{t}\left\Vert Y\right\Vert _{\widetilde{W}%
}^{2}\,ds\biggr\}  \notag \\
& \leq \mathbb{E}\,\exp \biggl\{\lambda \int_{0}^{t}\left\Vert Y\right\Vert
_{\widetilde{W}}^{2}\,ds-\frac{\lambda ^{2}}{2}\int_{0}^{t}Z_{2}(s)^{2}\,\,ds%
\biggr\}  \notag \\
& \leq H_{2}\mathbb{E}\exp \biggl\{\int_{0}^{t}\lambda Z_{2}(s)\,dW_{s}-%
\frac{1}{2}\int_{0}^{t}\left( \lambda Z_{2}(s)\right) ^{2}\,\,ds\biggr\},
\label{caw4}
\end{align}%
since $G(t,y),$ $U$ \ and the initial condition $Y_{0}$\ satisfy (\ref{g}), (%
\ref{bound}) and (\ref{inbound}), respectively. Here the constants \ $\theta
_{1},$ $\theta _{2},$ $\gamma _{1},$ $\gamma _{2}$ and $\ \tau $ are defined
by (\ref{con}).\

\bigskip

Taking into account that the expectations of the right hand sides in (\ref%
{caw3}) and in (\ref{caw4}) \ are equal to $1,$ we derive
\begin{equation*}
\mathbb{E}\,\exp \biggl\{\left( \lambda -\frac{\lambda ^{2}}{2}\theta
_{1}\right) \int_{0}^{t}\left\Vert Y\right\Vert _{\widetilde{W}}^{2}\,ds%
\biggr\}\leq H_{1}
\end{equation*}%
and
\begin{equation*}
\mathbb{E}\,\exp \biggl\{\left( \gamma _{2}\lambda -\frac{\lambda ^{2}}{2}%
\theta _{2}\right) \int_{0}^{t}\left\Vert Y\right\Vert _{\widetilde{W}%
}^{2}\,ds\biggr\}\leq H_{2}.
\end{equation*}%
Hence, we derive estimate (\ref{yy}), since one of conditions (\ref{g1}), (%
\ref{g2}) holds and
\begin{equation*}
A=\frac{1}{2\theta _{1}}=\max_{\lambda }\left( \lambda -\frac{\lambda ^{2}}{2%
}\theta _{1}\right) \qquad \text{and}\qquad B=\frac{\gamma _{2}^{2}}{2\theta
_{2}}=\max_{\lambda }\left( \gamma _{2}\lambda -\frac{\lambda ^{2}}{2}\theta
_{2}\right) .
\end{equation*}

$\bigskip \hfill \hfill \blacksquare $\newline

\bigskip

\bigskip

\subsection{Duality property}


\label{103}

\bigskip

In the next proposition we prove that the functions $Z_{n}$, constructed in
Proposition \ref{ex_uniq_lin}, and $(p_{n},q_{n})$, constructed in
Proposition \ref{est_pq}, satisfy the so-called \textit{duality} property.

\begin{proposition}
The solution $Z_{n}$ of \ \ system (\ref{pn}) and the solution $%
(p_{n},q_{n}) $ of system (\ref{est_pq}) satisfy the\textit{\ duality
equality}
\begin{equation}
\mathbb{E}\left( \sigma (Z_{n}(T)),\nabla h(Y_{T})\right) +\mathbb{E}%
\int_{0}^{T}\left( \nabla _{y}L\left( t,U,Y\right) ,Z_{n}\right) \ dt=%
\mathbb{E}\int_{0}^{T}\left( \Psi ,p_{n}\right) \,dt  \label{dual}
\end{equation}%
for any $\ \Psi \in L^{2}\left( \Omega \times (0,T)\times \mathcal{O}\right)
.$
\end{proposition}

\textbf{Proof. \ } Setting $\phi =h_{i}$ in equation (\ref{pn}), we obtain%
\begin{equation}
d\left( p_{n},h_{i}\right) _{V}=\left( \widetilde{f}(p_{n}),h_{i}\right) \,dt%
\vspace{2mm}-\left( R(q_{n}),h_{i}\right) \,dt+\left( \sigma \left(
q_{n}\right) ,h_{i}\right) \,dW_{t}  \label{ppp}
\end{equation}%
for%
\begin{eqnarray*}
\widetilde{f}(p_{n}) &=&-\nu \Delta p_{n}-\mathrm{curl}\,\sigma (Y)\times
p_{n}+\mathrm{curl}\left( \sigma \left( Y\times p_{n}\right) \right) , \\
R(q_{n}) &=&\nabla _{y}L\left( t,U,Y\right) +\nabla _{y}G(t,Y)^{T}q_{n}.
\end{eqnarray*}%
Let us remember that the solution
\begin{equation*}
Z_{n}(t)=\sum_{i=1}^{n}\zeta _{i}(t)h_{i}
\end{equation*}%
solves the stochastic differential equation (\ref{equation_etat_temps_n}).
Therefore multiplying (\ref{ppp}) by $\zeta _{i}(t)$ and summing over $%
i=1,...,n,$ we obtain
\begin{equation*}
\left( dp_{n},Z_{n}\right) _{V}=\left( \widetilde{f}(p_{n}),Z_{n}\right) \,dt%
\vspace{2mm}-\left( R(q_{n}),Z_{n}\right) \,dt+\left( q_{n},Z_{n}\right)
_{V}\,dW_{t}.
\end{equation*}%
Using property (\ref{bb1}) we get the following identity
\begin{align}
\left( dp_{n},Z_{n}\right) _{V}& =\{2\nu \left( Dp_{n}(t),DZ_{n}\right) \ dx
\notag \\
& -b\left( Z_{n},p_{n},\sigma (Y)\right) +b\left( p_{n},Z_{n},\sigma
(Y)\right)  \notag \\
& +b\left( p_{n},Y,\sigma (Z_{n})\right) -b\left( Y,p_{n},\sigma
(Z_{n})\right)  \notag \\
& -\left( \nabla _{y}L\left( t,U,Y\right) +\nabla
_{y}G(t,Y)^{T}q_{n},Z_{n}\right) \}\ dt  \notag \\
& +\left( q_{n},\sigma \left( Z_{n}\right) \right) \,dW_{t}.  \label{qq1}
\end{align}

By another hand \ since\ $Z_{n}$ is the solution of (\ref%
{equation_etat_temps_n}), then setting $\phi =h_{i}$ in this equation, we
get
\begin{equation*}
d\left( Z_{n},h_{i}\right) _{V}=\left( f(Z_{n}),h_{i}\right) \,dt\vspace{2mm}%
+\left( \nabla _{y}G(t,Y)Z_{n},h_{i}\right) \,dW_{t},
\end{equation*}%
where $f(Z_{n})=\nu \Delta Z_{n}-\mathrm{curl}\,\sigma (Z_{n})\times Y-%
\mathrm{curl}\,\sigma (Y)\times Z_{n}+\Psi .$ Recalling that $p_{n}$ is
defined by (\ref{pnqn}), if we multiply the last equality by $\mathfrak{p}%
_{i}(t)$ and sum over $i=1,...,n,$ we deduce%
\begin{equation*}
\left( dZ_{n},p_{n}\right) _{V}=\left( f(Z_{n}),p_{n}\right) \,dt\vspace{2mm}%
+\left( \nabla _{y}G(t,Y)Z_{n},p_{n}\right) \,dW_{t}.
\end{equation*}%
Applying identities (\ref{bb1}) and (\ref{bb2}), we derive%
\begin{eqnarray}
\left( dZ_{n},p_{n}\right) _{V} &=&\{-2\nu \left( Dp_{n}(t),DZ_{n}\right)
\notag \\
&&-b\left( p_{n},Y,\sigma (Z_{n})\right) +b\left( Y,p_{n},\sigma
(Z_{n})\right)  \notag \\
&&-b\left( p_{n},Z_{n},\sigma (Y)\right) +b\left( Z_{n},p_{n},\sigma
(Y)\right)  \notag \\
&&+(\Psi ,p_{n})\}\,dt+\left( \nabla _{y}G(t,Y)Z_{n},p_{n}\right) \,dW_{t}.
\label{dz}
\end{eqnarray}

Hence Ito's formula gives
\begin{eqnarray*}
d\left( Z_{n},p_{n}\right) _{V} &=&\left( dZ_{n},p_{n}\right) _{V}+\left(
Z_{n},dp_{n}\right) _{V}+\left( dZ_{n},dp_{n}\right) _{V} \\
&=&\left[ -\left( \nabla _{y}L\left( t,U,Y\right) ,Z_{n}\right) +(\Psi
,p_{n})\right] \ dt \\
&&+\left[ \left( q_{n},\sigma \left( Z_{n}\right) \right) +\left( \nabla
_{y}G(t,Y)Z_{n},p_{n}\right) \,\right] dW_{t},
\end{eqnarray*}%
combining formulaes (\ref{qq1}) and (\ref{dz}). The integration of the last
identity over the time interval $[0,T]$ gives%
\begin{align*}
\left( Z_{n}(T),p_{n}(T)\right) _{V}& =\int_{0}^{T}\left( \Psi ,p_{n}\right)
\,dt-\left( \nabla _{y}L\left( t,U,Y\right) ,Z_{n}\right) \ dt \\
& +\int_{0}^{T}\left[ \left( q_{n},\sigma (Z_{n})\right) +(\nabla
_{y}G(t,Y)Z_{n},p_{n})\right] \,dW_{t}.
\end{align*}%
Finally taking the expectation, \ we conclude that the pair $(p_{n},q_{n})$
satisfies relation (\ref{dual}).$\bigskip \hfill \hfill \blacksquare $%
\newline

\bigskip

\bigskip

In the next theorem we prove that $Z,$ constructed in Theorem \ref{the_2},
and $(p,q)$, constructed in Theorem \ref{the_23}, satisfy the duality
condition.

\begin{theorem}
\label{adjoint} Under assumptions of Proposition \ref{pop}, the solution
\begin{equation*}
Z\in L^{2}(\Omega ,L^{\infty }(0,T;W))
\end{equation*}%
of system (\ref{5.3}) and the solution $(p,q)\in L^{2}(\Omega ,L^{\infty
}(0,T;W))\times L^{2}(\Omega ,L^{2}(0,T;W))$ \ of \ the backward stochastic
equation (\ref{p_n}) satisfy the following duality property
\begin{equation}
\mathbb{E}\int_{0}^{T}\left( \nabla _{y}L\left( t,U,Y\right) ,Z\right) \ dt+%
\mathbb{E}\left( \sigma (Z(T)),\nabla h(Y_{T})\right) =\mathbb{E}%
\int_{0}^{T}\left( \Psi ,p\right) \,dt  \label{ww}
\end{equation}%
for any $\ \Psi \in L^{2}\left( \Omega \times (0,T)\times \mathcal{O}\right)
.$
\end{theorem}

\bigskip

\textbf{Proof.} \vspace{2mm} Firstly let us show that for subsequences of $%
\left\{ Z_{n}\right\} $ and $\ \left\{ \left( p_{n},q_{n}\right) \right\} $
we have
\begin{eqnarray}
\mathbb{E}\int_{0}^{T}\left( \nabla _{y}L\left( t,U,Y\right) ,Z_{n}\right) \
dt &\rightarrow &\mathbb{E}\int_{0}^{T}\left( \nabla _{y}L\left(
t,U,Y\right) ,Z\right) \ dt,  \notag \\
\mathbb{E}\left( \sigma (Z_{n}(T)),\nabla h(Y_{T})\right) &\rightarrow &%
\mathbb{E}\left( \sigma (Z(T)),\nabla h(Y_{T})\right) ,  \notag \\
\mathbb{E}\int_{0}^{T}\left( \Psi ,p_{n}\right) \,dt &\rightarrow &\mathbb{E}%
\int_{0}^{T}\left( \Psi ,p\right) \,dt\qquad \text{as }n\rightarrow \infty
\label{zzz}
\end{eqnarray}%
\ with $Z$ and the pair$\ \left\{ \left( p,q\right) \right\} $\ defined in (%
\ref{lim1}), (\ref{lim}).

By assumptions (\ref{lip1}) we know that $~\nabla _{y}L\left( t,U,Y\right)
\in L^{2}(\Omega \times (0,T)\times \mathcal{O}).$ Due to the absolute
continuity of the Lebesgue integral, for any fixed $~\varepsilon >0$, \
there exists $~\lambda =\lambda (\varepsilon )>0~$ such that
\begin{equation*}
\mathbb{E}\int_{0}^{T}||\mathbb{\chi}\nabla _{y}L\left( t,U,Y\right)
||^{2}dt<\varepsilon
\end{equation*}
where $~\mathbb{\chi}~$ is the characteristic function of the set $%
~\{(t,x,\omega ):~|\nabla _{y}L\left( t,U,Y\right) |>\lambda \}.~$ Hence
\begin{eqnarray}
|\mathbb{E}\int_{0}^{T}\mathbb{~}|( \mathbb{\chi}\nabla _{y}L\left(
t,U,Y\right) ,Z_{n}-Z) dt| &\leq& \left( \mathbb{E}\int_{0}^{T}||\mathbb{\chi%
}\nabla _{y}L\left( t,U,Y\right) ||^{2}dt\right) ^{1/4} \times \left(
\mathbb{E}\int_{0}^{T}\xi _{1}^{-2}(t)dt\right) ^{1/4}  \notag \\
&\times & \left( \mathbb{E}\int_{0}^{T}\xi
_{1}(t)(||Z_{n}||^{2}+||Z||^{2})dt\right) ^{1/4}\leq C\varepsilon
\label{dd1}
\end{eqnarray}%
as a direct consequence of (\ref{in}), (\ref{lim1}) and (\ref{yy}).

Since the function$~$ $\left( 1-\mathbb{\chi}\right) \mathbb{~}\nabla
_{y}L\left( t,U,Y\right) ~$ is bounded by $\lambda ,$ therefore, using again
(\ref{in}), (\ref{lim1}) and (\ref{yy}), we deduce%
\begin{equation}
\mathbb{E}\int_{0}^{T}\left( \left( 1-\mathbb{\chi}\right) \mathbb{~}\nabla
_{y}L\left( t,U,Y\right) ,Z_{n}-Z\right) dt\rightarrow 0\qquad \text{as }%
n\rightarrow \infty .  \label{dd2}
\end{equation}%
Combining (\ref{dd1}) and (\ref{dd2}), we derive (\ref{zzz})$_{1}$. \ \ The
justification of other two limits in (\ref{zzz}) can be done by a similar
way.

The function $Z_{n}$ \ and the pair $\left( p_{n},q_{n}\right) $ fulfill (%
\ref{dual}), \ therefore taking the limit transition in this equality we
show that $Z$ and the pair $(p,q)$ satisfy equality (\ref{ww}). $\hfill
\blacksquare $\newline

\bigskip

\bigskip

\section{Existence and optimality condition}

\bigskip

\label{81}

Let us state the main result of this article. Which essentially says that
the optimal control $U\in \mathcal{U}_{ad}^{b}$ exists and may be determined
through the optimality condition.

Let us recall the assumptions made along the article. We assume $U\in
\mathcal{U}_{ad}^{b}$ satisfying (\ref{bound0})-(\ref{bound}) and the
initial condition
\begin{equation*}
Y_{0}\in L^{p}(\Omega ,V)\cap L^{\infty }(\Omega ,\widetilde{W})\bigskip ,
\end{equation*}
moreover, we admit that (\ref{g}) and one of conditions (\ref{g1}), (\ref{g2}%
) hold. \ The functional $J$ \ satisfies conditions H1)-H3) \ (see (\ref%
{lip1})-(\ref{lip2})) and additionally it is lower semi-continuous with
respect to the weak topology of $\left( L^{2}(\Omega ,L^{\infty }(0,T;%
\widetilde{W}))\right) ^{2}.$ \

\begin{theorem}
\label{main_1} Under above mentioned conditions the control problem $(%
\mathcal{P})$ admits, at least, one optimal solution
\begin{equation*}
(\widehat{U},\widehat{Y})\in \mathcal{U}_{ad}^{b}\times \left( L^{2}(\Omega
,L^{\infty }(0,T;\widetilde{W}))\cap L^{p}(\Omega ,L^{\infty
}(0,T;V))\right) ,
\end{equation*}%
where $\widehat{Y}$ is the unique solution of \ system (\ref{equation_etat})
with $U=\widehat{U}.$

In addition, there exists a unique solution
\begin{equation*}
\left( \widehat{p},\widehat{q}\right) \in L^{2}(\Omega ,L^{\infty
}(0,T;W))\times L^{2}(\Omega ,L^{2}(0,T;W))
\end{equation*}%
of the backward stochastic differential equation
\begin{equation}
\left\{
\begin{array}{ll}
\displaystyle d\sigma (\widehat{p})=\left( -\nu \Delta \widehat{p}-\mathrm{%
curl}\,\sigma (\widehat{Y})\times \widehat{p}+\mathrm{curl}(\sigma (\widehat{%
Y}\times \widehat{p}))+\nabla {\pi }-f\right) \,dt\vspace{2mm} &  \\
\qquad \qquad \qquad \qquad -(\nabla _{y}L(t,\widehat{U},\widehat{Y})+\nabla
_{y}G(t,\widehat{Y})^{T}\widehat{q})\,dt+\sigma \left( \widehat{q}\right)
\,dW_{t}\mathbf{,} &  \\
\mathrm{div}\,\widehat{p}=0 & \mbox{in}\ (0,T)\times \mathcal{O},\vspace{2mm}
\\
\widehat{p}\cdot \mathrm{n}=0,\qquad (\mathrm{n}\cdot D\widehat{p})\cdot
\mathrm{\tau }=0 & \mbox{on}\ (0,T)\times \Gamma ,\vspace{2mm} \\
\widehat{p}(T)=h(\widehat{Y}(T)) & \mbox{in}\ \mathcal{O},%
\end{array}%
\right.  \label{adj_opt_eq_alpha}
\end{equation}%
such that if$~$ $\widehat{Z}$ \ is the solution of system (\ref{5.3}) for $%
~Y=\widehat{Y}~$ and $~\Psi $ replaced by $\Psi -\widehat{U},~$ then the
duality property
\begin{equation}
\mathbb{E}\int_{0}^{T}(\nabla _{y}L(t,\widehat{U},\widehat{Y}),\widehat{Z}%
)\,dt+\mathbb{E(}(\widehat{Z}(T)),\nabla h(\widehat{Y}(T)))_{V}=\mathbb{E}%
\int_{0}^{T}(\Psi -\widehat{U},\widehat{p})\,dt  \label{dual1}
\end{equation}%
and the following optimality condition%
\begin{equation}
\mathbb{E}\int_{0}^{T}\left\{ (\nabla _{u}L(t,\widehat{U},\widehat{Y}),\Psi -%
\widehat{U})+\int_{0}^{T}(\Psi -\widehat{U},\widehat{p})\,\right\} \,dt\geq 0
\label{opt1}
\end{equation}%
hold$.$
\end{theorem}

\bigskip

\textbf{Proof. \ }\textit{1st step: existence of an optimal control.} Let us
consider a minimizing sequence $\left\{ (U_{n},Y_{n})\right\} $ of the
functional $J$. \ Since the sequence $\left\{ U_{n}\right\} $ \ belongs to $%
\mathcal{U}_{ad}^{b}$ which is a closed convex and bounded subset of $%
L^{p}(\Omega ,L^{2}(0,T;H^{1}(\mathcal{O})))$, $4\leq p<\infty ,$ then there
exists a subsequence, still denoted by $\left\{ U_{n}\right\} _{n=1}^{\infty
}$, that converges weakly to $\widehat{U}\in L^{p}(\Omega ,L^{2}(0,T;H^{1}(%
\mathcal{O}))).$ As a consequence of Theorem \ref{the_1} the sequence $%
\left\{ Y_{n}\right\} $ \ is uniformly bounded in
\begin{equation*}
L^{2}(\Omega ,L^{\infty }(0,T;\widetilde{W}))\cap L^{p}(\Omega ,L^{\infty
}(0,T;V)).
\end{equation*}%
Therefore, there exists $\widehat{Y},$ \ such that
\begin{equation*}
Y_{n}\rightarrow \widehat{Y}\qquad \mbox{ weakly in  }L^{2}(\Omega
,L^{\infty }(0,T;\widetilde{W}))\cap L^{p}(\Omega ,L^{\infty }(0,T;V)).
\end{equation*}%
Using the same reasoning as in the proof of Theorem \ref{the_1} (see the
details in the article \cite{CC16}), we may verify that $\widehat{Y}$ is
solution of the stochastic differential system (\ref{equation_etat}). On the
other hand the functional $J$ \ is lower semi-continuous with respect to the
weak topology of $\ \left( L^{2}(\Omega ,L^{\infty }(0,T;\widetilde{W}%
))\right) ^{2}$
\begin{equation*}
J(\widehat{U},\widehat{Y})\leq \underline{\lim }_{n\rightarrow \infty
}J(U_{n},Y_{n}),
\end{equation*}%
which implies that $(\widehat{U},\widehat{Y})$ is an optimal pair.

\textit{2nd step: Proof of the necessary optimality condition (\ref{opt1}).}
Let $(\widehat{U},\widehat{Y})$ be the optimal control pair. Let us consider
arbitrary $\Psi \in \mathcal{U}_{ad}^{b}$. Then from Propositions \ref{Gat}, %
\ref{Gat1} the G\^{a}teaux derivative of the cost functional $J$ at the
point $\widehat{U}$ in the direction $\Psi -\widehat{U}$ is given by
\begin{equation*}
\begin{array}{ll}
\frac{dJ(\widehat{U}+\rho (\Psi -\widehat{U}))}{d\rho }\bigg\vert_{\rho =0}=
& \mathbb{E}\int_{0}^{T}\left\{ (\nabla _{u}L(t,\widehat{U},\widehat{Y}%
),\Psi -\widehat{U})+(\nabla _{y}L(t,\widehat{U},\widehat{Y}),\widehat{Z}%
)\right\} \,dt \\
& +\mathbb{E(}\nabla h(\widehat{Y}),\widehat{Z})_{V},%
\end{array}%
\end{equation*}%
where $\widehat{Z}$ \ is the G\^{a}teaux derivative of the control to state
map at point $\widehat{U}$ \ in the direction $\Psi -\widehat{U}$, being the
unique solution of the stochastic linearized system (\ref{5.3}) with $\Psi $
replaced by $\Psi -\widehat{U}$.

The set $\mathcal{U}_{ad}^{b}$ is convex, hence in the above calculation of $%
\frac{dJ(\widehat{U}+\rho (\Psi -\widehat{U}))}{d\rho }$ we could take $\rho
\in (0,1)$, that implies%
\begin{equation*}
\mathbb{E}\int_{0}^{T}\left\{ (\nabla _{u}L(t,\widehat{U},\widehat{Y}),\Psi -%
\widehat{U})+(\nabla _{y}L(t,\widehat{U},\widehat{Y}),\widehat{Z})\right\}
\,dt+\mathbb{E(}\nabla h(\widehat{Y}(T)),\widehat{Z}(T))_{V}\geq 0.
\end{equation*}%
Since here the duality property (\ref{ww}) can be written as (\ref{dual1}),
we easily verify that the solution $\widehat{p}$ of (\ref{adj_opt_eq_alpha}%
)\ fulfills the optimality condition (\ref{opt1}).$\hfill \hfill
\blacksquare $

\smallskip

\bigskip

\begin{remark}
Let us notice that a possible interpretation of the hypothesis (\ref{g1}) or
(\ref{g2}) concerning the result obtained in Theorem \ref{main_1} is the
following:

a) In the case of "small" influence of the noise, i.e. $L$ is taken small,
then condition (\ref{g1}) holds, implying the controllability of the motion
of the fluid for any viscosity values $\nu ;$

b) In the case of "small" time interval $[0,T],$ the condition (\ref{g1})
still holds, and we are able to control the fluid motion for any values of
the viscosity and for any intensity of the noise;

c) For non axisymmetric domains, the "very" viscous fluids satisfy (\ref{g2}%
) and its motion can be controlled for any intensity of the noise and any
interval of time.
\end{remark}

\bigskip

\bigskip

\textbf{Acknowledgment}

The work of N.V. Chemetov was supported by the Funda\c{c}\~{a}o para a Ci%
\^{e}ncia e a Tecnologia (Portuguese Foundation for Science and Technology),
project project UID/MAT/04561/2013.

The work of F. Cipriano was partially supported by the Funda\c{c}\~{a}o para
a Ci\^{e}ncia e a Tecnologia (Portuguese Foundation for Science and
Technology) through the project UID/MAT/00297/2013 (Centro de Matem\'{a}tica
e Aplica\c{c}\~{o}es).

\end{document}